\documentclass[12pt,leqno,oneside]{amsart}
\usepackage{mathrsfs,dsfont}
\usepackage{amsmath,amstext,amsthm,amssymb,bbm}
\usepackage{charter}
\usepackage{typearea}
\usepackage{pdfsync}
\usepackage[backref=page]{hyperref}

\usepackage{listings}
\usepackage{enumitem}
\usepackage{comment}

\usepackage{changepage}

\usepackage{color}

\usepackage{xfrac,float}

\usepackage[a4paper,left=2.3cm, right=2.3cm,
top=2cm,bottom=2cm,
footskip=1cm]{geometry}

\newtheorem{Theorem}{Theorem}[section]
\newtheorem{Proposition}[Theorem]{Proposition}
\newtheorem{Lemma}[Theorem]{Lemma}
\newtheorem{Corollary}[Theorem]{Corollary}
\newtheorem{Question}[Theorem]{Question}
\theoremstyle{definition}
\newtheorem{Definition}[Theorem]{Definition}
\newtheorem{Example}[Theorem]{Example}
\newtheorem{Remark}[Theorem]{Remark}

\newtheorem{Geometric Setting}[Theorem]
{Condition}

\numberwithin{equation}{section}

\newcommand{\bb}{\boldsymbol}

\newcommand{\dbar}{{\overline\partial}}

\DeclareMathOperator{\supp}{supp}
\DeclareMathOperator{\esupp}{ess.supp}

\DeclareMathOperator{\Vol}{Vol}
\DeclareMathOperator{\Div}{Div}
\DeclareMathOperator{\Bl}{Bl}

\DeclareMathOperator{\Ran}{Im}
\newcommand{\cali}[1]{\mathscr{#1}}

\newcommand{\cC}{\cali{C}}
\newcommand{\cL}{\cali{L}}\newcommand{\cH}{\cali{H}}
\newcommand{\field}[1]{\mathbb{#1}}

\newcommand{\R}{\field{R}}
\newcommand{\C}{\field{C}}
\newcommand{\N}{\field{N}}

\newcommand{\D}{\field{D}}
\newcommand{\F}{\field{F}}

\newcommand{\cLL}{\mathcal{L}^{2}}
\newcommand{\E}{\mathbb{E}}

\renewcommand{\P}{\mathbb{P}}

\begin{document}

\title{Toeplitz operators and zeros of square-integrable random holomorphic sections}

\author{Alexander Drewitz, Bingxiao Liu and George Marinescu} 
\address{Universit{\"a}t zu K{\"o}ln,  Department Mathematik/Informatik,
    Weyertal 86-90,   50931 K{\"o}ln, Germany}
    \email{adrewitz@uni-koeln.de}
\email{bingxiao.liu@uni-koeln.de}
\email{gmarines@math.uni-koeln.de}
\thanks{The authors are partially supported by the 
DFG Priority Program 2265 `Random Geometric Systems' (Project-ID 422743078).}
\thanks{G.\ M.\ is partially supported by DFG funded projects SFB/TRR 191
(Project-ID 281071066-TRR 191), 
and the ANR-DFG project QuaSiDy (Project-ID 490843120).
}

\date{\today}
\begin{abstract}
We use the theory of abstract Wiener spaces to construct a probabilistic model
for Berezin--Toeplitz quantization on a complete Hermitian complex manifold
endowed with a positive line bundle.
We associate to a function with compact support (a classical observable)
a family of square-integrable
Gaussian holomorphic sections. Our focus then is on the asymptotic distributions
of their zeros in the semi-classical limit, in particular, we prove
equidistribution results, large deviation estimates,
and central limit theorems of the random zeros on the support of the given function.
One of the key ingredients of our approach is the local asymptotic expansions of
Berezin--Toeplitz kernels with non-smooth symbols.

\end{abstract}

\maketitle

\tableofcontents

\section{Introduction}\label{introduction}

In this paper we prove several probabilistic
results on the action of a 
classical observable on random quantum states,
via Berezin--Toeplitz quantization. More precisely, we investigate 
the distribution of zeros, laws of large numbers,
large deviation estimates and central limit theorems. 
Given a symplectic manifold $X$ the Berezin--Toeplitz
quantization is a family of Hilbert spaces 
$\cH_\hbar$, where $\hbar$ is the Planck constant considered
as a small parameter, together with linear maps $T_\hbar$ from
$\cC^\infty(X)$ to
the space of bounded linear operators $\cL(\cH_\hbar)$.
From a physics point of view, the manifold $X$ can be seen as the phase space
of a physical system, a function $f\in\cC^\infty(X)$
as a classical observable, and the operator $T_{f,\hbar}=T_\hbar(f)$
the corresponding quantum observable.
A fundamental principle states that quantum mechanics
contains the classical one as the limiting case $\hbar\to0$.

Here our quantum spaces will be
$\cH_{1/p}:=H^0_{(2)}(X,L^p\otimes E)$, $p\in\N$ (thus $\hbar=1/p$), consisting
of square-integrable
holomorphic sections of the $p$-tensor powers of a positive holomorphic
line bundle $(L,h_L)\to X$ twisted with an auxiliary Hermitian holomorphic 
line bundle $(E,h_E)$.
The operators $T_{f,p}$ will be Toeplitz operators with symbol 
$f\in\cC^\infty(X)$, or more generally, also non-continuous symbols.
The Berezin--Toeplitz quantization and the underlying techniques
have many applications, ranging from symplectic topology \cite{Pol18},
asymptotics of the analytic torsion forms \cite{MR3615411},
topological quantum field theory \cite{MR2195137}, 
entanglement entropy \cite{MR4093864}, to
non-commutative geometry \cite{LaMcSa08}.

In this paper we focus on probabilistic aspects
of the Berezin--Toeplitz quantization.
For this purpose, we recall that in \cite{DrLM:2023aa} we introduced a general method 
to randomize the quantum states in $H^0_{(2)}(X,L^p\otimes E)$ 
by using Toeplitz operators and considering
random combinations of pure states.
For appropriate symbols $f$ such that $T_{f,p}$ is Hilbert-Schmidt
and injective, we consider random sections of the form
\begin{equation}
\bb{S}_{f,p}=\sum_{j=1}^{d_p}\eta^p_{j}\lambda^p_{j}S^p_{j}
\label{eq:1.1intropart}
\end{equation}
in
$H^0_{(2)}(X,L^p\otimes E)$, 
where $\{\eta^p_{j}\}_{j=1}^{d_p}$ denotes a sequence of independent and 
identically distributed (i.i.d.)\ standard
complex Gaussian variables, $(\lambda^p_{j})_{j=1}^{d_p}$
is the point spectrum of $T_{f,p}$ on $H^{0}_{(2)}(X,L^p\otimes E)$, 
and $\{S^p_{j}\}_{j=1}^{d_p}$ is an orthonormal basis 
of $H^0_{(2)}(X,L^p\otimes E)$ of $H^0_{(2)}(X,L^p\otimes E)$ such that
$T_{f,p}S^p_{j}=\lambda^p_{j}S_{j}$.
The rigorous definition of the probability
distribution on $H^{0}_{(2)}(X,L^p\otimes E)$ in \cite{DrLM:2023aa} proceeds by using
the machinery of constructing Gaussian measures on an abstract
Wiener space introduced  by Gross \cite{Gross67}. 
Our results concern the zero divisors of the random sections
$\bb{S}_{f,p}$. First, we describe how
the classical observable $f$ and its quantum counterpart $T_{f,p}$
influence statistical properties of $\bb{S}_{f,p}$ such as its expectation, variance, etc.
Subsequently, we consider the semiclassical limit 
$p\to\infty$ and obtain, as in many inverse problems, several
features of the geometric input and of the observable $f$. 
Note that on small length scales of order of the Planck scale $1/\sqrt{p}$, 
one loses track of the special features of the geometrical setting and
obtains a universal limiting behavior of random zeros 
\cite{AF:2022,MR1736335,MR1794066,DrLM:2023aa,MR1383056,MR3439098}. 
In our setting, we will observe a universal limiting behavior which (to leading order) is independent of the specific choice of our function $f$, as long as we restrict ourselves to the support of $f$; cf.\ Corollary
\ref{C:univ} and in Theorem \ref{thm:6.4number} for further details.

The distribution results from \cite{DrLM:2023aa} apply to functions $f$
which vanish up to order two. Their derivations are based on the
asymptotics of the kernels of Toeplitz operators
on the support of their symbol $f$ and the calculation
of the first several terms in the asymptotics.
In this article, we take a different approach and
prove a large deviation estimate from which the
distribution of the zeros follows, independent of the vanishing
order of the symbol. Moreover, the results hold
for a large class of non-smooth symbols $f$. For $X$ compact we will provide semiclassical
estimates for the lowest eigenvalues of $T_{f,p}$,
which are intimately linked to
the distribution of zeros.

\subsection{Geometric setting and square-integrable random holomorphic sections}
We now describe the geometric setting of our paper.
Let $(X,J,\Theta)$ be a connected complex Hermitian (paracompact) manifold of 
complex dimension $n$, where $J$ denotes the canonical complex 
structure of $X$, and $\Theta$ denotes a $J$-compatible Hermitian 
form. Then we have an induced
Riemannian metric $g^{TX}(\cdot,\cdot)=\Theta(\cdot,J\cdot)$ on $X$,
and the corresponding Riemannian volume form $\mathrm{dV}:=\Theta^n/n!$. 
With $\mathbb{F}=\C \text{ or }\R$, let $\mathcal{L}^{\infty}(X,\mathbb{F})$
denote the space of (essentially) bounded measurable functions on $X$,
and let $\mathcal{L}^{\infty}_{\mathrm{c}}(X,\F)$ denote the subspace of
$\mathcal{L}^{\infty}(X,\F)$ consisting of functions with compact (essential) support.

Let $L$, $E$ be two holomorphic line bundles on $X$ equipped with 
smooth Hermitian metrics $h_L, h_E$, and let $R^L$, $R^E$
denote their corresponding Chern curvatures. 
The first Chern form of $(L,h_L)$ is defined as
\begin{equation}
c_1(L,h_L):=\frac{\sqrt{-1}}{2\pi} R^L,
\end{equation}
which is a real $(1,1)$-form on $X$ representing
the first Chern class in both de Rham and Dolbeault cohomologies.

In the semiclassical setting, we assume $(L,h_L)$ to be positive and
consider the high tensor powers of $(L,h_L)$, that is, for $p\in\N_{\geq 1}$, 
the Hermitian line bundle
$(L^p\otimes E,h_p):=(L^{\otimes p}\otimes E, h_L^{\otimes p}\otimes h_E)$
on $X$. The space of square-integrable sections of $L^p\otimes E$
with respect to $h_p$ and $\mathrm{dV}$ is denoted by 
$\mathcal{L}^2(X,L^p\otimes E)$, endowed with the $\cLL$-norm
$\|\cdot\|^2_{\cLL}$.
The quantum space in this paper will be the space
of square-integrable holomorphic sections of $L^p\otimes E$,
\begin{equation}
H^0_{(2)}(X,L^p\otimes E):=\left\{s_p\in H^0(X,L^p\otimes E)\;: \; 
\|s_p\|^2_{\cLL}=\int_X |s_p(z)|^2_{h_p} \mathrm{dV}(z)<\infty\right\}.
\end{equation}
Then $H^0_{(2)}(X,L^p\otimes E)$ together with the $\cLL$-inner product
becomes a (separable) complex Hilbert space. 
We set
\begin{equation}
d_p:=\dim_\C H^0_{(2)}(X,L^p\otimes E) \in \N\cup\{\infty\}
\end{equation}
and denote by
\begin{equation}\label{eq:Bproj}
P_p:\mathcal{L}^2(X,L^p\otimes E)\to
H^0_{(2)}(X,L^p\otimes E)
\end{equation}
the orthogonal (Bergman) projection.
One fundamental method to study the sequence of Hilbert spaces 
$H^0_{(2)}(X,L^p\otimes E)$ is through their associated reproducing kernels
$P_p(x,y)$, called Bergman kernels. 
The asymptotic expansion of Bergman kernels as $p\to\infty$
was obtained in \cite[Section 6]{MM07} under the following
very general geometric conditions, which we will also
assume in this paper.

\begin{Geometric Setting}\label{condt:main}
The Riemannian metric $g^{TX}$ is complete on $X$, 
and there exist $\varepsilon_{0}>0$, 
$C_{0}>0$ such that
\begin{equation}
	\sqrt{-1}R^{L}\geq \varepsilon_{0} \Theta,\quad 
	\sqrt{-1}(R^{\mathrm{det}}+R^{E})>-C_{0}\Theta,\quad \text{and} \quad |\partial 
	\Theta|_{g^{TX}}< C_{0}.
	\label{eq:assumptionsintro}
\end{equation}
\end{Geometric Setting}

The first inequality in \eqref{eq:assumptionsintro} says that
$(L,h_L)$ is uniformly positive on $X$, and Condition \ref{condt:main}
also implies that there exist $C>0$, $p_0\in\N$, such that $d_p\geq Cp^n$ 
for $p\geq p_0$.
Condition \ref{condt:main} is a necessary assumption in Sections 
\ref{ss2nov} and  \ref{ss4Jan} -- \ref{section:statistics} 
which deal with semi-classical limits.

In \cite[Sections 2 and 3]{DrLM:2023aa}, we constructed a Gaussian random holomorphic section $\bb{S}_p$ in $H^0_{(2)}(X,L^p\otimes E)$ via the formula
\begin{equation}
\bb{S}_p:=\sum_{j=1}^{d_p}\eta^p_j S^p_j,
\label{eq:1.5intropart}
\end{equation}
where $\{S^p_j\}_j$ is an orthonormal basis of $H^0_{(2)}(X,L^p\otimes E)$, and $\{\eta^p_j\}_{j}$ is a family of i.i.d.\ standard complex Gaussian random variables. Moreover, we have uniqueness in the sense that the distribution of $\bb{S}_p$ does not depend on the choice of the orthonormal basis $\{S^p_j\}_j$. In \cite[Section 3]{DrLM:2023aa}, we have studied equidistribution results, large deviation estimates and hole probabilities for the zeros of $\bb{S}_p$ in the semi-classical limit. We also refer to \cite{DrLM:2023aa} for further discussions on  random zeros and random holomorphic sections in complex geometry.

However, when $d_p=\infty$, the random section $\bb{S}_p$ turns out to be almost surely not  square-integrable on $X$. 
Then the Berezin--Toeplitz quantization came into our construction in \cite[Section 4]{DrLM:2023aa} to define a Gaussian random $\cLL$-holomorphic sections.
Given $f\in\mathcal{L}^\infty_{\mathrm{c}}(X,\R)$ a real bounded measurable function with compact (essential) support, the Toeplitz operators associated to $f$ are defined for $p\in\N$,
 \begin{equation}
T_{f,p}:=P_p M_f \in\mathrm{End}(H^0_{(2)}(X,L^p\otimes E)),
\end{equation}
where $M_f$ denotes the pointwise multiplication by $f$. 
Moreover, $T_{f,p}$ is a self-adjoint Hilbert-Schmidt operator. 
We also set $T^2_{f,p}=T_{f,p}\circ T_{f,p}$ and denote by $T^2_{f,p}(x,y)$
the smooth integral kernel of $T^2_{f,p}$ with respect to the metric
$h_p$ and the volume form $\mathrm{dV}$.

Let $\Ran T_{f,p}$ denote the range of $T_{f,p}$ in $H^0_{(2)}(X,L^p\otimes E)$,
and let $H^0_{(2)}(X,L^p\otimes E, f):=\overline{\Ran T_{f,p}}
=(\ker T_{f,p})^{\perp}$
be the closure of the range of $T_{f,p}$, which itself is also a Hilbert space. 
As explained in \cite[Section 4, in particular Remark 4.15]{DrLM:2023aa}, 
regarding $T_{f,p}$ as an injective Hilbert-Schmidt operator on 
$H^0_{(2)}(X,L^p\otimes E, f)$ and applying the theory of abstract Wiener space, 
we get a unique Gaussian probability measure $\P_{f,p}$ on
$H^{0}_{(2)}(X,L^{p}\otimes E, f)$ which provides a model
for the action of $T_{f,p}$ on $\bb{S}_p$ defined in \eqref{eq:1.5intropart}.
\begin{Definition}\label{def:1.2BZ}
The \textbf{probabilistic Berezin--Toeplitz quantization} associated to the
symbol $f\in \mathcal{L}^\infty_{\mathrm{c}}(X,\R)$ is defined as the
sequence of Gaussian random $\cLL$-holomorphic sections
$\{\bb{S}_{f,p},p\in\N\}$, where each $\bb{S}_{f,p}$
denotes the random variable taking values in 
$H^0_{(2)}(X,L^p\otimes E, f)\subset H^0_{(2)}(X,L^p\otimes E)$ 
with law $\P_{f,p}$. 
An equivalent definition is given by formula \eqref{eq:1.1intropart}.
\end{Definition}
In this paper, we study the asymptotic behavior of the sequence of random
$(1,1)$-currents $\{[\mathrm{Div}(\bb{S}_{f,p})]\}_p$ on $X$
defined by the integration currents on the zero divisors of $\bb{S}_{f,p}$ as
$p\rightarrow +\infty$. In the case of compact K\"{a}hler manifolds 
and a smooth function $f$, such questions were also independently
studied by Ancona--Le Floch \cite{AF:2022} motivated by understanding the
Kodaira embedding twisted by Toeplitz operators.
Moreover, provided a smooth function $f$, in \cite{AF:2022} 
(for compact K\"{a}hler manifolds) and in \cite[Section 5]{DrLM:2023aa}, 
the equidistribution results of $[\mathrm{Div}(\bb{S}_{f,p})]$ as
$p\rightarrow+\infty$ were proved on the subset of the support of $f$, 
where $f$ only vanishes up to order $2$. 
The present article aims to contribute to the above body of work by providing a more profound understanding of the following natural questions:
\begin{itemize}
	\item[(i).] Do the above equidistribution and the large deviation results for $[\mathrm{Div}(\bb{S}_{f,p})]$ on the support of $f$ still hold true when $f$ has higher vanishing orders or a lower regularity?
	\item[(ii).] How are the random zeros distributed asymptotically outside the support of $f$? Can one quantify the difference between random zeros $\frac{1}{p}[\mathrm{Div}(\bb{S}_{f,p})]$ and the expected limit $c_1(L,h_L)$ on a subset where $f$ is supported on the most part of it?
\item[(iii).] Does the central limit type behavior of
$[\mathrm{Div}(\bb{S}_{f,p})]$ hold on the support of $f$? 
More concretely, following the work of Shiffman-Zelditch 
\cite{MR2465693,MR2742043} for $\bb{S}_p$ on 
a compact K\"{a}hler manifold, one can ask for the analogues 
of their results to $[\mathrm{Div}(\bb{S}_{f,p})]$ but probably 
only on the support of $f$.
\item[(iv).]  When $X$ is compact, a problem related to the above is to describe the asymptotic behavior of the
spectra of  $\{T_{f,p}\}_p$.  In particular, when $f$ is non-negative
and not fully supported on $X$, how does the lowest eigenvalues of
$T_{f,p}$ decay to $0$ as $p\rightarrow+\infty$?
\end{itemize}

\subsection{Asymptotic distribution of zeros of random $\cLL$-holomorphic sections}\label{sss1.3a}
At first, we introduce some notions on the regularity and the support of functions on $X$. For $f \in 
\mathcal{L}^{\infty}(X,\mathbb{F})$, let 
$\Vert f\Vert_{\mathcal{L}^{\infty}}=\inf\{C>0:\; |f|\leq C\; \text{a.e.}\}$ 
denote its essential supremum norm with respect 
to the measure $\mathrm{dV}$.
The essential support of $f$ on $X$ (with respect 
to $\mathrm{dV}$), denoted by $\esupp 
f$, is the smallest closed subset of $X$ such that $f$ vanishes almost everywhere
on its complement. When $f$ is also continuous, then 
$\esupp f$ coincides with the support 
$\mathrm{supp} f=\overline{\{x\in X:f(x)\neq0\}}$ of $f$. We will always call $\esupp 
f$ the support of $f$. Note that we say $f$ to be smooth (resp.\ $\mathscr{C}^k$) an open subset $U\subset X$ if 
there is a smooth (resp. $\mathscr{C}^k$) function $\widetilde{f}$ on $U$ such that 
$f|_{U}-\widetilde{f}=0$ almost everywhere in $U$ (with respect to 
$\mathrm{dV}$).

\begin{Definition}
With $\F = \C$ or $\R$,  let $\mathcal{L}^{\infty}_{\mathrm{const}}(X,\F)$ denote the subspace of $\mathcal{L}^{\infty}(X,\F)$ consisting of functions which are constant outside a compact subset, i.e.,
\begin{equation}
\mathcal{L}^{\infty}_{\mathrm{const}}(X,\F)=\big \{f\in \mathcal{L}^{\infty}(X,\F)\,|\, \text{there exists } c_f\in\F \text{ such that } f-c_f\in \mathcal{L}^{\infty}_{\mathrm{c}}(X,\F)\big\}.
\end{equation}
\end{Definition}

\begin{Definition}
For $f \in\mathcal{L}^\infty(X,\R)$, and $U$ an open subset of $X$, we say that $f$ is of class $\mathscr{C}^k$ ($k\in\N\cup\{\infty\}$) almost everywhere {\bf on $U$} if there exists a closed subset $D\subset U$ of Lebesgue measure $0$ such that $f|_{U\setminus D}\in \mathscr{C}^k(U\setminus D,\R)$. We also say that $f$ is of class $\mathscr{C}^k$ ($k\in\N\cup\{\infty\}$) almost everywhere {\bf near $U$} if it is so on an open neighbourhood of $\overline{U}$.
\end{Definition}

\begin{Example}
(i) A smooth function $f\in \mathscr{C}^k(X,\R)$ is always of class
$\mathscr{C}^k$ almost everywhere near any given open subset of $X$.

\noindent
(ii) Let $U$ be an open subset of $X$ such that
$\partial U:=\overline{U}\setminus U$ has Lebesgue measure zero in $X$. 
Then the characteristic function 
	\begin{equation}
	f(x)=\bb{1}_U(x):=\begin{cases}1&\text{ if } x\in U,\\ 0&\text{ if } x\not\in U, \end{cases}
	\end{equation}
	is smooth almost everywhere near $U$.
\end{Example}
We will identify the $(1,1)$-form $R^{L}$ with the Hermitian matrix
\begin{equation} \label{eq:Rdot}
	\dot{R}^{L}\in \mathrm{End}(T^{(1,0)}X)
\end{equation}
such that for $W,Y\in T^{(1,0)}X$,
\begin{equation*}
	R^{L}(W, \bar{Y})=\langle 
	\dot{R}^{L}W,\bar{Y}\rangle_{g^{TX}\otimes \C}.
\end{equation*}
We can now introduce the following quantity in order
to give a bound on the regularity of $f$, 
which is necessary in our methods to obtain 
the asymptotic results for random zeros.

\begin{Definition}\label{def:1.6m}
For any relatively compact subset $U\subset X$, set 
\begin{equation}
\kappa(R^L,U):=\sup\{\max \mathrm{spec}\,(\dot{R}^L_x),\; x\in U\}\geq \varepsilon_0,
\label{eq:1.09intro}
\end{equation}
and  define
\begin{equation}
m(U):=\left\lceil (6n+6)\frac{\kappa(R^L, U)}{\varepsilon_0}\right\rceil\in\N.
\label{eq:1.9intropart}
\end{equation}
In the prequantum case $\Theta=c_1(L,h_L)$, we can take 
$\kappa(R^L,U)=\varepsilon_0=2\pi$ and $m(U)=m(X):=6n+6$,
which is independent of the subset $U$. 
\end{Definition}
We now introduce a quantity that measures
the relative position of an open set with respect to
the the support of a function $f$.

\begin{Definition}\label{def:rfU}
Let $f \in\mathcal{L}^\infty_{\mathrm{c}}(X,\R)$
and let $U$ be an open subset of $X$. Define
\begin{equation}
r(f,U):=\sup\{r> 0\;:\; \text{geodesic ball } 
\mathbb{B}(x,r)\subset U\setminus \esupp\;f 
\},
\label{eq:4.3Jan2024}
\end{equation}
where the geodesic ball is taken with respect to $g^{TX}$.
Since $U$ is assumed to be open, if 
$U\setminus \esupp\;f\neq \varnothing$, 
we have $r(f,U)>0$. When $U\subset \esupp\, f$, 
we say that $f$ has full support on $U$. In this case, 
there is no nontrivial geodesic ball in $U\setminus \esupp\;f$, 
and we set $r(f,U)=0$.
\end{Definition}

The norm $\|\cdot\|_{U,-2}$ for the $(1,1)$-currents is defined
Definition \ref{def:3.5de} setting $\alpha=2$. For any two subsets 
$U$, $U'$ of $X$, the notation $U\Subset U'$ means that 
$\overline{U}$ is compact and contained in $U'$. 
Our first main result is a concentration estimate
which gives an upper bound on the deviation
of $\frac{1}{p}[\mathrm{Div}(\bb{S}_{f,p})]$ from $c_1(L,h_L)$
on an open set in the semi-classical limit in terms of
the relative position with respect to the support of $f$.
\begin{Theorem}\label{thm:4.6Jan2024}
Let $(X,J,\Theta)$ be a connected Hermitian complex manifold 
and let $(L,h_{L})$, $(E,h_{E})$ be two holomorphic line bundles
on $X$ with smooth Hermitian metrics.
Furthermore, assume that Condition \ref{condt:main} holds and fix 
a pair of nonempty open subsets $(U,U')$ of $X$ such that $U\Subset U'$. 
Then there exists a constant $r(U,U')>0$ such that if $f\in \mathcal{L}^{\infty}_{\mathrm{c}}(X,\R)$ is of $\mathscr{C}^{m(U')+1}$ almost everywhere on $U'$ with 
\begin{equation}
\delta_0(f):=\left(\frac{r(f,U')}{r(U,U')}\right)^{1/(2n+2)}<\frac{1}{2}\, , 
\end{equation}
then $\, U\cap \esupp\, f\neq\varnothing\,$, and for any 
	$\delta > \delta_0(f)$, there 
	exists a constant $C=C(U',f,\delta)>0$ and $p_0>0$ such that for all $p\geq p_0$ 
	we have
	\begin{equation}
		\P\Big( \ 
		\Big\|\frac{1}{p}[\Div(\bb{S}_{f,p})] 
		-  c_{1}(L,h_L)\Big\|_{U,-2}>\frac{\delta}{\pi} \  \Big)\leq 
		e^{-C\,p^{n+1}}.
		\label{eq:4.7Jan}
	\end{equation}
	As a consequence, we have
	\begin{equation}
		\P\Big( \ 
		\limsup_{p\rightarrow+\infty}\Big\|\frac{1}{p}[\Div(\bb{S}_{f,p})] 
		-  c_{1}(L,h_L)\Big\|_{U,-2}\leq \frac{\delta_0(f)}{\pi} \  \Big) =1.
		\label{eq:4.73Jan}
	\end{equation}
\end{Theorem}
\begin{Remark}
Observe that when $U\subset\esupp f$, 
we have $r(f,U)=0$, so we need the strict inequality 
$\delta > \delta_0(f)$
in the statement of Theorem \ref{thm:4.6Jan2024}.
\end{Remark}

In general, it is difficult to determine $r(U,U')$
precisely.
By the proof of
Proposition \ref{thm:4.3Jan} the constant $r(U,U')>0$
depends on the geometry of $g^{TX}$, the complex structure of $X$,
and several auxuliary constants.
However, we can still give a rough formula \eqref{eq:4.70Jan} for $r(U,U')$,
which is certainly not sharp.

The question of quantum ergodicity (mass distribution) for a 
sequence of holomorphic sections as tensor power $p$ tending
to infinity is a parallel problem to the asymptotic distributions
of their zeros as integration currents, whose central objects are
the following measures on $X$ defined by the holomorphic sections.

\begin{Definition}\label{def:4.10mass}
The mass distribution of a section $s_p\in H^0(X,L^p\otimes E)$ 
is defined as the measure
\begin{equation}
M_p(s_p):=\frac{1}{p^n}|s_p(z)|^2_{h_p}\mathrm{dV}
\end{equation}
on $X$, where $\mathrm{dV}=\Theta^n/n!$. 
If $s_p$ is square-integrable, then $M(s_p)$ is a finite measure.
\end{Definition}

 For Gaussian (or sub-Gaussian) holomorphic sections on compact K\"{a}hler manifolds or certain random polynomials on $\C^n$, such problems were investigated in \cite{NoVo:98, ShZ99, MR3742811, MR4159385}. In Section \ref{ss:mass}, we also consider the mass distribution of our random $\cLL$-holomorphic sections $\bb{S}_{f,p}$. In particular, we get a law of large numbers for $\int_X g(z)M_p(\bb{S}_{f,p})(z)$. Now we can state the result, as an analog of \cite[Theorem 1.4]{MR3742811}, whose proof will be given in Section \ref{ss:mass}.

\begin{Proposition}\label{prop:4.19p2}
Let $(X,J,\Theta)$ be a connected Hermitian complex manifold
and let $(L,h_{L})$, $(E,h_{E})$ be two holomorphic line bundles
on $X$ with smooth Hermitian metrics.
We assume that the Condition \ref{condt:main} holds. 
Let $U$ be a relative compact open subset of $X$, 
and fix a nontrivial 
$f\in\mathscr{C}^{m(U)+1}_{\mathrm{c}}(U,\R_{\geq 0})$. 
Then for any $g\in\mathscr{C}^0_{\mathrm{c}}(U)$, we have $\mathbb{P}$-a.s.\ that
\begin{equation}
\lim_{N\rightarrow +\infty}
\frac{1}{N}\sum_{1\leq p\leq N}\int_U 
g(z)M_{p}(\bb{S}_{f,p})(z)=
\int_U g(z)f(z)^2\mathrm{dV}^L(z),
\end{equation}
where the volume form $\mathrm{dV}^L:=c_1(L,h_L)^n/n!$ 
in the limit is defined independently from $\Theta$.
\end{Proposition}

\subsection{Large deviation and equidistribution on the 
support of {$f$}}\label{s1.4a}

As a special case of Theorem \ref{thm:4.6Jan2024}, 
we can give a type of large deviation estimates for the zeros of 
$\bb{S}_{f,p}$ on the support of $f$.
Such kind of estimates are also referred as concentration 
of measure. As a consequence, we obtain the equidistribution
of random zeros on the support of $f$, 
see also \cite[Theorem 1.1]{ShZ99} and 
\cite[Sections 3.6 and 5.2]{DrLM:2023aa}).

The following theorem generalizes the previous results
in \cite[Theorems 1.6 and 1.7]{DrLM:2023aa} by removing 
the conditions on the smoothness of $f$ and the vanishing
orders of $f$ on the given domain. Note that the quantity $m(U)$ 
is defined by \eqref{eq:1.9intropart}.

\begin{Theorem}\label{thm:6.2}
Let $(X,J,\Theta)$ be a connected complex Hermitian manifold and let $(L,h_{L})$, $(E,h_{E})$ be two holomorphic line bundles on $X$ with smooth Hermitian metrics.
We assume that the condition \ref{condt:main} holds. Fix $0\neq f\in \mathcal{L}^{\infty}_{\mathrm{c}}(X,\R)$. Let $U$ be an open subset of $X$ such that $U\subset \mathrm{ess.supp\,} f$ and $f$ is of class $\mathscr{C}^{m(U)+1}$ almost everywhere on $U$. Then for any 
	$\delta >0$, and $V\Subset U$, there 
	exists a constant $C=C(f,\delta, V)>0$ such that for all sufficiently large $p\in\N$,
	we have
	\begin{equation}
\P\Big( \ \Big\|\frac{1}{p}[\Div(\bb{S}_{f,p})] 
-  c_{1}(L,h_L)\Big\|_{V,-2}>\delta \  \Big)\leq 
e^{-C\,p^{n+1}}.
		\label{eq:6.0.7}
	\end{equation}
	As a consequence, we have $\P$-a.s.\ that
		\begin{equation}		
		\lim_{p\rightarrow+\infty}\Big\|\frac{1}{p}[\Div(\bb{S}_{f,p})] 
		-  c_{1}(L,h_L)\Big\|_{U,-2}=0.
		\label{eq:5.2Jan24}
	\end{equation}
\end{Theorem}

\begin{Remark}
In \cite[Section 3]{DrLM:2023aa} we proved 
large deviation estimates and equidistribution results for
$\langle\frac{1}{p}[\Div(\bb{S}_{f,p})]-  c_{1}(L,h_L),\varphi\rangle$
for a fixed $\varphi\in \Omega^{n-1,n-1}_{\mathrm{c}}(X)$,
in the case of Gaussian random holomorphic sections 
$\bb{S}_p$ (defined in \eqref{eq:1.5intropart}).
In \cite[Corollary 3.7]{DrLM:2023aa}, we proved the almost sure 
convergence \eqref{eq:5.2Jan24} under extra finiteness condition. 
But actually, Theorem \ref{thm:6.2} holds for the 
random holomorphic sections $\bb{S}_p$, since
all the asymptotic expansions
for $T^2_{f,p}(x,y)$ necessary to prove Theorem \ref{thm:6.2}
have analogoues for $P_p(x,y)$.
\end{Remark}

We now provide an interesting consequence of Theorem \ref{thm:6.2}.
For any Borel subset $U\subset X$ 
we set
 \begin{equation}
 \mathrm{Vol}^L_{2n}(U):=\int_U \frac{c_1(L,h_L)^n}{n!}>0.
 \label{eq:volomeU}
 \end{equation}
Analogously, if $Y$ is a complex submanifold of $X$ with complex codimension $1$, we define the $(2n-2)$-dimensional volume with 
respect to $c_{1}(L,h_{L})$ of $Y$ as 
 \begin{equation}
 \mathrm{Vol}^L_{2n-2}(Y):=\int_Y \frac{c_1(L,h_L)^{n-1}}{(n-1)!}|_Y.
 \label{eq:volomeofY}
 \end{equation}

For $s_{p}\in H^{0}(X,L^{p}\otimes E)\setminus\{0\}$, the $(2n-2)$-dimensional volume of the divisor $\Div(s_{p})$ (cf.\ \eqref{eq:divisor})
in an open subset $U\subset X$ as follows:
\begin{equation}
	\Vol^{L}_{2n-2}\big(\Div(s_{p})\cap U\big)=\sum_{Y\subset 
	Z(s_{p})}\mathrm{ord}_{Y}(s_{p})\mathrm{Vol}^L_{2n-2}(Y\cap U)\,\cdot
	\label{eq:6.1.9}
\end{equation}
If we use this volume to measure the size of the zeros of $s_{p}$ in $U$, 
then Theorem \ref{thm:6.2} leads to the following result.
\begin{Theorem}\label{thm:6.3}
We assume the same geometric conditions on $X, L, E$ as in Theorem \ref{thm:6.2}. Fix $0\neq f\in \mathcal{L}^{\infty}_{\mathrm{c}}(X,\R)$. Let $U$ be an open subset of $X$ such that $U\subset \mathrm{ess.supp\,} f$ and $f$ is of class $\mathscr{C}^{m(U)+1}$ almost everywhere on $U$.
	If $V$ is a nonempty relatively compact open subset of $U$ such that 
	$\partial V$ has zero measure in $X$, 
	then for any $\delta >0$, there 
	exists a constant $c=c_{f,\delta,V}>0$ such that for all sufficiently large $p$,
	we have
	\begin{equation}
		\P\Big( \ 
		\Big|\frac{1}{p}\Vol^{L}_{2n-2}(\Div(\bb{S}_{f,p})\cap V) 
		- n\Vol^{L}_{2n}(V) \Big|>\delta \ \Big)\leq 
		e^{-c\, p^{n+1}}.
		\label{eq:6.1.10}
	\end{equation}
	In addition, there exists a 
	constant $C_{f,V}>0$ such that for $p\gg 0$,
	\begin{equation}
		\P\big(\Div(\bb{S}_{f,p})\cap V=\varnothing\big)
		\leq e^{-C_{f,V}\,p^{n+1}}\,.
		\label{eq:1.6.14paris}
	\end{equation}
\end{Theorem}
The right-hand side of \eqref{eq:1.6.14paris}
is called the \textit {hole probability} for the random sections $\bb{S}_{f,p}$.

\subsection{Expectation of random zeros and pluripotential theory on $X$}\label{ss:1.4intro}

As a part of the equidistribution results for the zeros of $\bb{S}_{f,p}$, we also need to study the convergence of the expectation of $\frac{1}{p}[\Div(\bb{S}_{f,p})]$ as $(1,1)$-currents on $X$. In Theorem \ref{thm:3.6jan}, we show that
\begin{equation}
\frac{1}{p}\E[[\Div(\bb{S}_{f,p})]]=c_1(L,h_L)+\frac{\sqrt{-1}}{2\pi p} \partial\overline{\partial} \log {T_{f,p}^2(x,x)}+\frac{1}{p}c_1(E,h_E).
\label{eq:1.22intropart}
\end{equation}
Hence the main point is to study the current $\frac{\sqrt{-1}}{2\pi p} \partial\overline{\partial} \log {T_{f,p}^2(x,x)}$. Then by the observation from \eqref{eq:1.22intropart} that $\frac{1}{p}\E[[\Div(\bb{S}_{f,p})]]$ is a positive current on $X$, we can apply the techniques from the pluripotential theory, especially the theory of quasi-plurisubharmonic (quasi-psh) functions, to study the asymptotic properties of $\log {T_{f,p}^2(x,x)}$ as $p\rightarrow+\infty$. We will recall some basics for plurisubharmonic functions in Section \ref{Section:5.2march}.

As a consequence, we have the following theorem in a great generality.
\begin{Theorem}\label{thm:5.5jan24}
Let $(X,J,\Theta)$ be a connected complex Hermitian manifold and let $(L,h_{L})$, $(E,h_{E})$ be two holomorphic line bundles on $X$ with smooth Hermitian metrics.
We assume that the condition \ref{condt:main} holds. For $f\in \mathcal{L}^{\infty}_{\mathrm{c}}(X,\R)$, if there exists a small open ball $B\neq \varnothing$ such that $f|_B\in\mathscr{C}^1(B)$ and $f$ never vanishes on $B$. Let $U$ be a connected open subset of $X$ which is relatively compact and $U\cap B\neq\varnothing$, then there exist constants $C>0,\, C'\geq n$ depending only on $X$, $U$, $L$ and $f$ such that for all $p\geq 1$, any sequence of nonempty open subsets $\{A_p\}_{p\geq 1}$ of $U$, we have 
	\begin{equation}
	\frac{1}{\mathrm{Vol}(A_p)}\int_{A_p} T^2_{f,p}(x,x)\mathrm{dV}(x)\geq \exp\left(-\frac{Cp}{\mathrm{Vol}(A_p)}-C'\log{p}\right),
		\label{eq:4.11Toep}
	\end{equation}
	where $\mathrm{Vol}(A_p):=\int_{A_p}\mathrm{dV}(x)>0$.
	Moreover, with the same constant $C>0$ as above and for all $p\gg 0$, we have
	\begin{equation}
	\left\|\frac{1}{p}\mathbb{E}[[\Div(\bb{S}_{f,p})]]-c_1(L,h_L)\right\|_{U,-2}\leq \frac{C}{\pi}.
		\label{eq:5.13Toep24}
	\end{equation}
	There exists a subsequence $\{p_j\}_{j=1}^\infty\subset \N$ that is increasing to $+\infty$ and a quasi-psh function $\widehat{f}$ on $U$ such that we have the convergence of $(1,1)$-currents of order $2$ on $U$,
	\begin{equation}
	\lim_{j\rightarrow +\infty}\frac{1}{p_j}\mathbb{E}[[\Div(\bb{S}_{f,p_j})]]|_U=c_1(L,h_L)|_U+\frac{\sqrt{-1}}{2\pi} \partial\overline{\partial} \widehat{f}\geq 0.
	\label{eq:5.14limit}
	\end{equation}
	In particular, $\widehat{f}|_{B\cap U}\equiv 0$.
\end{Theorem}

A remark on Theorem \ref{thm:5.5jan24} for a compact Hermitian manifold $X$ is that when $f$ is non-negative (hence $T_{f,p}$ is injective) the constant $C$ in \eqref{eq:4.11Toep} and \eqref{eq:4.11Toep} can be determined by the lowest eigenvalues of $T_{f,p}$. Then a consequence, as we will explain in Sections \ref{ssintro:minimal} and \ref{ss:compactandeigen}, is that the above constant $C$ is related to the size of $\mathrm{ess.supp\,} f$ (or equivalently $X\setminus \mathrm{ess.supp\,} f$).

When we consider the case $U\subset \mathrm{ess.supp\,} f$, we get the convergence of $\frac{1}{p}\mathbb{E}[[\Div(\bb{S}_{f,p})]]$ on $U$ as $p\rightarrow+\infty$. The precise statement is given as follows.
\begin{Theorem}\label{thm:6.2equi}
Let $(X,J,\Theta)$ be a connected complex Hermitian manifold and let $(L,h_{L})$, $(E,h_{E})$ be two holomorphic line bundles on $X$ with smooth Hermitian metrics.
We assume that the condition \ref{condt:main} holds. Fix $0\neq f\in \mathcal{L}^{\infty}_{\mathrm{c}}(X,\R)$. Let $U$ be an open subset of $X$ such that $U\subset \mathrm{ess.supp\,} f$, and we assume that $f$ is of class $\mathscr{C}^1$ almost everywhere near $U$. Then we have, as $p\rightarrow +\infty$,
	\begin{equation}
\left\|\frac{1}{p} \log {T_{f,p}^2(x,x)}\right\|_{\mathcal{L}^1(U,\R)}\rightarrow 0.
		\label{eq:6.0.7equi}
	\end{equation}
Then we have, as $p\rightarrow +\infty$,
		\begin{equation}
		\left\|\frac{1}{p}\mathbb{E}[[\Div(\bb{S}_{f,p})]]- c_1(L,h_L)\right\|_{U,-2}\rightarrow 0.
		\label{eq:5.35feb24}
		\end{equation}
\end{Theorem}
Note that without any extra condition, we can not conclude \eqref{eq:5.35feb24} directly from the almost sure convergence \eqref{eq:5.2Jan24}, so that in the proof of Theorem \ref{thm:6.2equi}, the use of certain compactness result for quasi-psh functions is necessary in our method. The proofs of both theorems above are given in Section \ref{Section:5.2march}.

\subsection{Central limit theorem for random zeros on the support}\label{sss1.2a}

The following theorem extends \cite[Main Theorem]{STr} and \cite[Theorem 
1.2]{MR2742043} to our Toeplitz setting on possibly noncompact Hermitian manifolds. Moreover, as pointed out in \cite[Remark 3.17]{DrLM:2023aa}, such result also holds true for the Gaussian holomorphic sections $\{\bb{S}_p\}_p$ (defined by \eqref{eq:1.5intropart}) on noncompact Hermitian manifolds.

\begin{Theorem}\label{thm:3.5.1ss}
Let $(X,J,\Theta)$ be a connected complex Hermitian manifold and
let $(L,h_{L})$, $(E,h_{E})$ be two holomorphic line bundles on $X$
with smooth Hermitian metrics.
We assume that Condition \ref{condt:main} is satisfied. 
Fix $f\in \mathscr{C}^{\infty}_{\mathrm{c}}(X,\R)$ which is not identically zero,
and let $U$ be an open subset of $X$ such that $\overline{U}\subset \{f\neq 0\}$.
Let $\varphi$ be a real $(n-1,n-1)$-form on $X$ with $\mathscr{C}^3$-coefficients such that $\supp\varphi\subset U$ and $\partial\overline{\partial}\varphi\not\equiv 0$, set
	\begin{equation}
Z_{f,p}(\varphi):=\langle[\mathrm{Div}(\bb{S}_{f,p})],\varphi\rangle,
	\end{equation}
then as $p\rightarrow \infty$, the distribution of the random 
variables
	\begin{equation}
\frac{Z_{f,p}(\varphi)-\E[Z_{f,p}(\varphi)]}{\sqrt{\mathrm{Var}[Z_{f,p}(\varphi)]}}
	\end{equation}
converges weakly to $\mathcal{N}_{\R}(0,1)$.
\end{Theorem}
Note that the smoothness assumption 
$f\in \mathscr{C}^{\infty}_{\mathrm{c}}(X,\R)$
can be relaxed to $f\in \mathscr{C}^{m}_{\mathrm{c}}(X,\R)$
for sufficiently large $m$. 
In the case of compact K\"{a}hler manifolds, 
Shiffman and Zelditch \cite{MR2465693, MR2742043}
also computed explicitly the variance $\mathrm{Var}[Z_{f,p}(\varphi)]$. 
The same method also applies to our case due to our results for the 
normalized Berezin--Toeplitz kernels 
(in particular Theorem \ref{thm:5.1.1} below). 
The proof of Theorem \ref{thm:3.5.1ss} will be given in Section \ref{ss:6.2CLT}.

For a real $(n-1,n-1)$-form  $\varphi$ on $X$ with $\mathscr{C}^3$-coefficients, recall that $L(\varphi)\in \mathscr{C}^1(X,\R)$ is defined by
\begin{equation}
 \sqrt{-1}\partial\overline{\partial} \varphi = L(\varphi) \frac{c_1(L,h_L)^n}{n!}.
\end{equation}
If $U$ is as in Theorem \ref{thm:3.5.1ss} and $\supp\varphi\subset U$, then in Theorem \ref{thm:6.4number}, we prove that
for $p\gg 0$, 
\begin{equation}\label{eq:6.12intro}
\mathrm{Var}[Z_{f,p}(\varphi)]=p^{-n}\left(\frac{\zeta(n+2)}{4\pi^2}\int_U |L(\varphi)(z)|^2 \mathrm{dV}^L(z)+\mathcal{O}(p^{-1/2+\epsilon})\right),
\end{equation}
where
$$\zeta(n+2)=\sum_{k=1}^\infty \frac{1}{k^{n+2}}.$$ 
Note that $\int_U |L(\varphi)(z)|^2 \mathrm{dV}^L(z)$ is a positive quantity independent of the metric $\Theta$ on $X$.

With the same assumptions in Theorem \ref{thm:3.5.1ss}, 
we have $\frac{1}{p}\E[Z_{f,p}(\varphi)]\rightarrow \langle c_1(L,h_L),\varphi\rangle$ as $p\rightarrow+\infty$. Therefore, as a consequence of 
\eqref{eq:6.12intro}, Theorem \ref{thm:3.5.1ss}, and the definition of convergence in distribution as the pointwise convergence of the distribution functions towards the distribution function of the limiting random variable in all points of continuity,
we get the following universality result.
\begin{Corollary}\label{C:univ}
With the same assumptions in Theorem \ref{thm:3.5.1ss}, set
\begin{equation}
\sigma(U,h_L,\varphi):= \frac{\zeta(n+2)}{4\pi^2}
\int_U |L(\varphi)(z)|^2 \mathrm{dV}^L(z),
\end{equation}
then the distribution of the real random variable 
\begin{equation}\label{eq:Div_loc}
p^{n/2}\langle[\mathrm{Div}(\bb{S}_{f,p})]-p c_1(L,h_L),\varphi\rangle
\end{equation}
converges weakly to $\mathcal{N}_{\R}(0,\sigma(U,h_L,\varphi))$
as $p\rightarrow +\infty$. 
\end{Corollary}
This shows the local universality character of the distribution of
zeros of random $\cLL$-holomorphic sections in the sense that the
limiting distribution of the random variable \eqref{eq:Div_loc}
is independent of the function $f\in\mathscr{C}^\infty_{\mathrm{c}}(X,\R)$
under the condition that $\overline{U}\subset\{f\neq 0\}$.

\subsection{Berezin--Toeplitz kernel with non-smooth symbols}\label{sss1.4a}

The Berezin--Toeplitz kernel $T^2_{f,p}(x,y)$ is an essential element in all our proofs to the above results for the Gaussian $\cLL$-holomorphic sections $\bb{S}_{f,p}$ on $X$. For a smooth symbol $f$, the asymptotic expansions of $T^2_{f,p}(x,y)$ (in the non-compact setting) were given by the seminal works of Ma--Marinescu \cite{MM07, MM12} using the techniques of analytic localization. This method remain applicable for the non-smooth symbol $f$ as given by Barron, Ma, Marinescu and Pinsonnault in \cite{MR3822756}. Our results for $T^2_{f,p}(x,y)$ presented in Section \ref{ss2nov} for a non-smooth symbol $f$ can be regarded as the local versions of the results proved in \cite{MR3822756}, and our proofs are still built on the techniques of analytic localization developed in \cite{MM07}.

In this paper, our (local or global) regularity condition on the symbol $f$ will be among $\{\mathcal{L}^\infty,\mathscr{C}^0, \mathscr{C}^1, \mathscr{C}^{m+1}, \mathscr{C}^\infty\}$, depending on the different contexts, as is already alluded to by the theorems in the previous Sections. The smoothness on $f$ generally is not needed (in Sections \ref{ss4Jan} and \ref{section5:onsupport}) but we impose this condition in Section \ref{sss1.2a} and Section \ref{section:statistics} for the sake of simplicity. The $\mathscr{C}^{m+1}$-regularity (with $m=m(U)$) is necessary to have a proper near-diagonal asymptotics for the normalized Berezin--Toeplitz kernels $N_{f,p}$ (cf.\ Definition \ref{def:1.18intro}), which plays the role of correlation function of the holomorphic Gaussian field $\bb{S}_{f,p}$. Therefore, for the large deviation estimates and their consequences, such (local) regularity is assumed. Moreover, in most case, we only need the function $f$ to be $\mathscr{C}^{m+1}$ outside a negligible closed subset, so that we can include the interesting examples of $f$ such as the cut-off functions, indicator functions, etc, in our framework. In the cases where only the uniform bounds or the leading terms of $T^2_{f,p}(x,y)$ are needed, we can assume the function $f$ to be only locally $\mathscr{C}^0, \mathscr{C}^1$ or just $\mathcal{L}^\infty$.

We fix a real function 
$f\in\mathcal{L}^{\infty}_{\mathrm{const}}(X,\R)$. Then $T_{f,p}$ is a self-adjoint bounded operator, and we will mainly be concerned with the nonnegative operator $T_{f,p}^2:=T_{f,p}\circ T_{f,p}$.

\begin{Definition}\label{def:1.18intro}
	For $p\in\N^{*}$, the normalized Berezin--Toeplitz 
	kernel associated to the given $f$ is defined by
	\begin{equation}
		N_{f,p}(x,y):=\frac{|T_{f,p}^{2}(x,y)|_{h_{p,x}\otimes h^{*}_{p,y}}}{\sqrt{T_{f,p}^{2}(x,x)}\cdot 
		\sqrt{T_{f,p}^{2}(y,y)}}, \;\; x,y\in X,
		\label{eq:2.60nov}
	\end{equation}
	whenever $T_{f,p}^{2}(x,x)\neq 0$ and $T_{f,p}^{2}(y,y)\neq 0$.
\end{Definition}

For $x\in X$, let $\mathrm{inj}_{x}^{X}>0$ be the supremum of the radius 
$r>0$ such that the geodesic map $T_{x}X\ni v\mapsto \exp_{x}(v)\in 
X$ is a diffeomorphism when restricting to the open ball 
$B^{T_{x}X}(0,r)\subset T_{x}X$. Fix a relatively compact open subset 
$U\Subset X$, set
\begin{equation}
	\mathrm{inj}^{X}_{U}=\inf_{x\in U} \mathrm{inj}_{x}^{X} >0.
	\label{eq:1.37intro}
\end{equation}

The following theorem extends 
\cite[Theorem 5.1]{DLM:21} to this new setting of Berezin--Toeplitz operators. For $x\in X$, the distance 
function $\Phi_{x}$ will be defined in \eqref{eq:2.16July2}. The constant 
$\varepsilon_{0}$ was introduced in our assumption 
\eqref{eq:assumptionsintro}.

\begin{Theorem}\label{thm:5.1.1}
Let $(X,J,\Theta)$ be a connected complex Hermitian manifold and let $(L,h_{L})$, $(E,h_{E})$ be two holomorphic line bundles on $X$ with smooth Hermitian metrics.
We assume that the condition \ref{condt:main} holds. Fix $f\in\mathcal{L}^{\infty}_{\mathrm{const}}(X,\R)$.
	Let $U$ be a relatively compact open subset of $X$ such that $f$ 
	is smooth on an open neighbourhood of $\overline{U}$, closure of 
	$U$ in $X$, and $f$ does not vanish in $\overline{U}$. Then there 
	exist $0<\delta_{U}<\mathrm{inj}^{X}_{U}/4$ such that the 
	following uniform estimates on the normalized Berezin--Toeplitz kernel hold 
	for $x,y\in U$:
	For $k\geq 1$, there exists a constant $M_{k}>0$ (we may take $M_k=12k$) such that for 
	any fixed 
	$b\geq\sqrt{M_{k}/\varepsilon_{0}}$, we have for 
	all $p \gg 0\,$ with $b\sqrt{\frac{\log 
		p}{p}}\leq \delta_{U}$ that
	\begin{equation}
		N_{f,p}(x,y)=\begin{cases} 
		&\big(1+R_{p,x}(v')\big)\exp\Big(-\dfrac{p}{4}\,\Phi_{x}(0,v')^{2}\Big),  \\
		&\qquad\qquad\text{uniformly for } 
		\mathrm{dist}(x,y)\leq b\sqrt{\frac{\log 
		p}{p}}\,,\;\text{with\;}y=\exp_{x}(v'), v'\in T_{x}X; \\
		&\mathcal{O}(p^{-k}), \;\;   \text{uniformly for } 
		\mathrm{dist}(x,y) \geq 
		b\sqrt{\frac{\log p}{p}},
	\end{cases}
	\label{eq:5.2.4}
\end{equation}
where 
$$\sup_{x\in U,\, v'\in T_x X,\, \|v'\|\leq b\sqrt{\log{p}/p} } |R_{p,x}(v')|\rightarrow 0$$
 as $p\rightarrow +\infty$. More precisely, we have the following estimate, for any fixed $\epsilon\in\;]0,1/2]$, there exists $C=C(f,\epsilon,U)>0$ such that for any $x\in U, v'\in T_x X,\, \|v'\|\leq b\sqrt{\log{p}/p}$,
\begin{equation}
|R_{p,x}(v')|\leq C p^{-1/2+\epsilon}.
\label{eq:2.62feb24}
\end{equation}
\end{Theorem}

The proof of above theorem is given in Section \ref{ss2.4:BT}. Moreover, essentially by the same proof, we get a different version of Theorem \ref{thm:5.1.1} under a lower regularity assumption on $f$ as follows, a brief proof to it is also given in Section \ref{ss2.4:BT}.

Recall that $m(U)=\left\lceil (6n+6)\frac{\kappa(R^L, U)}{\varepsilon_0}\right\rceil$ is given in \eqref{eq:1.9intropart}, this definition actually follows from the choice $M_{n+1}=12(n+1)$ and $b=\sqrt{\frac{M_{n+1}}{\varepsilon_0}}$ in the proof of Theorem \ref{thm:5.1.1}. Then we have the following results.

\begin{Corollary}\label{cor:5.1.1}
	Let $(X,J,\Theta)$ be a connected complex Hermitian manifold and let $(L,h_{L})$, $(E,h_{E})$ be two holomorphic line bundles on $X$ with smooth Hermitian metrics.
We assume that the condition \ref{condt:main} holds. Fix $f\in\mathcal{L}^{\infty}_{\mathrm{const}}(X,\R)$.
	Let $U$ be a relatively compact open subset of $X$ such that $f$ 
	is of $\mathscr{C}^{m+1}$ with $m=m(U)$ on an open neighbourhood of $\overline{U}$ and $f$ does not vanish in $\overline{U}$. Then there 
	exist $0<\delta_{U}<\mathrm{inj}^{X}_{U}/4$ such that the 
	following uniform estimates on the normalized Berezin--Toeplitz kernel hold 
	for $x,y\in U$: set $b=\sqrt{\frac{12n+12}{\varepsilon_0}}$, we have for 
	all $p \gg 0\,$ with $b\sqrt{\frac{\log 
		p}{p}}\leq \delta_{U}$ that
	\begin{equation}
		N_{f,p}(x,y)=\begin{cases} 
		&\big(1+o(1)\big)\exp\Big(-\dfrac{p}{4}\,\Phi_{x}(0,v')^{2}\Big),  \\
		&\qquad\qquad\text{uniformly for } 
		\mathrm{dist}(x,y)\leq b\sqrt{\frac{\log 
		p}{p}}\,,\;\text{with\;}y=\exp_{x}(v'), v'\in T_{x}X; \\
		&\mathcal{O}(p^{-n-1}), \;\;   \text{uniformly for } 
		\mathrm{dist}(x,y) \geq 
		b\sqrt{\frac{\log p}{p}}.
	\end{cases}
	\label{eq:5.2.4Nov}
\end{equation}
\end{Corollary}

An advantage of the above results, Theorem \ref{thm:5.1.1}
and Corollary \ref{cor:5.1.1}, is that we only need the
local regularity of $f$ on the part where we want $N_{f,p}$
to have the asymptotic formula as in \eqref{eq:5.2.4} or 
\eqref{eq:5.2.4Nov}, and away from $U$, there is no 
regularity assumption on $f$.

\subsection{Lowest eigenvalue of Toeplitz operators on compact manifolds}
\label{ssintro:minimal}
Now we focus on a compact Hermitian manifold $(X,\Theta)$. 
In this case, $H^0_{(2)}(X,L^p\otimes E)=H^0(X,L^p\otimes E)$ is finite dimensional, 
and the Gaussian holomorphic section $\bb{S}_p$ (see \eqref{eq:1.5intropart}) 
can be regarded as the identity map on $H^0(X,L^p\otimes E)$ 
after equipping $H^0(X,L^p\otimes E)$ with the standard 
Gaussian probability measure associated to the $\cLL$-inner product. 
For a real bounded (measurable) function $f$ on $X$, the random section 
$\bb{S}_{f,p}$ in Definition \ref{def:1.2BZ} is equivalent to
\begin{equation}
\bb{S}_{f,p}=T_{f,p} \bb{S}_p.
\label{eq:5.35jan24}
\end{equation}

\begin{figure}[H]
\centering
\includegraphics[width=0.8\textwidth]{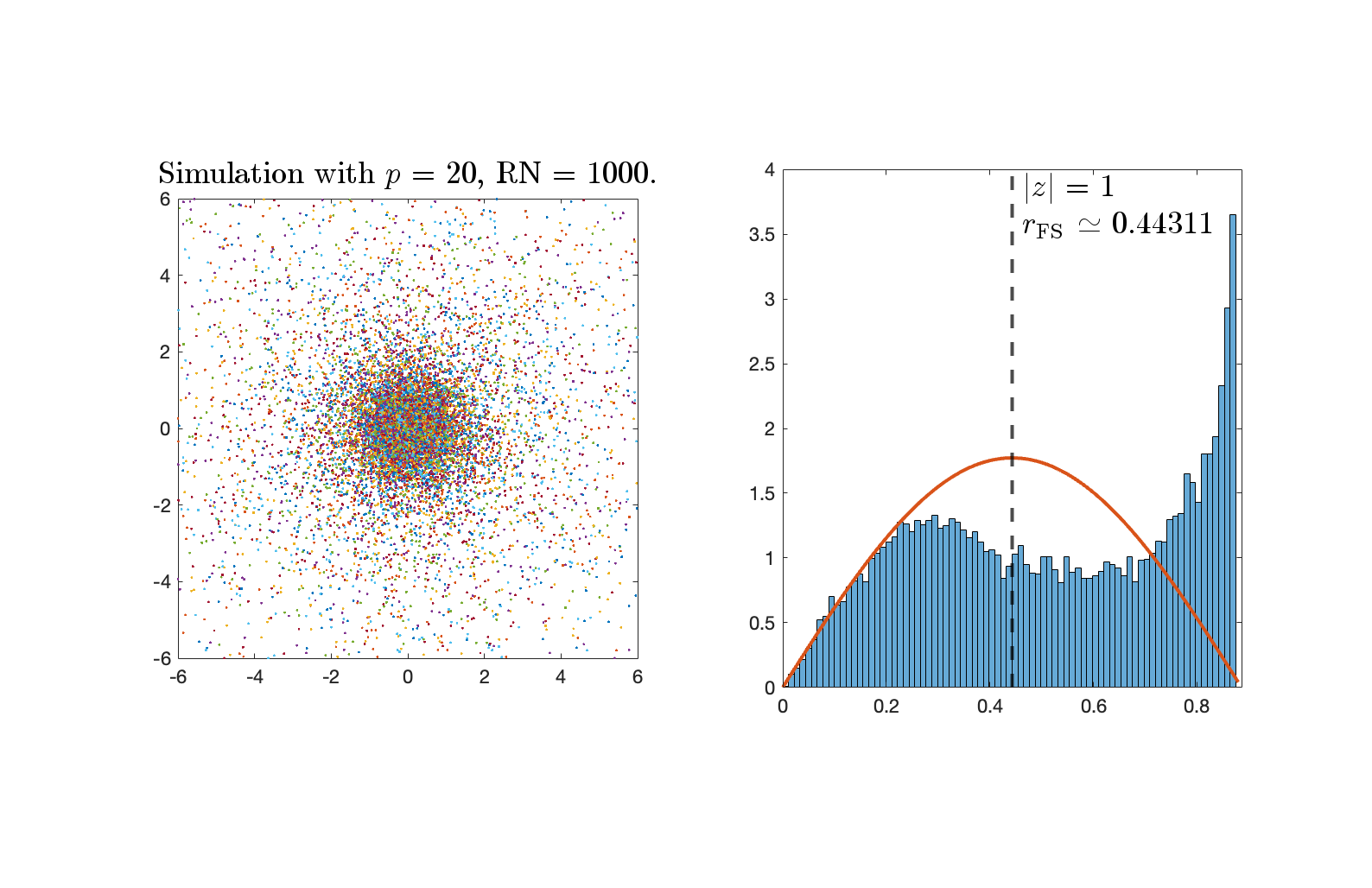}
\caption[Random zeros of $\bb{S}_{\bb{1}_\D , 20}$]
{Comparison of zeros of $\bb{S}_{\bb{1}_\D , 20}$
with $\omega_{\mathrm{FS}}$ on $\mathbb{CP}^1$. 
The density function $\psi(r_{\mathrm{FS}})$
(see \eqref{eq:FScurve}) is plotted as the red curve
in the right-hand side, and the region 
$\{|z|\leq 1\}=\{r_{\mathrm{FS}}\leq 0.44311\ldots\}$ 
gives the support of $\bb{1}_\D$.}
\label{fig:chi1p20}
\end{figure}

Even in this case, the problem about the asymptotic distribution of the random zeros 
$[\mathrm{Div}(\bb{S}_{f,p})]$ outside the support of $f$ 
remains open. In Section \ref{ss:5.4.3new}, 
we present simulations of the zeros of $T_{f,p}\bb{S}_p$
on the Riemann sphere $\mathbb{CP}^1$, where $\bb{S}_p$
is the $\mathrm{SU}(2)$-polynomial. More precisely, if $\D$
is a the unit disc of a standard local chart $U_0\simeq \C$,
which is a geodesic ball in $X=\mathbb{CP}^1$ of 
$g^{TX}_{\mathrm{FS}}$-radius 
$r_\mathrm{FS}=\frac{\sqrt{\pi}}{4}\simeq 0.44311\ldots\,$, 
a simulation for $p=20$ is shown in Figure \ref{fig:chi1p20}.
The left picture draws the $20000$ roots of $1000$ times of 
realizations of $\bb{S}_{\bb{1}_\D , 20}$ 
(that lie in that coordinate box), and the right picture is the 
density histogram according to the Fubini-Study distance of the 
zeros from the origin $z=0$, where the bell-shape curve represents
the density $c_1(\mathscr{O}(1),h_{\mathrm{FS}})=
\omega_{\mathrm{FS}}$.

From Figure \ref{fig:chi1p20}, with $p=20$, near the origin
(inside the support of $f=\bb{1}_\D$), we see the simulated zeros
approximate $c_1(\mathscr{O}(1),h_{\mathrm{FS}})$ quite well,
but outside the support
(the part $r_\mathrm{FS}> \frac{\sqrt{\pi}}{4}$), the simulated zeros behave
very differently from $c_1(\mathscr{O}(1),h_{\mathrm{FS}})$.
More simulation results will be given in Section \ref{ss:5.4.3new}
to illustrate the convergence results on the support of $f$
from Sections \ref{s1.4a} and \ref{ss:1.4intro}.

From \eqref{eq:1.22intropart}, the asymptotics of $\frac{1}{p}\log {T_{f,p}^2(x,x)}$ is a crucial term to study $\frac{1}{p}\E[[\Div(\bb{S}_{f,p})]]$. By the asymptotic expansion of $T_{f,p}^2(x,x)$ (see Theorem \ref{thm:2.4June}), we conclude that 
\begin{equation}
\log {T_{f,p}^2(x,x)} \leq C\log{p}.
\end{equation}
When $f$ is smooth and $(L,h_L)$ is prequantum (that is $\Theta=c_1(L,h_L)$), the lower bounds for $\log {T_{f,p}^2(x,x)}$ were obtained in \cite{Del2019, Del2020}, \cite{AF:2022}, \cite[Proposition 5.16]{DrLM:2023aa} under the assumption that $f$ vanishes only up to order $2$. But a proper lower bound for $\log {T_{f,p}^2(x,x)}$ in general case is missing. 

 The above expected lower bounds relate clearly to the lowest nonzero eigenvalues of $T^2_{f,p}$. Let us focus on the non-negative function $f$. For a nontrivial $f\in \mathcal{L}^{\infty}(X,\R_{\geq 0})$, $T_{f,p}$ is injective and positive, set $\lambda^p_{\min}(f):=\min \mathrm{Spec}(T_{f,p}) >0$, then on $X$, we get a lower bound for $\frac{1}{p}\log {T_{f,p}^2(x,x)}$:
\begin{equation}
\frac{1}{p}\log {T_{f,p}^2(x,x)}\geq \frac{2\log {\lambda^p_{\min}}(f) }{p} +\frac{\log {P_p(x,x)}}{p}.
\end{equation}
The easy case is that the essential infimum of $f$ is strictly positive on $X$, so that we can conclude $\frac{1}{p}\log {T_{f,p}^2(x,x)}\rightarrow 0$ uniformly on $X$ as $p\rightarrow +\infty$. 

When a nontrivial $f\geq 0$ has a $\mathscr{C}^1$-vanishing point in $X$, 
then $\lambda^p_{\min}(f)$ decays to $0$ as $p\rightarrow +\infty$
(see the analogous statements in Corollary \ref{cor:2.6norm} and
Remark \ref{lem:2.7norm}, see also \cite[Proposition 9.2.1]{Le_Floch_2018}). 
Then a first step to study $\lambda^p_{\mathrm{min}}(f)$ is to understand how 
fast it decays to $0$ when $f$ vanishes at some points in $X$.

In Section \ref{ss5.3Jan},  for the Hermitian line bundle $(L,h_L)=(\mathscr{O}(1),h_{\mathrm{FS}})$ on the Riemann sphere $\mathbb{CP}^1$, 
we have computed explicitly the Toeplitz spectra for three types of functions
and obtain three different asymptotic behavior for $\lambda^p_{\mathrm{min}}$. 
Let $U_0\simeq \C$ denote a standard complex chart for $\mathbb{CP}^1$.
\begin{itemize}
	\item For $k\in\N_{\geq 1}$, set $f_k(z):=\frac{|z|^{2k}}{(1+|z|^2)^k}$
	on $U_0\simeq \C$, then $f_k$ has only one vanishing point at $z=0$
	with vanishing order $2k$. We have
	\begin{equation}
\lambda^p_\mathrm{min}(f_k)=k!p^{-k}(1+\frac{k(k+3)}{2p}+\mathcal{O}(p^{-2})).
\label{eq:5.43JanIntro}
\end{equation}
	\item For $f(z):=e^{-\frac{1}{|z|^2}}$ on $U_0\simeq \C$, 
	then $f$ has only one vanishing point at $z=0$ with vanishing order $\infty$. We have
	\begin{equation}
\lambda^p_\mathrm{min}(f)=e^{-2\sqrt{p}(1+o(1))}.
\end{equation}
\item Let $\mathbb{B}\subset \mathbb{CP}^1$ be a geodesic ball, set 
$f(z):=\bb{1}_{\mathbb{B}}(z)$ the indicator function for $\mathbb{B}$. 
Assume $\overline{\mathbb{B}}\neq \mathbb{CP}^1$, hence $\mathrm{Vol}(\mathbb{B}) <1$. 
We have
	\begin{equation}
\lambda^p_\mathrm{min}(\bb{1}_{\mathbb{B}})=\mathrm{Vol}(\mathbb{B})^{p+1}.
\label{eq:1.46intro}
\end{equation}
\end{itemize}
In Question \ref{question13}, following the above examples, 
we summarize a question for the lowest Toeplitz eigenvalues
for the general compact Hermitian complex manifolds.

Now we present some partial results on $\lambda^p_\mathrm{min}$,
whose proofs will be given in Section \ref{ss:compactandeigen}. 
A complete answer still remains open.
\begin{Proposition}\label{lm:upperbounds}
Let $(X,J,\Theta)$ be a connected, compact Hermitian complex manifold 
and let $(L,h_{L})$, $(E ,h_E)$ be holomorphic line bundles on $X$ 
with smooth Hermitian metrics. Assume $h_L$ to be positive.
Fix $f\in \mathcal{L}^{\infty}(X,\R_{\geq 0})$ which is not identically zero.
\begin{itemize}
	\item[(i)] For $N\in\N_{\geq 1}$, if there exists $x_0\in X$
	such that $f$ is $\mathscr{C}^{2N+1}$ near $x_0$ and $f$ 
	vanishes at $x_0$ with vanishing order $2N$, then there exists
	$C>0$ such that for all $p>0$,
	\begin{equation}
	\min \mathrm{Spec}(T_{f,p})\leq Cp^{-N}.
	\label{eq:5.57march24}
	\end{equation}
	\item[(ii)] If there exists $x_0\in X$ such that $f$ is smooth
	near $x_0$ and $f$ vanishes at $x_0$ with vanishing order $+\infty$, 
	then for any $\ell\in\N$, there exists $C_\ell>0$ such that for all $p>0$,
		\begin{equation}
	\min \mathrm{Spec}(T_{f,p})\leq C_\ell p^{-\ell}.
	\label{eq:5.58march24}
	\end{equation}
\end{itemize}
\end{Proposition}

The following result provides a supportive evidence for the situation 
like \eqref{eq:1.46intro}, which also refines the lower bound in \eqref{eq:4.11Toep} 
in compact case.
\begin{Theorem}\label{thm:4.7Toeplitz}
Let $(X,J,\Theta)$ be a connected, compact Hermitian complex manifold and let 
$(L,h_{L})$, $(E ,h_E)$ be holomorphic line bundles on $X$ with smooth Hermitian metrics. 
Assume $h_L$ to be positive. For $f\in \mathcal{L}^{\infty}(X,\R_{\geq 0})$
which is not identically zero and is continuous near a nonvanishing point, 
there exist constants $C'>0$, $c'>0$ depending only on $X$, $L$, $E$ and 
$f$ such that for all $p\gg 0$,
\begin{equation}
 \min \mathrm{Spec}(T_{f,p})\geq C e^{-cp}
 \label{eq:5.48feb24}
 \end{equation}
If $\esupp f \neq X$, then for any $A>0$, there exists $C'=C'(f,A)>0$ such that for all $p\gg 0$,
\begin{equation}
 \min \mathrm{Spec}(T_{f,p})\leq C' e^{-A\sqrt{p\log{p}}}
 \label{eq:5.48m2}
 \end{equation}
\end{Theorem}


\subsection{Organization of the paper}
This paper is organized as follows:

In Section \ref{ss2nov}, we give the asymptotic expansions of the Berezin--Toeplitz kernels $T^2_{f,p}(x,y)$ under the local regularity assumption on $f$.

In Section \ref{section2}, we recall the definition of Gaussian $\cLL$-holomorphic sections given in \cite[Section 4]{DrLM:2023aa} and the related results.

In Section \ref{ss4Jan}, we prove Theorem \ref{thm:4.6Jan2024}, where the key intermediate result is the Proposition \ref{thm:4.3Jan}. In particular, we give the proof of Proposition \ref{prop:4.19p2} in Section \ref{ss:mass}.

In Section \ref{section5:onsupport}, we study the asymptotic distribution of random zeros on the support of $f$. In particular, the proofs of Theorems \ref{thm:6.2}, \ref{thm:6.3}, \ref{thm:5.5jan24} and \ref{thm:6.2equi} are given. The results on the lowest eigenvalues of $T_{f,p}$ for a compact Hermitian manifold $X$ are given in Section \ref{ss:compactandeigen}.

At last, in Section \ref{section:statistics}, we discuss the number variance
and the asymptotic normality of the zeros of $\bb{S}_{f,p}$
on $\supp f$ for a real smooth function $f$ with compact support.

\subsection*{Acknowledgments}
AD and BL thank NYU Shanghai for their hospitality.
We also thank Xiaonan Ma and St\'ephane Nonnenmacher
for useful discussions.

\section{Toeplitz operators and asymptotics of Toeplitz kernels}\label{ss2nov}
Let $(X,J,\Theta)$ be a connected complex Hermitian (paracompact) manifold of 
complex dimension $n$, where $J$ denotes the canonical complex 
structure of $X$, and $\Theta$ denotes a $J$-compatible Hermitian 
form. Then we have an induced
Riemannian metric $g^{TX}(\cdot,\cdot)=\Theta(\cdot,J\cdot)$ on $X$. We 
denote by $R^{\mathrm{det}}$ the curvature of the Chern connection 
$\nabla^{\mathrm{det}}$ on $K^{*}_{X}:=\det(T^{(1,0)}X)$ with respect 
to induced Hermitian metric by $g^{TX}$. Let $\operatorname{dist}(\cdot,\cdot)$ denote the Riemannian distance 
of $(X, g^{TX})$. Let $(L,h_{L})$, $(E,h_{E})$ be two Hermitian holomorphic line bundle 
on $X$, and let $\nabla^{L}$, $\nabla^{E}$ denote the corresponding Chern 
connections with the respective curvature forms $R^{L}$, $R^{E}$.
We always assume Condition \ref{condt:main} to hold.

\subsection{Bergman projections and the asymptotics of Bergman kernels}\label{ss2.1nov} 
The Riemannian volume form on $(X,g^{TX})$ is denoted by 
$\mathrm{dV}=\frac{\Theta^{n}}{n!}$. For $p\in\mathbb{N}_{>0}$, set $(L^{p}\otimes E, h_{p}):=(L^{\otimes 
p}\otimes E, h_{L}^{\otimes p}\otimes h_E)$. For $s, s'\in 
\mathscr{C}^{\infty}_{\mathrm{c}}(X,L^{p}\otimes E)$, the $\cLL$-inner 
product is defined as follows,
\begin{equation}
	\langle s,s'\rangle_{\cLL(X,L^{p}\otimes E)}:=\int_{X}\langle 
	s(x),s'(x)\rangle_{h_{p,x}} \mathrm{dV}(x).
\end{equation}
Let $\cLL(X, L^{p}\otimes E)$ be the completion of 
$\mathscr{C}^{\infty}_{\mathrm{c}}(X,L^{p}\otimes E)$ with respect 
to the above $\cLL$-inner product. Let $H^{0}(X, L^{p}\otimes E)$ 
denote the space of global holomorphic sections of $L^{p}\otimes E$ 
on $X$. We set
\begin{equation}
	H^{0}_{(2)}(X, L^{p}\otimes E):=\cLL(X,L^{p}\otimes E)\cap 
	H^{0}(X,L^{p}\otimes E).
	\label{eq:2.1.2June}
\end{equation}
Then it is a separable Hilbert subspace of $\cLL(X, L^{p}\otimes E)$. 
Set
\begin{equation}
	d_{p}:=\dim_{\C} H^{0}_{(2)}(X, L^{p}\otimes E) \in\N\cup\{\infty\}.	
	\label{eq:dimensionJune}
\end{equation}

Let $P_{p}: \cLL(X,L^{p}\otimes E)\rightarrow 
H^{0}_{(2)}(X,L^{p}\otimes E)$ denote the obvious orthogonal 
projection, which is called the Bergman projection. It has a smooth 
Schwartz integral kernel, denoted by $P_{p}(x,x')\in (L^{p}\otimes 
E)_{x}\otimes (L^{p}\otimes E)^{\ast}_{x'}$. Following the work of 
Ma-Marinescu \cite[Chapters 4 \& 6]{MM07}, we have the following 
results on the asymptotics of Bergman kernels (under the assumption 
\eqref{eq:mainassumptions}):
\begin{itemize}
	\item(Off-diagonal estimates) For any $\ell,m\in \mathbb{N}$, 
	$\varepsilon>0$, a 
	compact subset $K\subset X$, there exists 
	$C_{K,\ell,m,\varepsilon}>0$ such that for all $x,x'\in K$, 
	$\mathrm{dist}(x,x')\geq \varepsilon$, we have
	\begin{equation}
		|P_{p}(x,x')|_{\mathscr{C}^{m}}\leq C_{K,\ell,m,\varepsilon} 
		p^{-\ell}.
		\label{eq:2.1.3June}
	\end{equation}
	Here the $\mathscr{C}^{m}$-norm is induced by $\nabla^{L}$, 
	$\nabla^{E}$ and $h_{p}$, $g^{TX}$.
	\item(On-diagonal expansion) There exist coefficients 
	$\bb{b}_{r}\in\mathscr{C}^{\infty}(X,\C)$, $r\in\N$, such 
	that for any compact subset $K\subset X$, any $k,\ell\in \N$, 
	there exists $C_{k,\ell,K}>0$ such that for $p\in\N^{\ast}$,
	\begin{equation}
		\left|\frac{1}{p^{n}}P_{p}(x,x)-\sum_{r=0}^{k}\bb{b}_{r}(x)p^{-r}\right|_{\mathscr{C}^{\ell}(K)}\leq C_{k,\ell,K} p^{-k-1},
		\label{eq:2.1.4June}
	\end{equation}
where 
	$\bb{b}_{0}(x)=\det\left(\frac{\dot{R}^{L}}{2\pi}\right)$ (recall \eqref{eq:Rdot} for the definition of $\dot{R}^{L}$), and an explicit formula for
$\bb{b}_{1}$ was given as in 
\cite[(4.1.9)]{MM07}.
\end{itemize}
Furthermore, the near-diagonal expansion for $P_{p}(x,x')$ also holds 
uniformly on any given compact subset of $X$. To describe this 
expansion, we need to introduce some notation as follows.

Fix a point $x_{0}\in X$. Let $\{\bb{f}_{j}\}_{j=1}^{n}$ be an orthonormal basis of 
$(T_{x_{0}}^{1,0}X, g_{x_{0}}^{TX}(\cdot,\overline{\cdot}))$ such that
\begin{equation}
	\dot{R}^{L}_{x_{0}}\,\bb{f}_{j}=\mu_{j}(x_{0})\bb{f}_{j},
\end{equation}
where $\mu_{j}(x_{0})$, $j=1, \ldots, n$ are the eigenvalues of 
$\dot{R}^{L}_{x_{0}}$. We have
\begin{equation}\label{eq:b0}
	\mu_{j}(x_{0})\geq \varepsilon_{0},\; 
	\bb{b}_{0}(x_{0})=\prod_{j=1}^{n}\frac{\mu_{j}(x_{0})}{2\pi}.
\end{equation}
Set 
$\bb{e}_{2j-1}=\frac{1}{\sqrt{2}}(\bb{f}_{j}+\overline{\bb{f}}_{j})$, 
$\bb{e}_{2j}=\frac{\sqrt{-1}}{\sqrt{2}}(\bb{f}_{j}-\overline{\bb{f}}_{j})$, $j=1, \ldots, n$. 
Then they form an orthonormal basis of the (real) tangent vector 
space $(T_{x_{0}}X, g_{x_{0}}^{TX})$. Now we introduce the complex 
coordinate for $T_{x_{0}}X$ (with respect to the complex structure 
$J_{x_{0}}$). If $v=\sum_{j=1}^{2n} v_{j}\bb{e}_{j}\in T_{x_{0}}X$, we can write
\begin{equation}
	v=\sum_{j=1}^{n}(v_{2j-1}+\sqrt{-1}v_{2j})\frac{1}{\sqrt{2}} 
	\bb{f}_{j} + \sum_{j=1}^{n}(v_{2j-1}-\sqrt{-1}v_{2j})\frac{1}{\sqrt{2}} 
	\overline{\bb{f}}_{j}.
	\label{eq:5.1.8parisJune}
\end{equation}
Set $z=(z_{1},\ldots,z_{n})$ with $z_{j}=v_{2j-1}+\sqrt{-1}v_{2j}$, 
$j=1,\ldots, n$. We call $z$ the complex coordinate of $v\in T_{x_{0}}X$. 
Then by \eqref{eq:5.1.8parisJune},
\begin{equation}
	\frac{\partial}{\partial z_{j}} = \frac{1}{\sqrt{2}} 
	\bb{f}_{j} , \; \frac{\partial}{\partial \overline{z}_{j}} = \frac{1}{\sqrt{2}} 
	\overline{\bb{f}}_{j},
	\label{eq:5.1.9parisJune}
\end{equation}
so that
\begin{equation}
	v=\sum_{j=1}^{n}\Big(z_{j} \frac{\partial}{\partial z_{j}} 
	+\overline{z}_{j} \frac{\partial}{\partial \overline{z}_{j}}\Big).
	\label{eq:5.1.10parisJune}
\end{equation}
Note that $|\frac{\partial}{\partial z_{j}}|^{2}_{g_{x_{0}}^{TX}}=|\frac{\partial}{\partial 
\overline{z}_{j}}|^{2}_{g_{x_{0}}^{TX}}=\frac{1}{2}$.
For $v, v'\in T_{x_{0}}X$, let $z, z'$ denote the corresponding 
complex coordinates. Define
\begin{equation}
\mathcal{P}_{x_{0}}(v,v')=\prod_{j=1}^{n} 
\frac{\mu_{j}(x_{0})}{2\pi}
\exp\Big(\!-\frac{1}{4}\sum_{j=1}^{n}\mu_{j}(x_{0})(|z_{j}|^{2}
+|z'_{j}|^{2}-2z_{j}\overline{z}'_{j})\Big).
	\label{eq:5.1.11paris}
\end{equation}	
Define a weighted distance function $\Phi^{TX}_{x_{0}}(v,v')$ as follows,
\begin{equation}
\Phi^{TX}_{x_{0}}(v,v')^{2}=\sum_{j=1}^{n}\mu_{j}(x_{0})|z_{j}-z'_{j}|^{2}.
\label{eq:2.16July2}
\end{equation}
Then 
\begin{equation}
|\mathcal{P}_{x_{0}}(v,v')|=\prod_{j=1}^{n} 
\frac{\mu_{j}(x_{0})}{2\pi}\exp\Big(\!-\frac{1}{4}\Phi^{TX}_{x_{0}}(v,v')^{2}\Big).
\label{eq:2.17July}
\end{equation}
For sufficiently small $\delta_{0}>0$, we identify the small open 
ball
$B^{X}(x_{0},2\delta_{0})$ in $X$ with the ball $B^{T_{x}X}(0, 
2\delta_{0})$ in $T_{x_{0}}X$ via the geodesic coordinate. Let 
$\kappa(v)$ be the positive smooth function such that 
\begin{equation}
	\mathrm{dV}(\exp_{x_{0}}(v))=\kappa(v)\mathrm{dV}_{\mathrm{Eucl}}(v),
	\label{eq:2.18July}
\end{equation}
where $\mathrm{dV}_{\mathrm{Eucl}}$ denotes the Euclidean volume form 
on $T_{x_{0}}X$ with respect to $g^{TX}_{x_{0}}$. In particular, $\kappa(0)=1$.

There exists $C_{2}>0$ such that for $v,v'\in 
B^{T_{x_{0}}X}(0, 2\delta_{0})$, we 
have
\begin{equation}
	C_{2}\operatorname{dist}(\exp_{x_{0}}(v),\exp_{x_{0}}(v'))\geq 
	\Phi^{TX}_{x_{0}}(v,v')\geq 
	\frac{1}{C_{2}}\operatorname{dist}(\exp_{x_{0}}(v),\exp_{x_{0}}(v')).
	\label{eq:5.1.15DLM}
\end{equation}
In particular,
\begin{equation}
	\Phi^{TX}_{x_{0}}(0,v)\geq \varepsilon_{0}^{1/2} 
	\operatorname{dist}(x_{0},\exp_{x_{0}}(v))=\varepsilon_{0}^{1/2} 
	\Vert v\Vert.
	\label{eq:5.1.15paris}
\end{equation}
Moreover, if we consider a compact subset $K\subset X$, the 
constants $\delta_{0}$ and $C_{1}$ can be chosen uniformly for all 
$x_{0}\in K$. Similarly, by \eqref{eq:1.09intro}, on a compact subset $K$, we have
\begin{equation}
	\Phi^{TX}_{x_{0}}(0,v)^2\leq \kappa(R^L, K) 
	\Vert v\Vert^2, \; x_0\in K.
	\label{eq:2.21parisnew}
\end{equation}

We trivialize the line bundle $L$ on $B^{T_{x_{0}}X}(0, 2\delta_{0})$ 
using the parallel transport with respect to $\nabla^{L}$ along 
the curve $[0,1]\ni t\mapsto tv$, $v\in B^{T_{x_{0}}X}(0, 
2\delta_{0})$. Under this trivialization, for $v,v' \in B^{T_{x_{0}}X}(0, 
2\delta_{0})$, 
\begin{equation}
	P_{p}(\exp_{x_{0}}(v),\exp_{x_{0}}(v'))\in\mathrm{End}(L_{x_{0}}^{p}\otimes E_{x_{0}})=\C.
\end{equation}

By \cite[Theorems 4.2.1 \& 
6.1.1]{MM07}, for any compact subset $K\subset X$, there exists a constant $C'>0$ so that for any $\ell, 
m, N\in \mathbb{N}$, there 
exists $\delta>0$ and constant $C=C(K,\ell,m,N)>0$ such that for 
$x\in K$, $v,v'\in T_{x}X$, multi-indices $\alpha$, $\alpha'\in 
\mathbb{N}^{2n}$ with $|\alpha|+|\alpha'|\leq m$
$\Vert v\Vert,\,\Vert v'\Vert \leq 2\delta$, we have
\begin{equation}
	\begin{split}
		&\bigg|\frac{\partial^{|\alpha|+|\alpha'|}}{\partial 
		v^{\alpha}\partial (v')^{\alpha'}}\left(\frac{1}{p^{n}}P_{p}(\exp_{x}(v),\exp_{x}(v'))-\sum_{r=0}^{N}\mathcal{F}_{r}(\sqrt{p}v,\sqrt{p}v')\kappa^{-1/2}(v)\kappa^{-1/2}(v')p^{-r/2}\right)\bigg|_{\mathscr{C}^{\ell}(K)}\\
		&\leq 
		Cp^{-(N+1-m)/2}(1+\sqrt{p}\Vert v\Vert+\sqrt{p}\Vert 
		v'\Vert)^{2(N+n+\ell)+2+m}\exp(-C'\sqrt{p}\Vert v-v'\Vert)+\mathcal{O}(p^{-\infty}).
	\end{split}
	\label{eq:5.2.2paris}
\end{equation}
The functions $\mathcal{F}_{r}$, $r\in\N$ are given as follows,
\begin{equation}
	\mathcal{F}_{r}(v,v')=\mathcal{P}_{x}(v,v')\mathcal{J}_{r}(v,v'),
	\label{eq:5.1.19DLM}
\end{equation}
where $\mathcal{J}_{r}(v,v')$ is a polynomial in $v,v'$ of degree $\leq 3r$, whose 
coefficients are smooth in $x\in X$. In particular,
\begin{equation}
	\mathcal{J}_{0}=1.
\end{equation}
The notation $\mathcal{O}(p^{-\infty})$ means that this term is 
bounded uniformly by $C^{\prime\prime}p^{-k}$ for any given $k\in\N$ 
and with some constant $C^{\prime\prime}>0$ which is independent of the 
choices of
$x\in K$,  
$v, v'\in T_{x}X$ and $p$ involved in \eqref{eq:5.2.2paris}.

One crucial step in the approach of Ma-Marinescu \cite[Chapters 4 \& 
6]{MM07} for the Bergman kernel expansions is the localization the 
calculation of the Bergman kernel near one given point \cite[Section 
4.1]{MM07} due to the 
finite propagation speed of solutions of hyperbolic equations. For our 
convenience of explaining the proofs of Toeplitz kernels in next 
Section, we sketch this localization technique in the sequel.

Let $\dbar_{p}$ denote the $\dbar$-operator on 
$\mathscr{C}^{\infty}_{\mathrm{c}}(X, 
\Lambda^{\bullet}T^{*(0,1)}X\otimes L^{p}\otimes E)$, and let $\dbar^{\ast}_{p}$ 
denote its formal adjoint with respect to the $\cLL$-inner product. 
We always take the maximal extensions of $\dbar_{p}$, 
$\dbar^{\ast}_{p}$ as differential operators on 
$\cLL(X,\Lambda^{\bullet}T^{*(0,1)}X\otimes L^{p}\otimes 
E)$. Since $(X,g^{TX})$ is assumed to be complete, then the maximal 
extension of $\dbar^{\ast}_{p}$ coincides with the Hilbert adjoint of 
the maximal extension of $\dbar_{p}$. We still use the same notation to 
denote the above operators.

The Kodaira Laplacian 
$\Box^{E}_{p}=\dbar_{p}\dbar^{\ast}_{p}+\dbar^{\ast}_{p}\dbar_{p}$ is 
a densely defined, positive operator. In our setting, it has a unique 
self-adjoint extension, denoted also by $\Box^{E}_{p}$ and called Gaffney extension, whose domain is 
given by
\begin{equation}
	\mathrm{Dom}(\Box^{E}_{p})=\{s\in \mathrm{Dom}(\dbar_{p})\cap 
	\mathrm{Dom}(\dbar^{\ast}_{p})\;:\; \dbar_{p}s\in 
	\mathrm{Dom}(\dbar_{p}^{\ast}),\; \dbar^{\ast}_{p}s\in 
	\mathrm{Dom}(\dbar_{p})\}.
\end{equation}
Then we have
\begin{equation}
	H^{0}_{(2)}(X,L^{p}\otimes E)=\ker (\dbar_{p}|_{\cLL(X, L^{p}\otimes 
E)})=\ker(\Box^{E}_{p}|_{\cLL(X, L^{p}\otimes 
E)}).
\end{equation}

The 
assumptions in \eqref{eq:mainassumptions} implies a spectral gap 
\cite[(6.1.8)]{MM07} so that there exists $C>0$, $C'>0$ such that
\begin{equation}
	\mathrm{Spec}(\Box^{E}_{p})\subset \{0\}\cup [ Cp-C',+\infty[\,.
	\label{eq:2.27feb24}
\end{equation}
Moreover, the higher Dolbeault 
$\cLL$-cohomology groups vanish for $p\gg 0$.
We are mainly concerned with $\Box^{E}_{p}|_{\cLL(X, L^{p}\otimes 
E)}$, which will be denoted simply by $\Box^{E}_{p}$ in the sequel.

Recall that the injectivity radius $\mathrm{inj}^{X}_{U}$ is defined in \eqref{eq:1.37intro}. Fix $\delta\in [0,\mathrm{inj}^{X}_{U}/4[\;$. Let $h: \R\rightarrow [0,1]$ be an even 
smooth function such that
\begin{equation}	
	h(s)=\begin{cases}
	1 & \text{for}\quad |s|\leq \delta/2,\\
	0 & \text{for}\quad |s|\geq \delta.
	\end{cases}
	\label{eq:2.29June}
\end{equation}
Set
\begin{equation}\label{eq:defH(a)}
H(a)=\left(\int_{-\infty}^{+\infty}h(s)ds\right)^{-1}
\int_{-\infty}^{+\infty}e^{\sqrt{-1}sa}h(s)ds.
\end{equation}
Then $H(a)$ is an analytic even function that also lies in the 
Schwartz space $\mathcal{S}(\R)$ and $H(0)=1$.

We always consider the integer $p$ to be such that $Cp-C'\geq 0$ (the 
constants are from the spectral gap \eqref{eq:2.27feb24}). Set
\begin{equation}\label{eq:defphip}
	\phi_{p}(a)=\mathrm{1}_{[\sqrt{Cp-C'},+\infty[\,}(|a|)H(a).
\end{equation}
It is still an even function. Note that by the functional calculus, 
we have the operators $H(D_{p})$, $\phi_{p}(D_{p})$ well-defined as 
bounded operators on $\cLL(X,L^{p}\otimes E)$ with smooth integral 
kernels. In particular, if $p>C'/C$, we have
\begin{equation}
	H(D_{p})-P_{p}=\phi_{p}(D_{p}).
	\label{eq:2.32June}
\end{equation}

If $Q,Q'$ are two differential operators of order $m, m'$ with 
compact support respectively, then for any $\ell>0$, there exists 
$C_{\ell}>0$ such that for $p\geq C'/C$, we have for any 
$s\in\mathscr{C}^{\infty}_{\mathrm{c}}(X, L^{p}\otimes E)$,
\begin{equation}
	\Vert Q\phi_{p}(D_{p})Q's\Vert_{\cLL}\leq C_{\ell} p^{-\ell} 
	\Vert s\Vert_{\cLL}.
	\label{eq:2.33June}
\end{equation}
As a consequence, the Bergman projection $P_{p}$ can be approximated 
by $H(D_{p})$ up to an reminder of $\mathcal{O}(p^{-\ell})$ as well 
as in the level of their integral kernels. In particular, 
$H(D_{p})(x,x')$ only depends on the restriction of $D_{p}$ to 
$B^{X}(x,\delta)$, and we have
\begin{equation}
	H(D_{p})(x,x')=0, \;\text{for\;} \mathrm{dist}(x,x')\geq \delta.
	\label{eq:2.34June}
\end{equation}
As a consequence, the off-diagonal estimate \eqref{eq:2.1.3June} 
follows from \eqref{eq:2.32June}, \eqref{eq:2.33June} and 
\eqref{eq:2.34June}. For the on-diagonal expansion 
\eqref{eq:2.1.4June} and the near-diagonal expansion 
\eqref{eq:5.2.2paris}, the computation localizes to the small ball 
$B^{X}(x,\delta)$. The details are referred to \cite[Section 
4.1]{MM07}.

Finally, we recall the off-diagonal estimates for $P_p$ on a compact complex manifold.
\begin{Theorem}[{\cite[Theorem 1]{MM15} and \cite[Theorem 1]{MR3903323}}]\label{thm:offdiagonal}
Let $(X,J,\Theta)$ be a connected, compact Hermitian complex manifold and let $(L,h_{L})$, $(E,h_E)$ be holomorphic line bundles on $X$ with smooth Hermitian metrics. Assume $h_L$ to be positive. 
Then there exist constants $\bb{c}>0$ and $p_0\in\N$ such that for any $m\in \N$, there exists a constant $C_m>0$ such that for $p\geq p_0$, and for all $x,y\in X$,
\begin{equation}
|P_p(x,y)|_{\mathscr{C}^m}\leq C_m p^{n+m/2}e^{-\bb{c}\sqrt{p}\mathrm{dist}(x,y)}.
\label{eq:2.35feb24}
\end{equation}

Moreover, for any $\delta>0$, $A>0$, there exist constants $p_0\in \N$, $C_{\delta,A}>0$ such that for all $x,y\in X$ with $\mathrm{dist}(x,y)\geq \delta$, $p\geq p_0$, we have
\begin{equation}
|P_p(x,y)|_{h_{p,x}\otimes h^*_{p,y}}\leq C_{\delta,A} e^{-A\sqrt{p\log{p}}}.
\label{eq:2.36feb24}
\end{equation}
\end{Theorem}

The first part of the above theorem was proved by 
Ma--Marinescu \cite[Theorem 1]{MM15}, and their result 
holds for general (noncompact) complete manifolds with 
bounded geometry. Under the off-diagonal assumption 
$\mathrm{dist}(x,y)\geq \delta$, for the case where 
$m=0$, $\Theta=c_1(L,h_L)$ and $E=\C$ trivial line bundle, 
\eqref{eq:2.35feb24} is called Agmon estimate, 
and also follows from the parametix of the Szeg\H{o} kernels on
pseudo-convex domains (such as \cite{MR1656004, MR1828799}), 
see also \cite[Theorem 3.1]{MR3574652}. 
The sharper off-diagonal estimate \eqref{eq:2.35feb24} 
was proved by Christ \cite{MR3903323} by studying the 
near-diagonal estimate of the Green kernels for 
Kodaira Laplacians acting on $(0,1)$-forms.

\subsection{Berezin--Toeplitz quantization}\label{ss:2.2feb}
We have defined the following spaces of bounded measurable functions 
in the Introduction: 
$$\mathcal{L}^{\infty}_{\mathrm{c}}(X,\C)
\subset \mathcal{L}^{\infty}_{\mathrm{const}}(X,\C)
\subset \mathcal{L}^{\infty}(X,\C).$$

\begin{Definition}\label{def:4.4}
	For any $f\in \mathcal{L}^{\infty}(X,\C)$, set
	\begin{align}
T_{f,p}: H^{0}_{(2)}(X, L^{p}\otimes E) \rightarrow 
H^{0}_{(2)}(X,L^{p}\otimes E),\quad
		T_{f,p}:=P_{p}M_{f},
		\label{eq:2.1.1}
	\end{align}
where $M_f$ denotes the pointwise multiplication by $f$. 
The family of bounded operators $\{T_{f,p}\}_{p\in\N^{\ast}}$
is called Toeplitz operator associated 
to the symbol $f$, and we call $T_{f,p}$ the
Toeplitz operator of level $p$. Equivalently, $T_{f,p}$ can be seen
as an operator on $\mathcal{L}^{2}$-spaces,
$T_{f,p}:\mathcal{L}^{2}(X, L^{p}\otimes E) 
\rightarrow\mathcal{L}^{2}(X,L^{p}\otimes E)$, 
$T_{f,p}:=P_{p}M_{f}P_{p}$.
\end{Definition}

The map which associates to a function $f$ the operator $T_{f,p}$ on 
$\cLL(X, L^{p}\otimes E)$ is the Berezin--Toeplitz 
quantization of level $p$. The map 
$\mathcal{L}^{\infty}(X,\C)\ni f\mapsto 
\{T_{f,p}\}_{p\in\N^{\ast}}$ is called the Berezin--Toeplitz 
quantization \cite{BMS94,MM07,MM08,Pol18,MR1852228}. 
Clearly, we have
\begin{equation}
	\Vert T_{f,p} \Vert_{\mathrm{op}}\leq \Vert f 
	\Vert_{\mathcal{L}^{\infty}},
\end{equation}
where $\Vert T_{f,p} \Vert_{\mathrm{op}}$ denotes the operator norm 
of $T_{f,p}$.

For $f\in \mathcal{L}^{\infty}(X,\C)$ and $p\in 
\N^{\ast}$, $T_{f,p}$ always has a smooth integral kernel given by
\begin{equation}
	T_{f,p}(x,x')=\int_{X} 
P_{p}(x,x^{\prime\prime})f(x^{\prime\prime})P_{p}(x^{\prime\prime},x')
\mathrm{dV}(x^{\prime\prime}).
	\label{eq:2.1.2}
\end{equation}
Note also that the Hilbert adjoint of $T_{f,p}$ is 
$T_{\overline{f},p}$. If $f$ has compact support, 
an easy modification of the arguments in \cite[Proof of Proposition 
4.7]{DrLM:2023aa} shows that all $T_{f,p}$, $p\in \N^{*}$, are 
Hilbert-Schmidt. 

In \cite[Chapter 7]{MM07}, the asymptotic expansion of 
$T_{f,p}(x,x')$ as $p\rightarrow +\infty$ has been studied in 
detail with the assumption that $f$ is smooth on $X$. However, 
without the global smoothness of $f$, the same arguments
presented in \cite[Sections 7.2 \& 7.5]{MM07} can still be utilized to obtain 
the analogues of \cite[Lemmas 7.2.2 \& 7.2.4, Theorem 
7.5.1]{MM07}. Note that in \cite{BMMP14} 
Toeplitz operators with $\mathscr{C}^{k}$ symbol
were considered (see also \cite{CP17}), and the asymptotics of their kernels 
were established using the arguments from \cite[Sections 
7.2 \& 7.5]{MM07}. In particular, the first part of the following 
theorem for compact $X$ was already given in \cite[Lemma 
3.1]{BMMP14}.

\begin{Theorem}\label{thm:2.2June}
	Assume the geometric setting as in Condition \ref{condt:main}.
	Given $f\in \mathcal{L}^{\infty}_{\mathrm{const}}(X,\C)$. Then for 
	any compact subset $K\subset X$ and $m\in\N$, there exists 
$C_{K,m}>0$ such that for $x\in K$ and for all $p\geq 1$, we have 
the on-diagonal 
estimate
\begin{equation}
	\left|T_{f,p}(x,x)\right|_{\mathscr{C}^{m}(K)}\leq 
	C_{K,m} p^{n+\frac{m}{2}}.
	\label{eq:June2.38on}
\end{equation}
	For any $m,\ell\in\N$, $\varepsilon >0$, 
and a compact subset $K\subset X$, there exists 
$C_{K,m,\ell,\varepsilon}>0$ such that for $x,x'\in K$ with 
$\mathrm{dist}(x,x')\geq \varepsilon$ and for all $p\geq 1$, we have the off-diagonal 
estimate
\begin{equation}
	\left|T_{f,p}(x,x')\right|_{\mathscr{C}^{m}}\leq 
	C_{K,m,\ell,\varepsilon} p^{-\ell}.
	\label{eq:June2.38}
\end{equation}
Let $U\subset X$ be an open 
	subset such that that $f$ is smooth on $U$. There 
	exists a family of polynomials $\{Q_{r,x_{0}}(f)\in \C[v,v'], 
	v,v'\in T_{x_{0}}X\}_{x_{0}\in U}$ with the same parity as $r$, 
	and smooth in $x_{0}\in U$ such that for any compact subset 
	$K\subset U$, there exist $\delta_{K}>0$, $C'>0$ such that for any 
	$\ell,m,N\in\N$, $v,v'\in T_{x_{0}}X$, $\Vert v\Vert, \Vert v'\Vert 
	\leq \delta_{K}$, multi-indice $\alpha$, $\alpha'\in 
\mathbb{N}^{2n}$ with $|\alpha|+|\alpha'|\leq m$, there exist 
$C_{K,\ell,m,N}>0$, 
	$M_{K,\ell,m,N}\in\N$ such that
	\begin{equation}
		\begin{split}
			&\bigg|\frac{\partial^{|\alpha|+|\alpha'|}}{\partial 
			v^{\alpha}\partial (v')^{\alpha'}}\bigg(\frac{1}{p^{n}}T_{f,p}(\exp_{x_{0}}(v),\exp_{x_{0}}(v'))\\
			&\qquad\qquad\qquad-\sum_{r=0}^{N} 
			\left(Q_{r,x_{0}}(f)\mathcal{P}_{x_{0}}\right)(\sqrt{p}v,\sqrt{p}v')\kappa^{-1/2}(v)\kappa^{-1/2}(v')p^{-r/2}\bigg)\bigg|_{\mathscr{C}^{\ell}(K)}\\
		&\qquad\leq 
		C_{K,\ell,m,N}p^{-(N+1-m)/2}(1+\sqrt{p}\Vert v\Vert+\sqrt{p}\Vert 
		v'\Vert)^{M_{K,\ell,m,N}}\exp(-C'\sqrt{p}\Vert 
		v-v'\Vert)+\mathcal{O}(p^{-\infty}).
		\end{split}
		\label{eq:2.27June}
	\end{equation}
In particular, we have $Q_{0,x_{0}}(f)=f(x_{0})$. 
For $x_0\in U$, if $f$ vanishes at $x_0$ with vanishing order $N\in \N$, then we have 
$Q_{r,x_{0}}(f)\equiv 0$ for $r\leq N-1$. If 
	$x_{0}\not\in \esupp f$, then we have 
	$Q_{r,x_{0}}(f)\equiv 0$ for $r=0,1,\cdots$. 
\end{Theorem}
\begin{proof}
By considering the operator $\phi_{p}(D_{p})$ defined 
in previous Section (cf.\ \eqref{eq:defphip}, \eqref{eq:2.32June}), 
we get by \eqref{eq:2.32June}, 
\eqref{eq:2.33June} and \eqref{eq:2.1.1}, 
\begin{equation}
	T_{f,p}=H(D_{p})M_{f}H(D_{p})+\mathcal{O}(p^{-\infty}).
	\label{eq:2.38June}
\end{equation}
For any compact subset $K\subset 
X$, for $x,x'\in K$, we have
\begin{equation}
	\begin{split}
		T_{f,p}(x,x')&=\left(H(D_{p})M_{f}H(D_{p})\right)(x,x')+\mathcal{O}_{K}(p^{-\infty})\\
		&=\int_{X} 
H(D_{p})(x,x^{\prime\prime})f(x^{\prime\prime})H(D_{p})(x^{\prime\prime},x')
\mathrm{dV}(x^{\prime\prime})+\mathcal{O}_{K}(p^{-\infty}).
	\end{split}
	\label{eq:2.1.43June}
\end{equation}
Recalling the definitions \eqref{eq:2.29June}, 
\eqref{eq:defH(a)} of the functions $h$ and $H$, we have
by \eqref{eq:2.34June},
\begin{equation}
\left(H(D_{p})M_{f}H(D_{p})\right)(x,x')=0, \;\text{for\;} 
\mathrm{dist}(x,x')\geq 2\delta.
	\label{eq:2.39June}
\end{equation}
Then, using the same arguments from the proof of \eqref{eq:2.1.3June} 
from the estimate \eqref{eq:2.33June}, we get \eqref{eq:June2.38}.

Now let us consider a compact subset $K\subset X$ and
take take $0<\delta< \mathrm{inj}^{X}_{K}/4$ 
in the definition of function $h$ in \eqref{eq:2.29June}. 
Then for $x_{0}\in K$ and $x,x'\in 
B^{X}(x_{0},\delta)$ the asymptotic expansion of 
$T_{f,p}(x,x')$ is the same as $\left(H(D_{p})M_{f}H(D_{p})\right)(x,x')$, 
and as in \cite[Proof of Lemma 7.2.4]{MM07}, the computations of the expansion only depend 
on the values of $f$ on $B^{X}(x_{0}, 2\delta)\subset U$, then the 
lack of global smoothness of $f$ will not make any difference on this 
computations near $x_{0}$. More precisely, we consider the point 
$x_{0}\in K$, we trivialize the line bundles $L$, $E$ on the small ball 
$B^{X}(x_{0}, 3\delta)$ along the radial geodesics from the center 
$x_{0}$ with respect to their Chern connections, so that locally the 
line bundles $L$, $E$ are identified with the trivial line bundles 
given by $L_{x_{0}}$, $E_{x_{0}}$ respectively. Under this 
trivialization, for $v,v'\in T_{x_{0}}X$, $\Vert v \Vert,\Vert v'\Vert 
<3\delta$, we have
\begin{equation}
	P_{p}(\exp_{x_{0}}(v),\exp_{x_{0}}(v'))\in 
	\mathrm{End}(L_{x_{0}}^{p}\otimes E_{x_{0}})=\C.
\end{equation}
We regard $P_{p,x_{0}}(v,v'):= P_{p}(\exp_{x_{0}}(v),\exp_{x_{0}}(v'))$, ($v,v'\in T_{x_{0}}X, \Vert v \Vert,\Vert v'\Vert 
<3\delta$), as a smooth section function over $TX\times_{K} TX$. We 
refer to \cite[Sections 4.1.5 \& 4.2.1]{MM07} for more details.

Let $\rho:\R\rightarrow 
[0,1]$ be a smooth even function such that
\begin{equation}
	\rho(s)=\begin{cases}
		1 & \quad\text{if}\quad |s|<2, \\
		0 & \quad\text{if}\quad |s|>4.
	\end{cases}
\end{equation}
Then for $v,v'\in T_{x_{0}}X$ with $\Vert v\Vert, \Vert v'\Vert < 
\delta$, we have
\begin{equation}
	\begin{split}
			&T_{f,p}(\exp_{x_{0}}(v), \exp_{x_{0}}(v'))\\
			&=\int_{v''\in T_{x_{0}} 
	X }  
	P_{p,x_{0}}(v,v^{\prime\prime})\rho(\Vert 
	v^{\prime\prime}\Vert/\delta)f(\exp_{x_{0}}(v^{\prime\prime}))P_{p,x_{0}}(v^{\prime\prime}, v')\kappa_{x_{0}}(v^{\prime\prime})\mathrm{d}v^{\prime\prime}+\mathcal{O}_{K}(p^{-\infty}).
	\end{split}
	\label{eq:2.45July}
\end{equation}
where $\kappa_{x_{0}}$ is given as the function $\kappa$ in 
\eqref{eq:2.18July} but we put the subscript $x_{0}$
 to indicate the base point $x_{0}\in K$. Note that in 
 \eqref{eq:2.45July}, we do not need $f$ to be smooth near $x_{0}$.
 
To conclude \eqref{eq:June2.38on}, we can proceed as in the proof of \cite[Lemma 
3.1]{BMMP14}. Indeed, the derivatives on variable $x_{0}$ at a given 
point $x_{0}$ can be 
replaced by the derivatives on $v$ and $v'$ in \eqref{eq:2.45July}, 
which eventually makes derivatives on $v$ or $v'$ of 
$P_{p,x_{0}}(v,v')$, then the factor $p^{n+m/2}$ follows from 
\eqref{eq:5.2.2paris}.

Now we assume that $f$ is smooth on $U$ and $K$ is a compact subset 
of $U$. In the above arguments to obtain  \eqref{eq:2.45July}, we 
take $0<\delta< \min\{\mathrm{inj}^{X}_{K}, \mathrm{dist}(K, 
X\setminus U)\}/4$. 
Following exactly the same arguments in \cite[Proof of Lemma 
7.2.4]{MM07}, the identity \eqref{eq:2.45July} gives the expansion 
\eqref{eq:2.27June} with $Q_{0,x_{0}}(f)=f(x_{0})$ for the case $m=0$. 

For general $m\in\N$, we note that the expansions of 
$\frac{\partial^{|\alpha|+|\alpha'|}}{\partial v^{\alpha}\partial (v')^{\alpha'}}T_{f,p}(\exp_{x_{0}}(v),\exp_{x_{0}}(v'))$ is given by computing the following integration
\begin{equation}
	p^{n}\int_{v''\in T_{x_{0}} 
	X }  
	\frac{\partial^{|\alpha|}}{\partial v^{\alpha}}P_{p,x_{0}}(v,v^{\prime\prime})\kappa_{x_{0}}^{1/2}(v^{\prime\prime})\rho(\Vert 
	v^{\prime\prime}\Vert/\delta)f(\exp_{x_{0}}(v^{\prime\prime}))\frac{\partial^{|\alpha'|}}{\partial (v')^{\alpha'}}P_{p,x_{0}}(v^{\prime\prime}, v')\kappa_{x_{0}}^{1/2}(v^{\prime\prime})\mathrm{d}v^{\prime\prime},
	\label{eq:2.47July2}
\end{equation}
where we apply \eqref{eq:5.2.2paris} for 
$\frac{\partial^{|\alpha|}}{\partial 
v^{\alpha}}P_{p,x_{0}}(v,v^{\prime\prime})\kappa_{x_{0}}^{1/2}(v^{\prime\prime})$ and $\frac{\partial^{|\alpha'|}}{\partial (v')^{\alpha'}}P_{p,x_{0}}(v^{\prime\prime}, v')\kappa_{x_{0}}^{1/2}(v^{\prime\prime})$. Then \eqref{eq:2.27June} for general $m\geq 0$ follows from the analogous arguments as in \cite[Proof of Lemma 
7.2.4]{MM07} (also cf. \cite[Proof of Theorem 3.3]{BMMP14}) and a 
simple observation on \cite[(7.1.6) in Lemma 7.1.1]{MM07}: for 
polynomials $F, G\in \C[v,v']$ there exists a polynomial 
$\mathscr{K}[F,G]\in \C[v,v']$ such that 
\begin{equation}
	\begin{split}
			\int_{v^{\prime\prime}\in \R^{2n}} 
	F(v,v^{\prime\prime})\mathcal{P}(v,v^{\prime\prime})G(v^{\prime\prime},v')\mathcal{P}(v^{\prime\prime},v')\mathrm{d}v^{\prime\prime}
	=\mathscr{K}[F,G](v,v')\mathcal{P}(v,v'),
	\end{split}
	\label{eq:2.49march24}
\end{equation}
 and for multi-indice $\alpha$, $\alpha'\in 
\mathbb{N}^{2n}$, we have
\begin{equation}
	\begin{split}
			\int_{v^{\prime\prime}\in \R^{2n}} 
	\frac{\partial^{|\alpha|}}{\partial 
	v^{\alpha}}\left(F(v,v^{\prime\prime})\mathcal{P}(v,v^{\prime\prime})\right)\frac{\partial^{|\alpha'|}}{\partial (v')^{\alpha'}}\left(G(v^{\prime\prime},v')\mathcal{P}(v^{\prime\prime},v')\right)\mathrm{d}v^{\prime\prime}\qquad\qquad&\\
	=\frac{\partial^{|\alpha|+|\alpha'|}}{\partial v^{\alpha}\partial (v')^{\alpha'}}\left(\mathscr{K}[F,G](v,v')\mathcal{P}(v,v')\right).&
	\end{split}
	\label{eq:2.49feb24}
\end{equation}
This way, we complete the proof of \eqref{eq:2.27June}.

The formula for general 
$Q_{r,x_{0}}(f)(v,v')$ is given by \cite[(7.2.16)]{MM07} inductively when $f$ is $\mathscr{C}^r$ near $x_0$, more precisely, we have
\begin{equation}
Q_{r,x_{0}}(f)=\sum_{\ell_1+\ell_2+|\alpha|=r}\mathscr{K}[\mathcal{J}_{\ell_1, x_0}, \frac{\partial^\alpha f_{x_0}}{\partial v^\alpha}(0)\frac{v^\alpha}{\alpha!}\mathcal{J}_{\ell_2, x_0}],
\label{eq:2.52newMar}
\end{equation}
where $\mathcal{J}_{\ell, x_0}$ is the polynomial in \eqref{eq:5.1.19DLM} of degree $\leq 3\ell$.
Thus we conclude directly that $Q_{r,x_{0}}(f)$ is a polynomial whose coefficients are given in terms of the derivatives of $f$ at $x_0$ up to order $r$, by induction on $r$, we get $Q_{\ell,x_{0}}(f)\equiv 0$ for $\ell\leq r-1$ if $f$ vanishes at $x_0$ up to order $r$. If 
$x_{0}\not\in \esupp f$, we can modify $f$ to an 
equivalent function which is identically $0$ near $x_{0}$ so that 
they define the same Toeplitz operator, in particular, from again
\cite[(7.2.16)]{MM07}, we conclude $Q_{r,x_{0}}(f)\equiv 0$.
The proof is completed.
\end{proof}

\begin{Remark}\label{rm:2.3July}
	In fact, as in \cite[Lemma 3.1]{BMMP14} (cf. Theorem \ref{thm:offdiagonal}), from 
	\eqref{eq:2.45July}, we can improve the estimate 
	\eqref{eq:June2.38on} to the following result: for any compact 
	subset $K\subset X$, any $m\in \N$ and $\delta'>0$, there exist constants 
	$C,c>0$ such that for $x,x'\in K$, $\mathrm{dist}(x,x')\leq 
	\delta'$ and $p\geq 1$, we have
	\begin{equation}
		|T_{f,p}(x,x')|_{\mathscr{C}^{m}}\leq 
		Cp^{n+m/2}e^{-c\sqrt{p}\,\mathrm{dist}(x,x')}+\mathcal{O}(p^{-\infty}).
		\label{eq:2.49July2}
	\end{equation}
\end{Remark}

As in \cite[\S IV]{BMMP14}, instead of assuming $f$ to be smooth on $U$, we can consider an assumption of lower regularity. We have the following result, which is essentially a local version of \cite[Lemma 4.2]{BMMP14}.
\begin{Theorem}\label{thm:2.5Nov}
Assume the geometric setting as given in Condition \ref{condt:main}.
	Given $f\in \mathcal{L}^{\infty}_{\mathrm{const}}(X,\C)$. For $N\in\N$. Let $U\subset X$ be an open 
	subset such that that $f$ is $\mathscr{C}^{N+1}$ on $U$. For any compact subset 
	$K\subset U$, there exists $\delta_{K}>0$ such that for any 
	$m\in \{0,1,\cdots,N\}$, $v,v'\in T_{x_{0}}X$, $\Vert v\Vert, \Vert v'\Vert 
	\leq \delta_{K}$, there exist 
$C_{K,m,N}>0$, 
	$M_{K,m,N}\in\N$ such that
	\begin{equation}
		\begin{split}
			&\bigg|\frac{1}{p^{n}}T_{f,p}(\exp_{x_{0}}(v),\exp_{x_{0}}(v'))\\
			&\qquad\qquad\qquad-\sum_{r=0}^{m} 
			\left(Q_{r,x_{0}}(f)\mathcal{P}_{x_{0}}\right)(\sqrt{p}v,\sqrt{p}v')\kappa^{-1/2}(v)\kappa^{-1/2}(v')p^{-r/2}\bigg|_{\mathscr{C}^{0}(K)}\\
		&\qquad\leq 
		C_{K,m,N}p^{-(m+1)/2}(1+\sqrt{p}\Vert v\Vert+\sqrt{p}\Vert 
		v'\Vert)^{M_{K,m,N}}\exp(-C'\sqrt{p}\Vert 
		v-v'\Vert)+\mathcal{O}(p^{-\infty}),
		\end{split}
		\label{eq:2.27June2nd}
	\end{equation}
	where the polynomials $Q_{r,x_0}(f)$ ($r=0,\cdots, N$) are the same as given in Theorem \ref{thm:2.2June} but their coefficients are of class $\mathscr{C}^{N+1-r}$ in variable $x_0\in U$.
\end{Theorem}
\begin{proof}
We proceed as in the proof of Theorem \ref{thm:2.2June}, notice that \eqref{eq:2.45July} still holds uniformly for $x_0\in K$, this way, we only need the local computation to conclude \eqref{eq:2.27June2nd}, and the local computation is the same as in \cite[Proof of Lemma 4.2]{BMMP14}.
\end{proof}

As a consequence of Theorem \ref{thm:2.5Nov}, analogous to 
\cite[Lemma 7.4.2]{MM07}, \cite[Theorem 5.1 and Remark 5.7]{BMMP14} and employing the same arguments in their proof, we also have the following result:
\begin{Corollary}\label{cor:2.6norm}
Assume the geometric setting as given in Condition \ref{condt:main}.
	For $f\in \mathcal{L}^{\infty}_{\mathrm{c}}(X,\C)$, if there 
	exists a point $x_{0}\in X$ such that $|f(x_{0})|=\Vert 
	f\Vert_{\mathcal{L}^{\infty}}$ and $f$ is $\mathscr{C}^1$ near $x_{0}$, 
	then 
	\begin{equation}
		\lim_{p\rightarrow +\infty}\Vert 
		T_{f,p}\Vert_{\mathrm{op}}=\Vert 
		f\Vert_{\mathcal{L}^{\infty}}.
		\label{eq:2.54ap24}
	\end{equation}
	More precisely, there exists a constant $C>0$ such that for $p\gg 0$,
		\begin{equation}
		\Vert 
		T_{f,p}\Vert_{\mathrm{op}}\geq \Vert 
		f\Vert_{\mathcal{L}^{\infty}}-\frac{C}{\sqrt{p}}.
	\end{equation}
\end{Corollary}
\begin{Remark}\label{lem:2.7norm}
 	In fact, the result \eqref{eq:2.54ap24} still holds true if we only assume $f$ 
 	to be continuous at its (essential) maximum point $x_{0}$, which was proved in 
 	\cite[Theorem 5.1]{BMMP14} for the case of compact manifolds.
\end{Remark}

The following result is an extension of \cite[Theorem 3.7]{BMMP14},
where it is proved for compact case, and the same proof applies
in our setting since the (essential) support of $f$ is assumed to be compact.
\begin{Theorem}[cf. {\cite[Theorem 3.7]{BMMP14}}]\label{thm:2.8jan}
	For $f\in \mathcal{L}^{\infty}_{\mathrm{c}}(X,\C)$, then $T_{f,p}$ is a trace class. We have the 
	following expansion as $p\rightarrow +\infty$: for any $N\in \N$,
	\begin{equation}
		\mathrm{Tr}[T_{f,p}]=\sum_{r=0}^N \bb{t}_{r,f}p^{n-r}+\mathcal{O}(p^{n-N-1}),
	\end{equation}
	where $\bb{t}_{r,f}=\int_X f(x)\bb{b}_{r}(x) \mathrm{dV}(x)$, the functions $\bb{b}_{r}$ ($r=0, 1, \ldots $) are given in \eqref{eq:2.1.4June}.
\end{Theorem}

Let us discuss a bit more on the asymptotics of the spectrum of $T_{f,p}$.
We restrict us to the case 
$f\in\mathscr{L}^\infty_{\mathrm{c}}(X,\R_{\geq 0})$. 
Then $T_{f,p}$ is a self-adjoint compact operator on 
$H^{0}_{(2)}(X,L^p\otimes E)$. When $d_{p}=\infty$, 
the residual spectrum of $T_{f,p}$ contains only $0$,
and each nonzero eigenvalue in the point 
spectrum of $T_{f,p}$ always has finite multiplicity.
 When $f$ is given as the indicator function of a (Borel) 
 subset of $X$, the asymptotic statistics of the eigenvalues 
 of $T_{f,p}$ were studied by 
 \cite{MR2016088,MR4093864,MR3433287,MR1828799} 
 in various settings for compact or noncompact $X$. 
 Moreover, for a continuous function $f\geq 0$ with compact support, 
 the spectral densitiy measure 
$\mu_{f,p}$ of $T_{f,p}$ is defined as the sum of the Dirac masses
at all the eigenvalues (counted with multiplicities) of $T_{f,p}$,
which is locally finite. 
A result of \cite{ILMM24} shows that as $p\rightarrow +\infty$,
we have the weak convergence of measures on $(0, |f|_{\mathscr{C}^0}]$,
\begin{equation}
p^{-n}\mu_{f,p}\rightarrow f_*\left(\frac{1}{n!}c_1(L,h_L)^n\right).
\end{equation}
This extends the results for compact K\"{a}hler manifolds or domains in $\C^n$,
such as \cite{MR2016088, MR1828799}.

In Section \ref{ss:compactandeigen}, we will give more results
for the asymptotics of the lowest eigenvalues of $T_{f,p}$
on a compact Hermitian manifold.

\subsection{Compositions of Toeplitz operators}\label{ss5.1a}
Now we consider the composition $T_{f,p}\circ T_{g,p}$ (or simply $T_{f,p}T_{g,p}$) of Toeplitz 
operators for two 
functions $f,g\in \mathcal{L}^{\infty}_{\mathrm{const}}(X,\C)$. Based at 
the previous Section, the 
following theorem is an easy extension of some results presented by 
Ma and Marinescu in \cite[Theorems 7.4.1 \& 
7.5.1]{MM07} and \cite[Theorem 0.2]{MM12}.

\begin{Theorem}\label{thm:2.4June}
Assume the geometric setting as in Condition \ref{condt:main}.
	Given $f, g\in \mathcal{L}^{\infty}_{\mathrm{const}}(X,\C)$. Then for any compact subset $K\subset X$, there exists 
$C_{K,m}>0$ such that for $x\in K$ and for all $p\geq 1$, we have 
the on-diagonal 
estimate
\begin{equation}
	\left|(T_{f,p}T_{g,p})(x,x)\right|_{\mathscr{C}^{m}}\leq 
	C_{K,m} p^{n+\frac{m}{2}}.
	\label{eq:June2.44on}
\end{equation}
	
For any $m,\ell\in\N$, $\varepsilon >0$, 
and a compact subset $K\subset X$, there exists 
$C_{K,m,\ell,\varepsilon}>0$ such that for $x,x'\in K$ with 
$\mathrm{dist}(x,x')\geq \varepsilon$ and for all $p\geq 1$, we have the off-diagonal 
estimate
\begin{equation}
	\left|(T_{f,p}T_{g,p})(x,x')\right|_{\mathscr{C}^{m}}\leq 
	C_{K,m,\ell,\varepsilon} p^{-\ell}.
	\label{eq:July2.38}
\end{equation}

	Let $U\subset X$ be an open 
	subset such that that both $f$ and $g$ are smooth on $U$. There 
	exists a family of polynomials $\{Q_{r,x_{0}}(f,g)\in \C[v,v'], 
	v,v'\in T_{x_{0}}X\}_{x_{0}\in U}$ with the same parity as $r$, 
	and smooth in $x_{0}\in U$ such that for any compact subset 
	$K\subset U$, there exist $\delta_{K}>0$, $C'>0$ such that for any 
	$\ell,N,m\in\N$, $v,v'\in T_{x_{0}}X$, $\Vert v\Vert, \Vert v'\Vert 
	\leq \delta_{K}$ and multi-indices $\alpha$, $\alpha'\in 
\mathbb{N}^{2n}$ with $|\alpha|+|\alpha'|\leq m$, there exist $C=C_{K,\ell,m,N}>0$, 
	$M_{K,\ell,m,N}\in\N$ such that
	\begin{equation}
		\begin{split}
			\bigg|\frac{\partial^{|\alpha|+|\alpha'|}}{\partial 
			v^{\alpha}\partial (v')^{\alpha'}}\bigg(\frac{1}{p^{n}}(T_{f,p}T_{g,p})(\exp_{x_{0}}(v),\exp_{x_{0}}(v'))\qquad\qquad\qquad\qquad\qquad\qquad\qquad\qquad\qquad\qquad\\
			-\sum_{r=0}^{N} 
			\left(Q_{r,x_{0}}(f,g)\mathcal{P}_{x_{0}}\right)(\sqrt{p}v,\sqrt{p}v')\kappa^{-1/2}(v)\kappa^{-1/2}(v')p^{-r/2}\bigg)\bigg|_{\mathscr{C}^{\ell}(K)}&\\
		\leq 
		Cp^{-(N+1-m)/2}(1+\sqrt{p}\Vert v\Vert+\sqrt{p}\Vert 
		v'\Vert)^{M_{K,\ell,m,N}}\exp(-C'\sqrt{p}\Vert 
		v-v'\Vert)+\mathcal{O}(p^{-\infty}).&
		\end{split}
		\label{eq:2.43June}
	\end{equation}
	In particular, we have $Q_{0,x_{0}}(f,g)=f(x_{0})g(x_{0})$. If 
	$x_{0}\in U$ and
	$x_{0}\not\in \esupp f \cap \esupp g$, then we have 
	$Q_{r,x_{0}}(f,g)\equiv 0$ for $r=0,1,\cdots$.
\end{Theorem}

\begin{proof}
	The proof of this theorem follows from the analogous arguments in 
	the proof of \cite[Theorem 7.5.1]{MM07} and of Theorem 
	\ref{thm:2.2June} in previous Section. In particular, for 
	$x_{0}\in X$, for $v,v'\in T_{x_{0}}X$ with $\Vert v\Vert$, 
	$\Vert v'\Vert \leq \delta$ ($\delta>0$ is chosen properly), then 
	we have (cf. \cite[(7.4.5)]{MM07})
	\begin{equation}
		\begin{split}
				&(T_{f,p}T_{g,p})(\exp_{x_{0}}(v),\exp_{x_{0}}(v'))\\
				&\qquad=\int_{v^{\prime\prime}\in T_{x_{0}}X} 
				T_{f,p}(\exp_{x_{0}}(v),\exp_{x_{0}}(v^{\prime\prime}))\rho(\Vert v^{\prime\prime}\Vert/\delta)T_{f,p}(\exp_{x_{0}}(v^{\prime\prime}),\exp_{x_{0}}(v^{\prime}))\kappa_{x_{0}}(v^{\prime\prime})\mathrm{d}v^{\prime\prime}\\
				&\qquad\qquad\qquad\qquad+\mathcal{O}(p^{-\infty}).
		\end{split}
		\label{eq:2.53July2}
	\end{equation}
	Then combining \eqref{eq:2.53July2} with \eqref{eq:June2.38on} 
	and \eqref{eq:2.49July2}, we conclude the estimate 
	\eqref{eq:June2.44on}. Then similarly, by \eqref{eq:June2.38}, we 
	get \eqref{eq:July2.38}. When $f$, $g$ are smooth near $x_{0}$, the 
	expansion \eqref{eq:2.43June} follows from the combination of the expansion \eqref{eq:2.27June}, \eqref{eq:2.49feb24} and \eqref{eq:2.53July2} by the arguments given in \cite[Section 
7.4]{MM07}. Note that an explicit formula for $Q_{r,x_{0}}(f,g)$ is 
given in \cite[(7.4.7)]{MM07}(see also \eqref{eq:2.53march}). The proof is complete.
\end{proof}

\begin{Lemma}\label{Lm:2.11new}
The polynomial $Q_{r,x_0}(f,g)(v,v')$ in $v,v'$ has degree $\leq\,3r$, and its coefficients are given by polynomials in terms of the derivatives of $f$, $g$ at $x_0$ up to order $r$.
\end{Lemma}
\begin{proof}
We use the notation in \cite[Lemma 7.1.1]{MM07} (see also \eqref{eq:2.49march24}), we have
\begin{equation}
Q_{r,x_{0}}(f,g)=\sum_{\ell_1+\ell_2=r}\mathscr{K}[Q_{\ell_1,x_0}(f),Q_{\ell_2,x_0}(g)].
\label{eq:2.53march}
\end{equation}
By \eqref{eq:2.49march24} and \eqref{eq:2.52newMar}, we have $\deg Q_{\ell_1,x_0}(f) \leq 3\ell_1$, $\deg Q_{\ell_2,x_0}(g) \leq 3\ell_2$, hence the degree of $Q_{r,x_{0}}(f,g)$ (in $v,v'$) $\leq 3r$. The rest part for the coefficients is also deduced from \eqref{eq:2.52newMar}. The proof is completed.
\end{proof}

\begin{Remark}
	Analogous to Remark \ref{rm:2.3July}, we can improve the estimate 
	\eqref{eq:June2.44on} to the following result: for any compact 
	subset $K\subset X$, any $m\in \N$ and $\delta'>0$, there exist constants 
	$C,c>0$ such that for $x,x'\in K$, $\mathrm{dist}(x,x')\leq 
	\delta'$ and $p\geq 1$, we have
	\begin{equation}
		|(T_{f,p}T_{g,p})(x,x')|_{\mathscr{C}^{m}}\leq 
		Cp^{n+m/2}e^{-c\sqrt{p}\,\mathrm{dist}(x,x')}+\mathcal{O}(p^{-\infty}).
		\label{eq:2.55July2}
	\end{equation}
\end{Remark}

As an analog of Theorem \ref{thm:2.5Nov}, combining the localized computation \eqref{eq:2.53July2} with the arguments in \cite[Proof of Theorem 4.1]{BMMP14}, we conclude the following near-diagonal expansions for the kernel of $T_{f,p}T_{g,p}$.
\begin{Theorem}\label{thm:2.5Novfg}
	Assume the geometric setting as given in Condition \ref{condt:main}. Given $N\in\N$ and $f, g\in \mathcal{L}^{\infty}_{\mathrm{const}}(X,\C)$.
	Let $U\subset X$ be an open 
	subset such that that both $f|_U,\; g|_U\in\mathscr{C}^{N+1}(U)$. Then for any compact subset 
	$K\subset U$, there exists $\delta_{K}>0$ such that for any 
	$m\in\{0,1,\cdots,N\}$, $v,v'\in T_{x_{0}}X$, $\Vert v\Vert, \Vert v'\Vert 
	\leq \delta_{K}$, there exist $C_{K,m,N}>0$, 
	$M_{K,m,N}\in\N$ such that
	\begin{equation}
		\begin{split}
			\bigg|\frac{1}{p^{n}}(T_{f,p}T_{g,p})(\exp_{x}(v),\exp_{x}(v'))\qquad\qquad\qquad\qquad\qquad\qquad\qquad\qquad\qquad\qquad\\
			-\sum_{r=0}^{m} 
			\left(Q_{r,x}(f,g)\mathcal{P}_{x}\right)(\sqrt{p}v,\sqrt{p}v')\kappa^{-1/2}(v)\kappa^{-1/2}(v')p^{-r/2}\bigg|_{\mathscr{C}^{0}(K)}&\\
		\leq 
		C_{K,m,N}p^{-(m+1)/2}(1+\sqrt{p}\Vert v\Vert+\sqrt{p}\Vert 
		v'\Vert)^{M_{K,m,N}}\exp(-C'\sqrt{p}\Vert 
		v-v'\Vert)+\mathcal{O}(p^{-\infty}),&
		\end{split}
		\label{eq:2.43Nov}
	\end{equation}
	where the coefficients $Q_{r,x}(f,g)$ are the polynomials given in Theorem \ref{thm:2.4June}. 
\end{Theorem}

\subsection{Normalized Berezin--Toeplitz kernels;
proof of Theorem \ref{thm:5.1.1}}\label{ss2.4:BT}
In this Section, we give the proof of Theorem \ref{thm:5.1.1}.
Due to the expansion presented in Theorem \ref{thm:2.4June}, the 
	proof of this theorem is an easy modification of the proofs of 
	\cite[Theorems 1.8 and 5.1]{DLM:21}. We include the details as 
	follows. 
\begin{proof}[Proof of Theorem \ref{thm:5.1.1}]

	Note that since $f$ is smooth on an open neighbourhood of 
$\overline{U}$, so that we can apply the asymptotic expansion 
\eqref{eq:2.43June} for all points $x,y\in U$. In particular, we have 
the uniform expansion on $\overline{U}$,
\begin{equation}
	T_{f,p}^{2}(x,x)=p^{n}f(x)^{2}\bb{b}_{0}(x)+\mathcal{O}(p^{n-1}),
	\label{eq:2.72Nov}
\end{equation}
where $\bb{b}_{0}(x)$ is as given in \eqref{eq:2.1.4June}. Since 
we assume $f$ to be nonvanishing on $\overline{U}$, then there exists 
$c_{U}>0$ such that for all sufficiently large $p>0$, we have
\begin{equation}
	T_{f,p}^{2}(x,x)\geq c_{U}p^{n}.
\end{equation}

	We start by proving the second estimate of the 
theorem. One way to see this estimate is from \eqref{eq:2.55July2}, 
where the constant $M_{k}$ is determined by the constant $c$ that appears in 
the exponential term. In the sequel, we prove it
by the arguments in \cite[Section 2.3]{DLM:21}.

Note that $U$ is relatively compact in $X$, so 
$\overline{U}$ is compact and all the results of Theorem 
\ref{thm:2.4June} are 
applicable for the points in $U$. Let $\delta_{K}>0$ be the sufficiently small quantity 
stated in the last part of Theorem 
\ref{thm:2.4June} with $K=\overline{U}$. Then by 
\eqref{eq:July2.38}, if 
$x,y\in U$ is such that $\operatorname{dist}(x,y)\geq \delta_{K}$, we have
\begin{equation}
	|T_{f,p}^{2}(x,y)|_{h_{p,x}\otimes h^{*}_{p,y}}\leq C_{U,0,k,\delta_{K}}\, p^{-k}.
	\label{eq:3.1.1}
\end{equation}

For the given $k$, we will determine a constant $M_{k}$ later on.
We fix a large enough $p_{0}\in\mathbb{N}$ such that
\begin{equation}
	b\,\sqrt{\frac{\log{p_{0}}}{p_{0}}}\leq \frac{\delta_{K}}{2}.
	\label{eq:3.1.2}
\end{equation}
For $p>p_{0}$, if $x,y\in U$ is such that 
$b\sqrt{\frac{\log{p}}{p}}\leq\operatorname{dist}(x,y)< \delta_{K}$, then 
we take advantage of the expansion in \eqref{eq:2.43June} with $N=2k$, $x_{0}=x$, 
$v=0$, 
$y=\exp_{x}(v')$, and $v'\in T_{x}\Sigma$, in order to obtain
\begin{equation}
		\begin{split}
			\bigg|\frac{1}{p^{n}}T^{2}_{f,p}(x,y)-\sum_{r=0}^{2k} 
			\left(Q_{r,x}(f,f)\mathcal{P}_{x}\right)(0,\sqrt{p}v')\kappa^{-1/2}(v')p^{-r/2}\bigg|&\\
		\leq 
		Cp^{-(k+1/2)}(1+\sqrt{p}\Vert 
		v'\Vert)^{M_{K,0,0,2k}}\exp(-C'\sqrt{p}\Vert 
		v'\Vert)+\mathcal{O}(p^{-\infty}).&
		\end{split}
		\label{eq:2.60June}
	\end{equation}

Now for $k \ge 1$, by Lemma \ref{Lm:2.11new}, we have 
\begin{equation}
	\ell(r):=\max\{\deg_{v'}Q_{r,x}(f,f)(0,v'), x\in U\} \leq 3r.
	\label{eq:2.61July}
\end{equation}

Note that $\Vert v'\Vert =\operatorname{dist}(x,y)$. By 
\eqref{eq:2.17July}, \eqref{eq:5.1.15paris}, \eqref{eq:2.61July},
 and the properties of 
$Q_{r,x}(f,f)$ in Theorem \ref{thm:2.4June}, together with the fact 
that $\delta_{K}>\Vert v'\Vert \geq 
b\sqrt{\frac{\log{p}}{p}}$, for $r=0,1,\cdots, 2k$, we get that
\begin{equation}
	|\left(Q_{r,x}(f,f)\mathcal{P}_{x}\right)(0,\sqrt{p}v')\kappa^{-1/2}(v')|\leq C_{K,r} 
	p^{\ell(r)/2}\exp \Big (-\frac{\varepsilon_{0}}{4} 
	b^{2}\log{p} \Big ),
\end{equation}
where the constant $C_{K,r}>0$ does not depend on $x\in U$. If
\begin{equation}
	\frac{\varepsilon_{0}}{4}b^{2}\geq 
	3k,
\end{equation}
then we have for $r=0,\ldots, 2k$,
\begin{equation}
	|\left(Q_{r,x}(f,f)\mathcal{P}_{x}\right)(0,\sqrt{p}v')\kappa^{-1/2}(v')p^{n-r/2}|\leq 
	C_{K,r}p^{n-k}.
	\label{eq:2.64July}
\end{equation}
We may take $M_{k}= 12k$ in 
our constraint for $b$.
Finally, combining \eqref{eq:2.72Nov}--\eqref{eq:2.64July}, we get the 
desired estimate for any $p>p_{0}$.

We next prove the first part of our theorem. For 
this purpose, we only need to consider sufficiently large $p$ such that 
$b\sqrt{\frac{\log{p}}{p}}\leq 
\frac{\delta_{K}}{2}$, or $p\geq p_{0}$.

Recall that $\kappa(R^L, U)>0$ is defined in \eqref{eq:1.09intro}. We set
\begin{equation}
m(U,b):=\left\lceil \frac{b^2}{2} \kappa(R^L, U)\right\rceil.
\end{equation}
By \eqref{eq:2.21parisnew}, for $\Vert v'\Vert =|z'|\leq 
b\sqrt{\frac{\log{p}}{p}}$, we have
\begin{equation}
\exp \Big(\frac{p}{4}\Phi^{TX}_x(0,v')^{2}\Big ) \leq p^{m(U,b)/2}.
\end{equation}

In the expansion \eqref{eq:2.43June}, we take $x_{0}=x, y=\exp_{x}(v'), 
N=m(U,b)$, so $\operatorname{dist}(x,y)=\Vert v'\Vert =|z'|\leq 
b\sqrt{\frac{\log{p}}{p}}$, where $z'\in\C$ is the complex coordinate 
for $v'$. We infer
\begin{equation}
		\begin{split}
			\bigg|\frac{1}{p^{n}}T^{2}_{f,p}(x,y)-\sum_{r=0}^{m(U,b)} 
			\left(Q_{r,x}(f,f)\mathcal{P}_{x}\right)(0,\sqrt{p}v')\kappa^{-1/2}(v')p^{-r/2}\bigg|
		\leq 
		Cp^{-(m(U,b)+1)/2}.
		\end{split}
		\label{eq:3.1.7}
	\end{equation}
Since $\Vert v'\Vert \leq b\sqrt{\frac{\log{p}}{p}}$, by \eqref{eq:2.61July} we infer that
$|Q_{r,x}(f,f)(0,\sqrt{p}v')|\leq C_{f,U}|\log{p}|^{\ell(r)/2}$ for some constant $C_{f,U}>0$. Note that $|\log{p}|^{\ell(r)/2} p^{-r/2}=\mathcal{O}(p^{-1/2+\epsilon})$ for $r\geq 1$ as $p$ grows.

Note that $f$ is a real function which does not vanish near $\overline{U}$. The expansion \eqref{eq:2.72Nov} in combination with \eqref{eq:3.1.7} then supplies us with
\begin{equation}
	\begin{split}
		\frac{\exp \Big(\frac{p}{4}\Phi^{TX}_x(0,v')^{2}\Big )T^2_{f,p}(x,y)}{\sqrt{T_{f,p}^2(x,x)}\sqrt{T_{f,p}^2(y,y)}}&=\frac{\bb{b}_0(x)f(x)^2\kappa^{-1/2}(v')+\mathcal{O}(p^{-1/2+\epsilon})}{\sqrt{f(x)^2\bb{b}_0(x)+\mathcal{O}(p^{-1})}\sqrt{f(y)^2\bb{b}_0(y)+\mathcal{O}(p^{-1})}}\\
		&=1+\mathcal{O}(\Vert v'\Vert +p^{-1/2+\epsilon})\\
		&=1+o(1), \text{ as }p\rightarrow +\infty.
	\end{split}
	\label{eq:3.1.8bis}
\end{equation}
In \eqref{eq:3.1.8bis}, the small term $o(1)$ in the last line represents our function $R_{p,x}(v')$, then \eqref{eq:2.62feb24} follows clearly.

Since in the asymptotic expansion \eqref{eq:3.1.7} we 
have trivialized the line bundle near $x$ using the Chern 
connections, we have
\begin{equation}
	|T_{f,p}^2(x,y)|_{h_{p,x}\otimes 
	h^\ast_{p,y}}=|T_{f,p}^2(x,y)|.
	\label{eq:3.1.9bis}
\end{equation}
Combining \eqref{eq:3.1.8bis} and \eqref{eq:3.1.9bis}, we get the first part of \eqref{eq:5.2.4}. This completes the proof of Theorem 
\ref{thm:5.1.1}.	
\end{proof}

Now, when $f$ has a lower regularity such as $\mathscr{C}^{m+1}$ with $m=m(U,b)$, the above arguments still hold with the bounded choices of $k$. We give the proof of Corollary \ref{cor:5.1.1} as follows.
\begin{proof}[Proof of Corollary \ref{cor:5.1.1}]

The inequality \eqref{eq:3.1.1} does not require a specific regularity on $f$ due to \eqref{eq:July2.38}. Then by Theorem \ref{thm:2.5Novfg}, when $f|_U$ is of $\mathscr{C}^{m(U)+1}(U)$, the arguments \eqref{eq:2.60June} -- \eqref{eq:2.64July} holds true with $k=n+1$ and $b=\sqrt{\frac{12n+12}{\varepsilon_0}}$, so that the second part of \eqref{eq:5.2.4Nov} holds. The first part of \eqref{eq:5.2.4Nov} follows from the same arguments as in last step of the proof of Theorem \ref{thm:5.1.1}, i.e., \eqref{eq:3.1.7} -- \eqref{eq:3.1.9bis}. This way, we complete our proof.
\end{proof}

Note that the upper bound \eqref{eq:2.62feb24} is not optimal for $R_{p,x}(v')$, since we have $R_{p,x}(0)=0$ and $\nabla R_{p,x}(0)=0$ (here $\nabla$ denotes the coordinate derivatives in $v'$).
The following estimate is an analog of \cite[Proposition 2.8]{MR2465693} in our Berezin--Toeplitz setting.
\begin{Lemma}\label{lm:2.14feb24}
With the same assumptions in Theorem \ref{thm:5.1.1} (in particular, $f$ is smooth near $U$), the term $R_{p,x}(v')$ satisfies the following estimate: there exists $C_1=C_1(f,\epsilon, U)>0$ such that for all sufficiently large $p$, $x\in U$, $v'\in T_x X$ with $\|v'\|\leq b\sqrt{\log{p}}$,
\begin{equation}
|R_{p,x}(v'/\sqrt{p})|\leq C_1\|v'\|^2 p^{-1/2+\epsilon}
\label{eq:2.75feb24}
\end{equation}

For given $k, \ell\in \N$, there exist a sufficiently large $b>0$ such that there exist a constant $C_2>0$ such that for all $x,y\in U$, $\mathrm{dist}(x,y)\geq b\sqrt{\log{p}/p}$, we have
\begin{equation}
|\nabla^\ell_{x,y}N_{f,p}(x,y)|\leq C_2 p^{-k}.
\label{eq:2.76feb24}
\end{equation}
\end{Lemma}
\begin{proof}
Since $f$ is smooth on a neighborhood of $U$, then the expansion \eqref{eq:2.43June} holds on $\overline{U}$ with $\mathscr{C}^\ell$-norm, this way, \eqref{eq:2.76feb24} follows by repeating the same arguments in the first part of the proof of Theorem \ref{thm:5.1.1}. 

By \eqref{eq:3.1.7}, we conclude that for all $x\in U, v'\in T_xX$, $\|v'\|\leq \delta_K$ (where $K=\overline{U}$), we have
\begin{equation}
|\nabla^2 R_{p,x}(v')|\leq Cp^{1/2}\left(1+(\sqrt{p}\|v'\|)^{M}\right),
\end{equation}
for some large integer $M\in\N$.

Set $H_{p,x}(v')=R_{p,x}(v'/\sqrt{p})$. Then for $\|v'\|\leq b\sqrt{\log{p}}$,
\begin{equation}
|\nabla^2 H_{p,x}(v')|\leq C_{\epsilon,U} p^{-1/2+\epsilon}.
\end{equation}
Then 
\begin{equation}
|\nabla H_{p,x}(v')|\leq \sup_{t\in [0,1]} |(\nabla^2 H_{p,x})(tv')|\cdot\|v'\|\leq  C_{\epsilon,U} p^{-1/2+\epsilon} \|v'\|.
\label{eq:2.78feb24}
\end{equation}
Then \eqref{eq:2.75feb24} follows after a similar computation as in \eqref{eq:2.78feb24}.
\end{proof}

\section{Gaussian $\mathcal{L}^2$-holomorphic sections via Toeplitz 
operators}\label{section2}
In this section we recall the Gaussian $\cLL$-holomorphic 
section of a Hermitian line bundle on a complex manifold defined 
through a given Hilbert-Schmidt Toeplitz operator, this kind of 
random holomorphic sections were constructed in \cite[Section 4]{DrLM:2023aa}, taking
advantage of the theory of abstract Wiener spaces.

\subsection{Gaussian $\mathcal{L}^{2}$-holomorphic 
sections}\label{ss:4.3b}

We now recall the construction of Gaussian $\mathcal{L}^{2}$-holomorphic 
sections given in \cite[Sections 4.3 and 4.4]{DrLM:2023aa}. Let $(F,h_F)$ be a holomoprhic line bundle on $X$ with smooth Hermitian metric $h_F$. Furthermore, by $H^{0}_{(2)}(X,F)$ we denote the space of $\mathcal{L}^{2}$-holomorphic sections of $F$ on $X$, which is a separable Hilbert space with the Hilbert metric given by the $\mathcal{L}^{2}$-inner product induced by $h_F$ and $\mathrm{dV}$. Set
\begin{equation}
d_F:=\dim_\C H^{0}_{(2)}(X,F)\in \N\cup\{\infty\}.
\end{equation}
In this Section, we always assume that $d_F>0$ and by $P: \cLL(X,F)\rightarrow H^{0}_{(2)}(X,F)$ we denote the Bergman projector.

Let $\mathcal{L}^\infty_{\mathrm{c}}(X,\R)$
be the space of measurable essentially bounded real functions on $X$
with compact essential support.
For $f\in\mathcal{L}^\infty_{\mathrm{c}}(X,\R)$
the Toeplitz operator 
$T_f=PM_f: H^{0}_{(2)}(X,F)\rightarrow H^{0}_{(2)}(X,F)$
is Hilbert-Schmidt and self-adjoint.
We will consider in the sequel only non-trivial functions 
$f\neq0\in\mathcal{L}^\infty_{\mathrm{c}}(X,\R)$.

When $d_F<\infty$, for $f\geq 0$ with compact support, 
$T_f$ is a linear isomorphism of $H^{0}_{(2)}(X,F)$
with strictly positive eigenvalues. When $d_F=\infty$, 
the residual spectrum of $T_f$ contains only $0$, 
and each nonzero eigenvalue in the point 
spectrum of $T_{f}$ always has finite multiplicity. 
In particular, the nonzero eigenvalues of $T^2_{f}$ form
a decreasing sequence of strictly positive real numbers, 
\begin{equation}
	\lambda^2_{1}\geq \lambda^2_{2}\geq \ldots\geq 
	\lambda^2_{m}\geq \ldots\rightarrow 0,
	\label{eq:2.1.9}
\end{equation}
where $\lambda_1,\lambda_2,\ldots $ denote the 
nonzero values in the point spectrum of $T_{f}$, 
repeated according to their multiplicities. If furthermore $f\geq 0$, 
the eigenvalues $\lambda_1,\lambda_2,\ldots $ are strictly positive.

Next, we introduce
	\begin{equation}
		H^{0}_{(2)}(X,L,f):=(\ker T_{f})^{\perp}=\overline{T_{f}H^{0}_{(2)}(X,L)}\subset 
		H^{0}_{(2)}(X,L),
	\end{equation}
which is a Hilbert space, and the sections in $H^{0}_{(2)}(X,L,f)$ 
are the $\cLL$-holomorphic 
sections of $L$ \textit{detected} by $f$. We consider the (self-adjoint)
Hilbert-Schmidt operator
	\begin{equation}
		T_{f}^{\sharp}:=T_{f}|_{H^{0}_{(2)}(X,L,f)}: 
		H^{0}_{(2)}(X,L,f)\rightarrow 
		H^{0}_{(2)}(X,L,f).
	\end{equation}
Note that for a notrivial $f\geq 0$ with compact support, 
$T_f$ is injective and $H^{0}_{(2)}(X,L,f)=H^{0}_{(2)}(X,L)$.
We furthermore set $d_F(f):=\dim_\C H^{0}_{(2)}(X,L,f)\in \N\cup \{\infty\}$,
and it is immediate that $d_F(f)\leq d_F$. 

The following lemma now is elementary.
\begin{Lemma}\label{lm:3.1}
For $f\in\mathcal{L}^\infty_{\mathrm{c}}(X,\R)\setminus\{0\}$
we have
\begin{equation}
d_F(f)\geq \frac{1}{|f|^2_{\mathcal{L}^\infty(X)}}
\int_X T^2_f(x,x)\mathrm{dV}(x).
\end{equation}
\end{Lemma}

Suppose that $T_f\neq 0$ (that is, $d_F(f)\geq 1$). In this case, we can choose
 an orthonormal basis $\{S_{j}\}^{d_F(f)}_{j=1}$  of 
$H^{0}_{(2)}(X,F,f)$ with respect to the $\cLL$-metric 
 such that
\begin{equation}
	T_{f}S_{j}=\lambda_{j}S_{j}.
	\label{eq:2.1.10}
\end{equation}

On $H^{0}_{(2)}(X,L,f)$ with $d_F(f)\geq 1$, the operator $T_{f}^{\sharp}$ is one-to-one and Hilbert-Schmidt, and
$\|\cdot\|_{f}:=\|T^\sharp_{f}\cdot\|$ defines a Hermitian measurable norm on 
$H^{0}_{(2)}(X,F,f)$ due to \cite[Proposition 4.2]{DrLM:2023aa}. We furthermore denote by $\mathcal{B}_{f}(X,F)$ the 
completion of $H^{0}_{(2)}(X,F,f)$ with respect to 
$\|\cdot\|_{f}$ and set
\begin{equation}
	\ell^{2}_{f}(\C)=\Big\{(a_{j}\in\C)_{j\geq 1}\;:\;\sum_{j\geq 
	1}\lambda_{j}^{2}|a_{j}|^{2}<\infty \Big\}.
	\label{eq:2.1.16}
\end{equation}
Endowed with the induced norm this is clearly a separable Hilbert space, and using the basis as in \eqref{eq:2.1.10}, we have
\begin{equation}
	\mathcal{B}_{f}(X,F)\simeq \ell^{2}_{f}(\C).
	\label{eq:2.1.17}
\end{equation}

\begin{Proposition}\label{prop:2.10}
	Assume $0\neq f\in \mathcal{L}^\infty_{\mathrm{c}}(X,\R)$ and $T_f\neq 0$. 
	Then the operator $T_{f}$ extends uniquely to an isomorphism of Hilbert spaces 
	\begin{equation}
		\widehat{T}_{f}: \big(\mathcal{B}_{f}(X,F),\|\cdot\|_{f} \big)\rightarrow 
		\big (H^{0}_{(2)}(X,F,f),\|\cdot\|_{\cLL(X,F)} \big).
		\label{eq:5.30}
	\end{equation}
\end{Proposition}
Given $0\neq f\in \mathcal{L}^\infty_{\mathrm{c}}(X,\R_{\geq 0})$, if $d_F(f)<\infty$, 
we set
\begin{equation}
\big(\mathcal{B}_{f}(X,F),\|\cdot\|_{f} \big)= 
\big(H^{0}_{(2)}(X,F,f),\|\cdot\|_{f} \big),\quad \text{and} 
\quad \widehat{T}_{f}:=T^\sharp_{f}.
\end{equation}
This unifies the notation for both cases $d_F<\infty$ and 
$d_F=\infty$.

Let $\mathcal{B}_{f}(X,F)^{*}$ be the topological dual space of $\mathcal{B}_{f}(X,F)$. 
If $\alpha\in \mathcal{B}_{f}(X,F)^{*}$, then it is uniquely determined by the  continuous linear functional on $\alpha|_{H^{0}_{(2)}(X,F,f)}$ on $H^{0}_{(2)}(X,F,f)$. This 
way, we regard $\mathcal{B}_{f}(X,F)^{*}$ as a (dense) subspace of 
$H^{0}_{(2)}(X,F,f)^{*}$, where $H^{0}_{(2)}(X,F,f)^{*}$ can be identified with $H^{0}_{(2)}(X,F,f)$ via the $\cLL$-inner product. By a slight abuse of notation we denote by 
$\mathcal{S}_f$ the Borel $\sigma$-algebra of $\mathcal{B}_{f}(X,F)$. Then each 
$\alpha\in \mathcal{B}_{f}(X,F)^{*}$ is a Borel-measurable function from 
$\mathcal{B}_{f}(X,F)$ to $\C$.
For $V\subset \mathcal{B}_{f}(X,F)^{*}\subset H^{0}_{(2)}(X,F,f)$ an arbitrary finite 
dimensional subspace we introduce the notation
\begin{align}
\begin{split}
	\phi_{V}: \mathcal{B}_{f}(X,F)&\rightarrow V,\;
	\phi_{V}(b)=\sum_{j=1}^{\dim_{\C} V} (b,v_{j})v_{j},
\end{split}
	\label{eq:2.1.12s}
\end{align}
where $\{v_{j}\}$ is an orthonormal basis of $(V, \|\cdot\|_{\cLL(X,F)})$.
Gross \cite{Gross67} proved the following result. 
\begin{Theorem}[{L.\ Gross \cite{Gross67}, see also \cite[Theorem 4.3]{DrLM:2023aa}}]\label{thm:Gross}
	Let $\|\cdot\|_f$ denote the measurable norm on $H^{0}_{(2)}(X,F,f)$ as introduced after \eqref{eq:2.1.10}. 
	There exists a unique probability measure $\mathcal{P}_f$ on 
	$(\mathcal{B}_{f}(X,F),\mathcal{S}_f)$ such that for  $V\subset \mathcal{B}_{f}(X,F)^{*}$ any finite dimensional 
	subspace,
	\begin{equation}
		\mathcal{P}_f(\phi^{-1}_{V}(U))=\mu_{V,\;\|\cdot\|_{\cLL(X,F)}}(U),
		\label{eq:2.1.13}
	\end{equation}
	for all Borel subset  $U$ of $V$, where $\mu_{V,\;\|\cdot\|_{\cLL(X,F)}}$ denotes the standard 
	Gaussian measure on $V$ with respect to the Hermitian metric 
	associated with $\|\cdot\|_{\cLL(X,F)}$. The triple 
	$(\mathcal{B}_{f}(X,F),\mathcal{S}_f,\mathcal{P}_f)$ is called an abstract 
	Wiener space.
\end{Theorem}

\begin{Definition}\label{def:4.11}
Let $\P_{f}$ be the
Gaussian probability measure on $H^{0}_{(2)}(X,F,f)$ given 
by the pushforward of the probability measure $\mathcal{P}_{f}$ from Theorem \ref{thm:Gross} through the 
isomorphism in \eqref{eq:5.30}. This way, we randomize the sections 
in $H^{0}_{(2)}(X,F,f)$. A Gaussian (random) $\cLL$-holomorphic section $\bb{S}_f$ of $F$ associated to a nonzero $f\in \mathcal{L}^\infty_{\mathrm{c}}(X,\R)$ is a random variable valued in $H^{0}_{(2)}(X,F,f)$ with law $\P_{f}$, i.e., $\bb{S}_f\sim \P_f$. 

\end{Definition}

	When $0<d_F(f)<\infty$, then $\mathcal{B}_{f}(X,F)=H^{0}_{(2)}(X,F,f)$, 
	and $\mathcal{P}_{f}=\P_{\mathrm{st}}$ is exactly the standard Gaussian 
	probability measure on $H^{0}_{(2)}(X,F,f)$ with respect to the 
	$\cLL$-inner product; so $\P_{f}$ is the pushforward of $\P_{\mathrm{st}}$
	via the isomorphism $T_{f}$. For $0<d_F(f)<\infty$, let $\bb{S}$ denote the random holomorphic section valued in $H^{0}_{(2)}(X,F,f)$ with law $\P_{\mathrm{st}}$, that is
	\begin{equation}
	\bb{S}=\sum_{j=1}^{d_F(f)}\eta_j S_j,
	\label{eq:3.10nov}
	\end{equation}
	where $\{S_j\}_j$ is an orthonormal basis of $H^0_{(2)}(X,F,f)$ as in \eqref{eq:2.1.10} and $\{\eta_j\}$ is a sequence of i.i.d.\ standard complex Gaussian random variables. In this setting, the random section $\bb{S}_f$ associated to $f$ is given equivalently by 
	\begin{equation}
	\bb{S}_f=T^\sharp_f \bb{S}=T_f \bb{S}.
	\label{eq:3.11nov}
	\end{equation}

The following lemma is an easy modification of \cite[Lemma 4.12]{DrLM:2023aa}
\begin{Lemma}\label{lm:4.13ss}
Assume $d_F(f)\in \N_{\geq 1}\cup\{\infty\}\,$, 
$0\neq f\in  \mathcal{L}^\infty_{\mathrm{c}}(X,\R)$. 
For any nonzero $S\in H^{0}_{(2)}(X,F)$, the random variable on 
$(H^{0}_{(2)}(X,F,f), \P_{f})$ defined via 
$$H^{0}_{(2)}(X,F,f)\ni 
s\mapsto \langle s, S\rangle_{\cLL(X,F)}\in\C$$
is a centered complex Gaussian variable with 
variance $\|T_{f}S\|^{2}_{\cLL(X,F)}$. 
\end{Lemma}

In the case $d_F(f)=\infty$, we consider the orthonormal basis consisting of eigensections of $T_f$ as in \eqref{eq:2.1.10}. Setting
\begin{equation}
\eta_j=\frac{1}{\lambda_j}\langle \bb{S}_f,S_j\rangle_{\cLL(X,F)}
\end{equation}
the $\{\eta_j\}_{j=1}^\infty$ form a sequence of i.i.d.\ standard complex Gaussian random variables. Then our random section $\bb{S}_f$ associated to $f$ is given equivalently by the formula
	\begin{equation}
	\bb{S}_f=\sum_{j=1}^{\infty}\eta_j \lambda_j S_j.
	\label{eq:3.12nov}
	\end{equation}
The well-definedness of the sum in \eqref{eq:3.12nov} is already discussed in \cite[Proposition 2.1]{DrLM:2023aa}. Also note that \eqref{eq:3.12nov} then is consistent with \eqref{eq:3.10nov} and \eqref{eq:3.11nov} which treated the  case $d_F(f)<\infty$.

\subsection{Currents and the Poincar\'e-Lelong formula}\label{sectionC}
The zero-set of a holomorphic section is a complex analytic set which 
is in general singular. The analytic tool used to deal with 
singularities in complex geometry is the theory of currents, introduced
by de Rham \cite{dR55} (see \cite{DemAG,GH94} and especially 
\cite{Fed69} for complete expositions).  

Let $X$ be a complex manifold of dimension $n$ and let $E$
be a complex vector bundle on $X$. The space
of smooth sections of $E$ is denoted by $\cC^\infty(X,E)$
and is endowed with the $\cC^\infty$-topology of uniform
convergence of all derivatives on compact sets.
The space
of smooth sections of $E$ with compact support 
is denoted by $\cC^\infty_{\rm c}(X,E)$
and is endowed with the topology of inductive limit
of spaces of smooth sections with support on a given compact set.
In particular, we denote by 
$\Omega^{n-1,n-1}_{\mathrm{c}}(X)$ the space
of smooth $(n-1,n-1)$-forms with compact support.

The space of $(1,1)$-currents on $X$ 
is the topological dual of the space $\Omega^{n-1,n-1}_{\mathrm{c}}(X)$ 
(called test forms in this context). In the sequel, we 
let $\langle T, \varphi\rangle$ be the pairing between a $(1,1)$-current $T$ 
and a test form $\varphi\in \Omega^{n-1,n-1}_{\mathrm{c}}(X)$.
A $(1,1)$-current is called of order $k\in\N_0$ if it is
continuous in the $\cC^k$-topology, equivalently, it
extends as a linear continuous functional to the space
of $(n-1,n-1)$-forms of class $\cC^k$ with compact support.

\begin{Definition}\label{def:3.5de}
For an open subset $U\subset X$, if $T$ is a $(1,1)$-current of order $0$ (for example, a poisitive $(1,1)$-current), we define the following norm of $T$ on $U$ for $\alpha\in\N$,
\begin{equation}
\|T\|_{U,-\alpha}:=\sup|\langle T,\varphi\rangle|,
\end{equation}
where $\varphi$ runs over all test forms in $\Omega^{n-1,n-1}_{\mathrm{c}}(U)$ with $|\varphi|_{\mathscr{C}^\alpha}\leq 1$.
\end{Definition}

For any analytic hypersurface $V\subset X$, we define
the current of integration $[V]$ on $V$ by 
$$\varphi\mapsto \int_{V}\varphi:=\int_{V_{\mathrm{reg}}}\varphi,
\quad \varphi\in\Omega^{n-1,n-1}_{\mathrm{c}}(X),$$ 
where 
$V_{\mathrm{reg}}$ is the regular set of $V$ (a complex submanifold 
of codimension $1$). By a theorem 
of Lelong (\cite[p.\ 32]{GH94} \cite[III-2.7]{DemAG}) the current 
of integration on $V$ is a closed positive $(1,1)$-current (hence a stronly positive $(1,1)$-current due to \cite[III-1.9]{DemAG}).
It is clear that $[V]$ is a current of order $0$ on $X$.

Let $F$ be a holomorphic line bundle on $X$.
For a holomorphic section $s\in H^{0}(X,F)\setminus\{0\}$ the divisor 
$\Div(s)$ of $s$ is defined as the formal sum
\begin{equation}
	\Div(s)=\sum_{V \subset Z(s)} \mathrm{ord}_{V}(s) V,
	\label{eq:divisor}
\end{equation}
where $V$ runs over all the irreducible analytic hypersurfaces 
contained in $Z(s)$, and $\mathrm{ord}_{V}(s)$ denotes the 
vanishing order of $s$ along $V$.
Let $Z(s)$ denote the 
set of zeros of $s$, which is 
a purely $1$-codimensional analytic subset of $X$. The current
of integration (with multiplicities) on the divisor $\Div(s)$
is defined by 
\begin{equation}
[\Div(s)]=\sum_{V \subset Z(s)} \mathrm{ord}_{V}(s) [V],\quad
\big\langle[\Div(s)],\varphi\big\rangle=
\int_{\Div(s)}\varphi:=
\sum_{V \subset Z(s)} \mathrm{ord}_{V}(s)\int_V\varphi\,.
	\label{eq:divisorint}
\end{equation}
Assume that $F$ is endowed with a smooth Hermitian
metric $h_F$. By the Poincar\'{e}-Lelong formula \cite[Theorem 2.3.3]{MM07} we have
\begin{equation}
	[\Div(s)]=\frac{\sqrt{-1}}{2\pi}
	\partial\bar{\partial}\log|s|^{2}_{h_{F}}+ c_{1}(F,h_{F}),
	\quad\text{for $s\in H^{0}(X,F)$}.
	\label{eq:2.25}
\end{equation}
This important formula is crucial for our purposes. It
links the zero-divisor to the curvature and to the
logarithm of the pointwise norm of a section, 
which is analytically easier to tackle and allows the introduction
of the Bergman kernel or Berezin-Toeplitz kernel into the picture.   

\subsection{Expectation of zeros of Gaussian {$\cLL$}-holomorphic sections}\label{section2.2}
We use the same notation as in Section \ref{ss:4.3b}. In this Section we always assume $d_F(f)\in \N_{\geq 1}\cup\{\infty\}\,$ with $0\neq f\in  \mathcal{L}^\infty_{\mathrm{c}}(X,\R)$. Let $T_f^2(x,y)$ denote the smooth kernel of the operator $T_f^2$. We consider the orthonormal basis of $H^0_{(2)}(X,F,f)$ consisting of eigensections of $T_f$ as in \eqref{eq:2.1.10}, then we have the nonnegative function on $X$,
\begin{equation}
T_f^2(x,x)=\sum_{j=1}^{d_F(f)} \lambda_j^2 |S_j(x)|^2_{h_F}.
\end{equation}

Similar to the proof of \cite[Lemma 4.13]{DrLM:2023aa}, we get that the $(1,1)$-current $\partial\overline{\partial} \log {T_f^2(x,x)}$ is well-defined on $X$. If we proceed as in \cite[Section 2.3]{DrLM:2023aa}, in particular, as in the proof of \cite[Lemma 2.6]{DrLM:2023aa}, we conclude the following result.
\begin{Lemma}[Definition of positive current {$\gamma_f(X,F)$}]\label{lem:currentf}
Assume $d_F(f)\in \N_{\geq 1}\cup\{\infty\}\,$ with $0\neq f\in  \mathcal{L}^\infty_{\mathrm{c}}(X,\R)$. The following current is a closed positive $(1,1)$-current on $X$,
\begin{equation}
\gamma_f(F,h_F):=c_1(F,h_F)+\frac{\sqrt{-1}}{2\pi} \partial\overline{\partial} \log {T_f^2(x,x)}.
\end{equation}
\end{Lemma}

The base locus of $H^{0}_{(2)}(X,F,f)$ is the proper analytic set
\begin{equation}
	\Bl_f(X,F):=\big\{x\in X\;: \; \text{$s(x)=0$ for all 
	$s\in H^{0}_{(2)}(X,F,f)$}\big\}.
	\label{eq:intro2.1}
\end{equation}
Then $\Bl_f(X,F)=\{x\in X\,:\,T^2_f(x,x)=0\}$. Hence $\gamma_f(F,h_{F})$
is a smooth form if $\Bl_f(X,F)=\varnothing$. Moreoever, when $f\geq 0$, we also have $\Bl_f(X,F)=\{x\in X\,:\,P(x,x)=0\}=: \Bl(X,F)$. 

\begin{Theorem}\label{thm:3.6jan}
Assume $d_F(f)\in \N_{\geq 1}\cup\{\infty\}\,$ with $0\neq f\in  \mathcal{L}^\infty_{\mathrm{c}}(X,\R)$. Let $\bb{S}_f$ be the random $\cLL$-holomorphic section of $F$ associated to $f$ defined in Definition \ref{def:4.11}. Then $[\mathrm{Div}(\bb{S}_f)]$ is random variable valued in the space of $(1,1)$-currents on $X$ (equipped with weak topology). Moreover, the expectation $\E[[\mathrm{Div}(\bb{S}_f)]]$ exists as a closed positive $(1,1)$-current on $X$ and we have
\begin{equation}
\E[[\mathrm{Div}(\bb{S}_f)]]=\gamma_f(F,h_F).
\label{eq:3.20Nov}
\end{equation}
\end{Theorem}
\begin{proof}
For $m\in\N_{>0}$, consider the random $(1,1)$-current $\mu_m(\bb{S}_f)$ as follows, for any test form $\varphi\in \Omega^{n-1,n-1}_{\mathrm{c}}(X)$,
\begin{equation}
\langle \mu_m(\bb{S}_f),\varphi\rangle:=\left\langle \frac{\sqrt{-1}}{2\pi}
	\partial\bar{\partial}\log \big(|\bb{S}_f|^{2}_{h_{F}}+\frac{1}{m}\big)+ c_{1}(F,h_{F}), \varphi\right\rangle.
\end{equation}
It is clear that $\langle \mu_m(\bb{S}_f),\varphi\rangle$ is a random variable valued in $\C$ (i.e., a measurable function on the underlying probability space). Moreover, by \eqref{eq:2.25} and dominated convergence, we have
\begin{equation}
\lim_{m\rightarrow \infty} \langle \mu_m(\bb{S}_f),\varphi\rangle=\langle [\mathrm{Div}(\bb{S}_f)],\varphi\rangle.
\end{equation}
Therefore, $\langle [\mathrm{Div}(\bb{S}_f)],\varphi\rangle$ is a random variable, and then $[\mathrm{Div}(\bb{S}_f)]$ is a well-defined random variable valued in the space of $(1,1)$-currents with respect to the weak topology. 

Analogous to \cite[Theorem 1.4]{DrLM:2023aa}, using the computations as in \cite[Proof of Theorem 1.1]{DrLM:2023aa}(cf.\ \cite[Proof of 
Proposition 4.2]{CM11}), we can conclude that for $\varphi\in\Omega^{n-1,n-1}_{\mathrm{c}}(X)$,
$$ \E[\langle [\mathrm{Div}(\bb{S}_f)],\varphi\rangle]=\langle \gamma_f(F,h_F),\varphi\rangle\in\C, $$
then we get \eqref{eq:3.20Nov}. This way, we complete our proof.
\end{proof}

\section{Asymptotics of random zeros of {$\cLL$}-holomorphic sections}
\label{ss4Jan}
Let us recall the main assumption in the semi-classical limit -- Condition \ref{condt:main} -- as follows: we 
always assume $(X,g^{TX})$ to be a complete Riemannian manifold, and 
that there exist $\varepsilon_{0}>0$, 
$C_{0}>0$ 
such that
\begin{equation}
	\sqrt{-1}R^{L} \geq \varepsilon_{0} \Theta,\; 
	\sqrt{-1}(R^{\mathrm{det}}+R^{E})>-C_{0}\Theta,\;\; |\partial 
	\Theta|_{g^{TX}}< C_{0}.
	\label{eq:mainassumptions}
\end{equation}
If $X$ is compact and $(L,h_{L})$ is positive, then the above 
conditions always hold true. Due to our assumptions on $(X,\Theta)$ and $L, E$, for $p\gg 0$, we have $d_p:= \dim_\C H^0_{(2)}(X,L^p\otimes E)\geq 1$. We may assume that $d_p\geq 1$ for all $p\geq 1$.

\begin{Definition}\label{def:4.1nov}
For $p\geq 1$, let $\bb{S}_{f,p}$ denote the random $\cLL$-holomorphic section of $L^p\otimes E$ associated to the nontrivial $f\in \mathcal{L}^{\infty}_{\mathrm{c}}(X,\R)$ defined in Definition \ref{def:4.11}.
\end{Definition}

The goal of this section is to study the asymptotic behavior of the zeros of random section $\bb{S}_{f,p}$ as $p\rightarrow \infty$, where the asymptotic behavior of $T_{f,p}(x,x')$ is a crucial ingredient. As we saw in Section \ref{ss2nov}, in order to analyze $T_{f,p}(x,x')$ in the semi-classical limit, we always need to assume certain regularities on $f$.

We start with the following lemma, which enable us to apply the results from Section \ref{section2} for the line bundles $L^p\otimes E$. Recall that the function $\boldsymbol{b}_0(x)\neq 0$ is given in \eqref{eq:2.1.4June}. 

\begin{Lemma}\label{lemma:4.2jan}
For $f\in \mathcal{L}^{\infty}_{\mathrm{c}}(X,\R)$, if there exists an open subset $U\neq \varnothing$ such that $f|_U\in\mathscr{C}^1(U)$ and $f$ never vanishes on $U$, then for all sufficiently large $p$, we have $T_{f,p}\neq 0$, hence $d_p(f)\geq 1$.

Moreover, fix a proper open subset $V$ of $U$, we have
\begin{equation}
d_p(f)\geq \frac{p^n}{|f|^2_{\mathcal{L}^\infty(X)}}\int_V f^2(x)\bb{b}_0(x)\mathrm{dV}(x)+\mathcal{O}(p^{n-1}).
\label{eq:4.3parisfeb24}
\end{equation}
\end{Lemma}
\begin{proof}
If $T_{f,p}=0$, then for all $x\in X$, $T_{f,p}(x,x)\equiv 0$. Since $f|_U\in\mathscr{C}^1(U)$, we can apply Theorem \ref{thm:2.5Nov} to $f$ on $U$ so that fix $x\in U$, we have
\begin{equation}
T_{f,p}(x,x)=p^n f(x)\boldsymbol{b}_0(x)+\mathcal{O}_x(p^{n-1}).
\end{equation}
By our assumption $f(x)\neq 0$, then for all sufficiently large $p$, $T_{f,p}(x,x)\neq 0$. Therefore, we have $T_{f,p}\neq 0$.

Note that
\begin{equation}
\int_X T_{f,p}^2(x,x)\mathrm{dV}(x)\geq \int_V T_{f,p}^2(x,x)\mathrm{dV}(x).
\label{eq:4.5parisfeb24}
\end{equation}
Similarly, we can apply Theorem \ref{thm:2.5Novfg} to $f$ on $U$, then \eqref{eq:4.3parisfeb24} follows from Lemma \ref{lm:3.1} and \eqref{eq:4.5parisfeb24}.
\end{proof}

\subsection{Bounded measurable functions and their supports}
This Section is a preparation for the proof of Theorem \ref{thm:4.6Jan2024}. Recall that the quantity $r(f,U)$ for $f\in \mathcal{L}^{\infty}_{\mathrm{c}}(X,\R)$ is defined in Definition \ref{def:rfU}. Now we give some examples of function $f$ such that $r(f,U)$ has a good geometric sense.
\begin{Example}
\begin{itemize}
	\item[i)] For $t\in \R$, set $\sin^+(t\pi):=\max\{0,\sin(t\pi)\}$. Now fix $N\in\N_{\geq 1}$ and $R\gg 0$, let $ \mathbb{B}(0,R)\subset \C$ be the open ball in $\C$ with radius $R$, and define for $z=x+\sqrt{-1}y\in\C$,
	\begin{equation}
	f_{N,R}(z):=\bb{1}_{\mathbb{B}(0,R)}(z) \sin^+(Nx\pi)\sin^+(Ny\pi).
	\end{equation}
	Then $f_{N,R}$ is smooth almost everywhere on $\C$, and we have for $R\gg 0$,
	\begin{equation}
	r(f_{N,R}, \mathbb{B}(0,R))=\frac{\sqrt{2}}{2N}.
	\end{equation}
	In Figure \ref{fig:supportf}, it shows the support in black blocks of the function $f_{8,2}(z)$ (that is, $N=8$, $R=2$) on $\C$.
	\begin{figure}[h]
\centering
\includegraphics[width=0.4\textwidth]{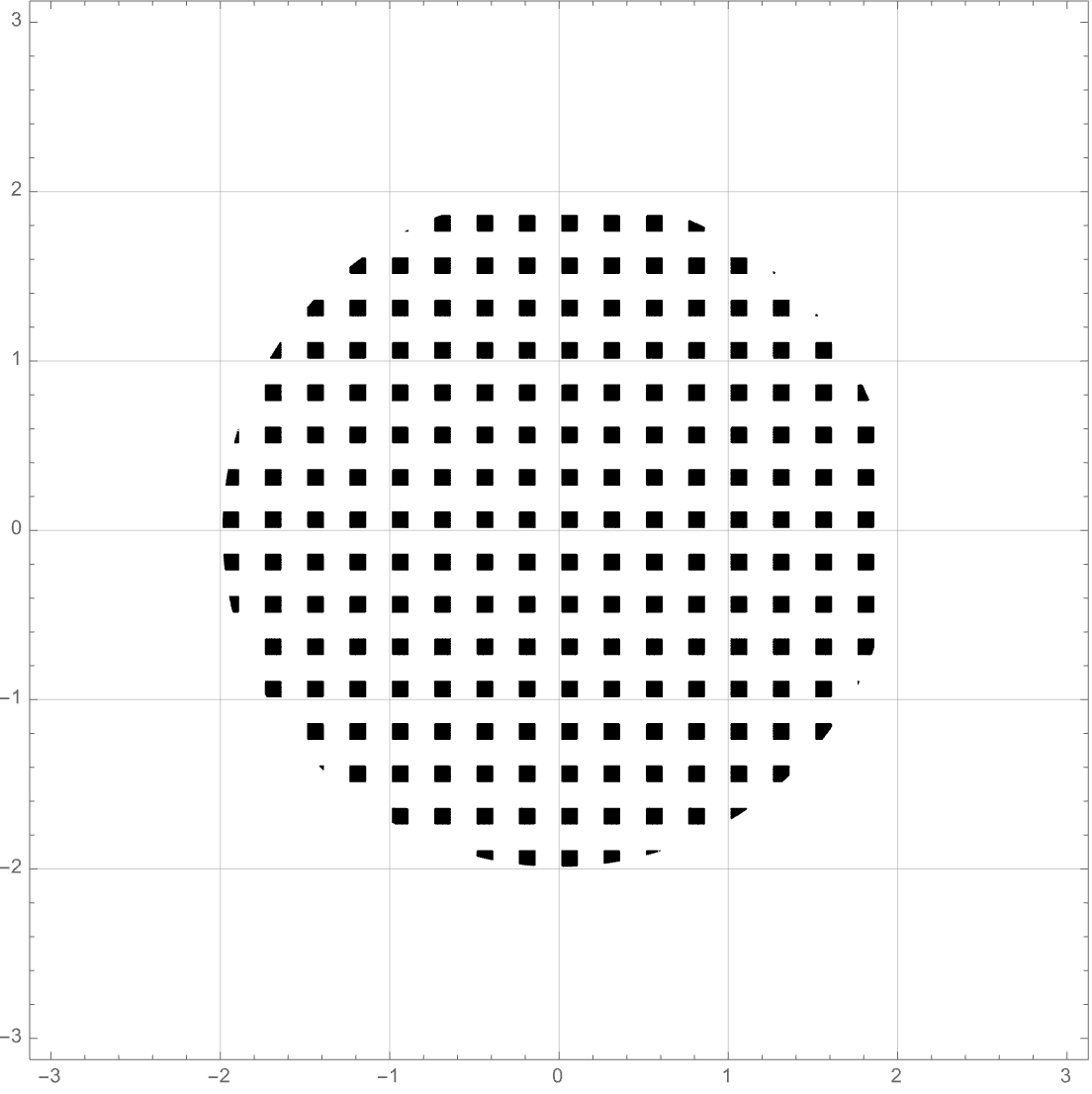}
\caption[Support of $f_{8,2}(z)$]{Support of $f_{8,2}(z)$ on $\C$ plotted as black blocks}
\label{fig:supportf}
\end{figure}
\item[ii)] Let $\{B_j\}_{j=1}^N$ be a finite collection of disjoint geodesic open ball of respective radius $r_j>0$ in a compact Hermtian manifold $(X,\Theta)$, set $V:=X\setminus\cup_j B_j$, then $f:=\bb{1}_V$ is a function which is smooth almost everywhere on $X$, and $r(f,X)=\max\{ r_j\;:\; j=1,\ldots,N\}$.
\item[iii)] Let $(\Sigma,\omega)$ be a compact Riemann surface with a K\"{a}hler form $\omega$. Let $K$ be the $1$-skelton of a (smooth) triangluation of $\Sigma$. Then for a sufficient small $\varepsilon>0$, let $K(\varepsilon)$ be the closed subset defined as the $\varepsilon$-neighbourhood of $K$ in $\Sigma$ (which is still a proper subset). This is a thick graph embbeded in $\Sigma$. Then we can use the quantity $r(\bb{1}_{K(\varepsilon)},\Sigma)$ to quantify the approximation of the thick graph $K(\varepsilon)$ to $\Sigma$. Taking a sequence of strict refinements of $K$ and a suitable decreasing sequence of $\varepsilon$, then the sequence of the quantities $r(\bb{1}_{K(\varepsilon)},\Sigma)$ will converge to $0$.
\end{itemize}
\end{Example}

Let $g^{T\R^{2n}}$ denote the standard Riemannian
metric on $\C^n\simeq \R^{2n}$ given by the 
standard Euclidean inner product on $\R^{2n}$. 
The following proposition is elementary.
\begin{Proposition}\label{prop:4.3jan}
For every compact subset $K\subset X$ there exist 
constants $t_K>0$, $c_K\geq 1$ such that for every point
$x\in K$, the open geodesic ball $\mathbb{B}(x,t_K)$ in 
$(X,g^{TX})$ is always included in a connected, 
open holomorphic local chart $(V_x\subset X,\phi_x)$ 
centered at $x$ such that $\phi_x$ is a biholomophism 
between the open subset $V_x$ and $\phi_x(V_x)\subset \C^n$ 
with the following properties:
\begin{itemize}
\item  $\phi_x(x)=(0,\ldots,0)$, and $g^{TX}_x=
\phi_x^\ast(g^{T\R^{2n}}_0)$;
\item $\phi(V_x)$ is a convex domain of $\C^n$;
\item let $\mathrm{dist}^{\phi_x}(\cdot,\cdot)$ 
be the distance function on $V_x\times V_x$ 
associated to the Riemannian metric 
$\phi_x^\ast(g^{T\R^{2n}})$, then for any $y,y'\in V_x$,
\begin{equation}
\frac{1}{c_K} \mathrm{dist}^{\phi_x}(y,y')\leq 
\mathrm{dist}(y,y')\leq c_K \mathrm{dist}^{\phi_x}(y,y').
\end{equation}
\end{itemize}
As a consequence, for a given $x\in K$, and for 
$f \in\mathcal{L}^\infty(X,\R_{\geq 0})$, we have
\begin{equation}
\frac{r(f,V_x)}{c_K}\leq r_{\C^n}(f\circ\phi_x^{-1},
\phi_x(V_x))\leq c_K r(f,V_x),
\label{eq:4.6Feb}
\end{equation}
where $r_{\C^n}(f\circ\phi_x^{-1},\phi_x(V_x))$ 
is the number defined as in \eqref{eq:4.3Jan2024}
with respect to the standard Riemannian metric $g^{T\R^{2n}}$ on $\C^n$.
\end{Proposition}

\begin{Remark}
The choices of $(t_K, c_K)$ for a given $K$ in 
Proposition \ref{prop:4.3jan} are not unique, 
in fact we can always make $t_K$ as small as we want. 
Moreover, if we take a sufficiently small $t_K>0$, 
then the constant $c_K$ will be very close to $1$.
\end{Remark}

\subsection{A concentration estimate; proof of Theorem \ref{thm:4.6Jan2024}}
In this section, we prove Theorem \ref{thm:4.6Jan2024}
which provides a general concentration estimate for the difference between
the zeros of the random section $\bb{S}_{f,p}$ and the expected 
limit $c_1(L,h_L)$ on a given domain, provided that the function $f$ 
is supported on the large part of this domain. 
Recall that the norm $\|\cdot\|_{U,-2}$ for the $(1,1)$-currents
is given Definition \ref{def:3.5de} with $\alpha=2$, and that 
the quantity $m(U)\in\N$ defined in \eqref{eq:1.9intropart},
then Theorem \ref{thm:4.6Jan2024} is a consequence 
of the following proposition. 

\begin{Proposition}\label{thm:4.3Jan}
With the same assumption in Theorem \ref{thm:4.6Jan2024} for $X$ and for the line bundles $L,E$. Fix a pair of nonempty open subsets $(U,U')$ of $X$ with $U\subset \subset U'$, then there exists a constant $r(U,U')>0$ such that $f\in \mathcal{L}^{\infty}_{\mathrm{c}}(X,\R)$ which is of $\mathscr{C}^{m(U')+1}$ almost everywhere on $U'$ with 
\begin{equation}
\left(\frac{r(f,U')}{r(U,U')}\right)^{1/(2n+2)}<\frac{1}{2}\, ,
\end{equation}
then $\, U\cap\, \esupp\, f\neq\varnothing\,$, and for any 
	$\delta > \left(\frac{r(f,U')}{r(U,U')}\right)^{1/(2n+2)}$, we have a constant $C=C_{U',f,\delta}>0$ such that for all sufficiently large $p$
\begin{equation}
		\P\Big(
		\int_{U}\Big|\log{\big|\bb{S}_{f,p}(x)\big|_{h_{p}}}\Big| 
		\,\mathrm{dV}(x) >
		\delta p \Big)\leq e^{-Cp^{n+1}}.
		\label{eq:4.67feb24}
	\end{equation}
\end{Proposition}

Now we explain how to get Theorem \ref{thm:4.6Jan2024}
from Proposition \ref{thm:4.3Jan}, and the proof of Proposition 
\ref{thm:4.3Jan} is defered to Section \ref{Section3.3}.
\begin{proof}[Proof of Theorem \ref{thm:4.6Jan2024}]
Combining the first part of Proposition \ref{thm:4.3Jan} with Lemma \ref{lemma:4.2jan}, we get that for all sufficiently large $p$, $T_{f,p}\neq 0$, that is $d_{p}(f)\geq 1$. Then the results from Section \ref{section2} apply. In particular, the random section $\bb{S}_{f,p}$ is not the zero variable.

The Poincar\'{e}-Lelong formula \eqref{eq:2.25} shows that 
	\begin{equation}
\frac{\sqrt{-1}}{\pi}\partial\overline{\partial}
\log{|\bb{S}_{f,p}|_{h_{p}}}=[\Div(\bb{S}_{f,p})]-pc_{1}(L,h_L)-c_1(E,h_E)
	\label{eq:3.4.13Jan}
	\end{equation}
as an identity of $(1,1)$-currents on $X.$	Let $U$ be the open subset as assumed in Theorem \ref{thm:4.6Jan2024}. Now fix $\varphi\in \Omega^{n-1,n-1}_{\mathrm{c}}(X)$ with $\mathrm{supp}\ 
\varphi\subset U$ (that is, $\varphi\in \Omega^{n-1,n-1}_{\mathrm{c}}(U)$). Note that $\mathrm{supp}\,\varphi\subset U\subset U'$.
Then
\begin{equation}
\begin{split}
\Big\langle\frac{1}{p}[\Div(\bb{S}_{f,p})],\varphi\Big\rangle
-\int_{X}c_{1}(L,h_L)\wedge \varphi=\frac{\sqrt{-1}}{p\pi}\int_{X}
\log{|\bb{S}_{f,p}|_{h_{p}}}\,\partial\overline{\partial}\varphi +\frac{1}{p}\langle c_1(E,h_E),\varphi\rangle.
		\end{split}
		\label{eq:3.4.14Jan}
	\end{equation}
Since $\varphi$ has a compact support in $U$, so has 
$\partial\overline{\partial}\varphi$. 
Then
	\begin{equation}
	\begin{split}
	\left|\frac{\sqrt{-1}}{p\pi}\int_{X}\log{|\bb{S}_{f,p}|_{h_{p}}}\,
		\partial\overline{\partial}\varphi\right|\leq 
		\frac{|\varphi|_{\mathscr{C}^2(U)}}{p\pi}\int_{U}
		\big|\log{|\bb{S}_{f,p}(x)|_{h_{p}}}\big|\,\mathrm{dV}(x).
		\end{split}
		\label{eq:3.4.17Jan}
	\end{equation}
Applying Proposition \ref{thm:4.3Jan}
to get a constant $r(U,U')>0$, we fix a $f$ as in the theorem 
with $$\left(\frac{r(f,U')}{r(U,U')}\right)^{1/(2n+2)} <\frac{1}{2},$$ 
then $\, U\cap\, \esupp\, f\neq\varnothing\,$. 
When we take $$\delta > \left(\frac{r(f,U')}{r(U,U')}\right)^{1/(2n+2)},$$ 
we can always fix a sufficiently small $\varepsilon>0$ such that
$$\delta-2\varepsilon > \left(\frac{r(f,U')}{r(U,U')}\right)^{1/(2n+2)}.$$ 
Since the term $\frac{1}{p} c_1(E,h_E)$ converges to $0$ as 
$p\rightarrow \infty$, then there exists an integer $p_0\in\N$
(depending on $(E,h^E)$) such that for any 
$\varphi\in \Omega^{n-1,n-1}_{\mathrm{c}}(U)$ and all $p\geq p_0$,
\begin{equation}
\left|\frac{1}{p}\langle c_1(E,h_E),\varphi\rangle\right|
\leq \frac{\varepsilon |\varphi|_{\mathscr{C}^2(U)}}{\pi}\,\cdot
\end{equation} 
Applying Proposition \ref{thm:4.3Jan} to right-hand side of 
\eqref{eq:3.4.17Jan} with $\delta-\varepsilon$, we get
\begin{equation}
		\P\Big(\frac{1}{p}
		\int_{U}\Big|\log{\big|\bb{S}_{f,p}(x)\big|_{h_{p}}}\Big| 
		\,\mathrm{dV}(x)>
		\delta-2\varepsilon \Big)\leq e^{-Cp^{n+1}}.
		\label{eq:4.67Jan}
	\end{equation}
For $p\geq p_0$, except the event from \eqref{eq:4.67Jan} of probability 
$\leq e^{-Cp^{n+1}}$, we have for all $\varphi\in \Omega^{n-1,n-1}_{\mathrm{c}}(U)$,
\begin{equation}
\begin{split}
		&\Big|\Big\langle\frac{1}{p}[\Div(\bb{S}_{f,p})] 
		-  c_{1}(L,h_L),\varphi\Big\rangle\Big|\\
		&\leq   
		\frac{|\varphi|_{\mathscr{C}^2(U)}}{p\pi}\int_{U}
		\big|\log{|\bb{S}_{f,p}(x)|_{h_{p}}}\big|\,\mathrm{dV}(x)+ |\frac{1}{p}\langle c_1(E,h_E),\varphi\rangle|\\
		&\leq \frac{1}{\pi}\left(|\varphi|_{\mathscr{C}^2(U)}(\delta-2\varepsilon)+\varepsilon |\varphi|_{\mathscr{C}^2(U)}\right) \leq |\varphi|_{\mathscr{C}^2(U)}\frac{\delta-\varepsilon}{\pi},
		\end{split}
	\end{equation}
Equivalently, except the event in \eqref{eq:4.67Jan} of probability 
$\leq e^{-Cp^{n+1}}$, we have 
\begin{equation}
\left\|\frac{1}{p}[\Div(\bb{S}_{f,p})] 
		-  c_{1}(L,h_L)\right\|_{U,-2}\leq \frac{\delta-\varepsilon}{\pi}.
\end{equation}
This way, we get \eqref{eq:4.7Jan}. 
To prove \eqref{eq:4.73Jan}, applying Borel-Cantelli lemma to 
\eqref{eq:4.7Jan}, we get
\begin{equation}
		\P\Big(\limsup_{p\rightarrow +\infty} \ 
		\Big\|\frac{1}{p}[\Div(\bb{S}_{f,p})] 
		-  c_{1}(L,h_L)\Big\|_{U,-2}\leq \frac{\delta}{\pi} \  \Big)=1.
		\label{eq:4.77Jan}
	\end{equation}
By taking a sequence of $\delta$'s that decreases to $\delta_0(f)$,
we infer \eqref{eq:4.73Jan} from \eqref{eq:4.77Jan}. 
\end{proof}

\subsection{Local sup-norms of random $\cLL$-holomorphic sections}

Before proving Proposition \ref{thm:4.3Jan}, we need to investigate 
the probabilities for $\mathcal{M}^{U}_{p}(\bb{S}_{f,p})$
taking both atypically large and small values, respectively. 
We follow the same strategy as in \cite[Section 3.3]{DrLM:2023aa}.

Let $U\subset X$ be a relatively compact open subset. For $s_{p}\in H^{0}(X, L^{p}\otimes E)$, we set
\begin{equation}
	\mathcal{M}^{U}_{p}(s_{p})=\sup_{x\in 
	U}|s_{p}(x)|_{h_{p}}<+\infty.
	\label{eq:5.2.1}
\end{equation}

\begin{Proposition}\label{prop:3.14ss}
 Assume Condition \ref{condt:main} to hold. Fix $0\neq f\in \mathcal{L}^{\infty}_{\mathrm{c}}(X,\R_{\geq 0})$. Let $U\subset X$ be a relatively compact open subset, then for any $\delta>0$, there exists a constant $C_{U,\delta}>0$ such 
	that for $p\in\N_{>1}$, 
	\begin{equation}
		\P\big (\mathcal{M}^{U}_{p}(\bb{S}_{f,p})\geq 
		e^{\delta p}\big)\leq e^{-\delta 
		p^{n+1}+C_{U,\delta}p^{n}\log{p}}\,.
		\label{eq:5.2.2}
	\end{equation}
\end{Proposition}
\begin{proof}
	We fix $\delta>0$ and let $r>0$ be sufficiently small so that we can 
	choose a finite set of points $\{x_{j}\}_{j=1}^{\ell}\subset U$ such 
	that the geodesic open balls $B^{X}(x_{j},r)$, $j=1,\ldots, \ell$, 
	form an open covering of $\overline{U}$. Since $r$ is 
	sufficiently small, then we can assume that each larger ball
	$B^{X}(x_{j},2r)$ lies in a complex chart (hence viewed as an 
	open subset of $\C^{n}$), and that for each $j$, we can fix a 
	local holomorphic frame $e_{L,j}$ (resp. $e_{E,j}$) of $L$ (resp. $E$) on a neighborhood of $B^{X}(x_{j}, 2r)$ 
	with $\sup_{x\in B^{X}(x_{j},2r)} |e_{L,j}(x)|_{h_{L}}=1$ and $\sup_{x\in B^{X}(x_{j},2r)} |e_{E,j}(x)|_{h_{E}}\leq 1$.
	Set
	\begin{equation}
		\nu=\min\big\{ \inf_{x\in 
		B^{X}(x_{j},2r)}|e_{L,j}(x)|_{h_{L}} \,:\,j=1 , \ldots, \ell\big\}.
		\label{eq:3.14}
	\end{equation}
	It is clear that $0<\nu\leq 1$. By fixing $r$ small enough, we 
	can and do assume that 
	 \begin{equation} \label{eq:nuBd}
	 -\log{\nu}\leq 
	\frac{\delta}{4}\,\cdot
	\end{equation} 	
	Then there 
	exists a constant $C>0$ such that for each $j=1,\ldots, \ell$, if 
	$f$ is a holomorphic function on a neighborhood of 
	$B^{X}(x_{j},2r)$, then
	\begin{equation}
		\sup_{x\in B^{X}(x_{j},r)}|f(x)|\leq 
		C\|f\|_{\mathcal{L}^{2}(B^{X}(x_{j},2r))},
		\label{eq:3.16}
	\end{equation}
	where the volume form $\mathrm{dV}(x)$ on $X$ is used in the norm 
	$\|\cdot\|_{\mathcal{L}^{2}(B^{X}(x_{j},2r))}$.
	Note that the choices of $x_{j}$, $r$, $\ell$, and the constants 
	$\nu$, $C$
	are independent of the tensor power $p$. Set 
		$\widetilde{U}=\cup_{j} B^{X}(x_{j},2r)\supset U$.
	For $p\in\N, s_{p}\in H^{0}(X,L^{p}\otimes E)$, on each $B^{X}(x_{j},2r)$, we write
	\begin{equation}
		s_{p}|_{B^{X}(x_{j},2r)}=f_{j}e_{L,j}^{\otimes p}\otimes e_{E,j},
		\label{eq:3.15}
	\end{equation}
	where $f_{j}$ is a holomorphic function on the chart in $\C^{n}$ 
	corresponds to $B^{X}(x_{j},2r)$. Then we have
	\begin{equation}
		\begin{split}
			\mathcal{M}^{U}_{p}(s_{p})=\sup_{x\in U}|s_{p}(x)|_{h_{p}}&\leq \max_{j} \sup_{x\in 
			B^{X}(x_{j},r)}|f_{j}(x)|\\
			&\leq C\max_{j} 
			\{\|f_{j}\|_{\mathcal{L}^{2}(B^{X}(x_{j},2r))}\}\\
			&\leq 
			\frac{C'}{\nu^{p}}\max_{j}\{\|s_{p}\|_{\mathcal{L}^{2}(B^{X}(x_{j},2r),L^{p}\otimes E)}\}\\
			&\leq 
			\frac{C'}{\nu^{p}}\|s_{p}\|_{\mathcal{L}^{2}(\widetilde{U},L^{p}\otimes E)},
		\end{split}
		\label{eq:3.17bis}
	\end{equation}
where the constant $C'$ independent of $p$ is determined by 
$$C/C'=\min\big\{ \inf_{x\in 
B^{X}(x_{j},2r)}|e_{E,j}(x)|_{h_{E}} \,:\,j=1 , 
\ldots, \ell\big\}=: \nu_E.
$$
The next step is to estimate the quantity 
$\E[\|\bb{S}_{f,p}\|^{2p^{n}}_{\mathcal{L}^{2}(\widetilde{U},L^{p})}]$ 
for $p\geq 2$. Applying H\"{o}lder's inequality with 
$\frac{1}{p^{n}}+\frac{p^{n}-1}{p^{n}}=1$, we get 
	\begin{equation}
\E\big[\|\bb{S}_{f,p}\|^{2p^{n}}_{\mathcal{L}^{2}(\widetilde{U},L^{p}\otimes E)}\big]
\leq \Vol(\widetilde{U})^{p^{n}-1}\E\Big[\int_{\widetilde{U}}
|\bb{S}_{f,p}(x)|^{2p^{n}}_{h_{p}}(x)\mathrm{dV}\Big].
		\label{eq:holderpn}
	\end{equation}
For $p\geq 2$, let $\{S^p_m\}_{m=1}^{d_p(f)}$ be an orthonormal basis of $H^0_{(2)}(X,L^p\otimes E,f)\subset H^0_{(2)}(X,L^p\otimes E)$ with respect to $\cLL$-inner product and consisting of eigensections of $T_{f,p}$ with nonzero eigenvalues $\{\lambda^p_m\neq 0\}_m$. As explained in \eqref{eq:3.12nov}, we can construct a sequence of i.i.d.\ standard complex Gaussian variable $\{\eta^p_m\}_m$ such that we can write
\begin{equation}
\bb{S}_{f,p}=\sum_{m=1}^{d_p(f)} \eta^p_m \lambda^p_m S^p_m.
\label{eq:4.34series}
\end{equation}
As in \eqref{eq:3.16}, on a neighborhood of $B^{X}(x_{j}, 2r)$, write
\begin{equation}
		S^{p}_{m}=f^{p}_{m}e_{L,j}^{\otimes p}\otimes e_{E,j}.
		\label{eq:3.18bis}
	\end{equation}
	If $x\in B^{X}(x_{j},2r)$, set
	\begin{equation}
		F_{j}(x)=\sum_{m=1}^{d_{p}(f)} \eta^p_{m} \lambda^p_m f^{p}_{m}(x).
		\label{eq:3.19bis}
	\end{equation}
	Then $F_{j}(x)$ is a complex Gaussian random variable with 
	(total) variance $\sum_{m=1}^{d_{p}(f)}|\lambda^p_m f^{p}_{m}(x)|^{2}$. 
	By our 
	assumption on the local frame $e_{L,j}$, we get
	\begin{equation}
		\sum_{m=1}^{d_{p}(f)}|\lambda^p_m f^{p}_{m}(x)|^{2}\leq \frac{1}{\nu^{2p}\nu_E^2}T^2_{f,p}(x,x).
		\label{eq:3.20bis}
	\end{equation}
	Then we 
	have
	\begin{equation}
\E\big[|F_{j}(x)|^{2p^{n}}\big]=p^{n}!
\Big(\sum_{m=1}^{d_{p}(f)}|\lambda^p_m f^{p}_{m}(x)|^{2}\Big)^{p^{n}}.
		\label{eq:3.21}
	\end{equation}
As a consequence, we get that for $x\in \widetilde{U}$,
\begin{equation}
	\begin{split}
\E\big[|\bb{S}_{f,p}(x)|_{h_{p}}^{2p^{n}}\big]\leq \E\big[|F_{j}(x)|^{2p^{n}}\big]
\leq\frac{p^{n}!}{\nu^{2p^{n+1}}\nu_E^{2p^n}}(T^2_{f,p}(x,x))^{p^{n}}.
	\end{split}
\end{equation}
Since we are in the context of $\sigma$-finite measures and 
	the integrands are non-negative, Tonelli's Theorem applies, so that
	\begin{equation}
		\E\Big[\int_{\widetilde{U}}|\bb{S}_{f,p}(x)|^{2p^{n}}_{h_{p}}\mathrm{dV}(x)\Big]\leq \frac{p^{n}!}{\nu^{2p^{n+1}}\nu_E^{2p^n}}\int_{\widetilde{U}} (T^2_{f,p}(x,x))^{p^{n}}\mathrm{dV}(x).
		\label{eq:3.24}
	\end{equation}
Moreover, by the estimate for $T^2_{f,p}(x,x)$ given in \eqref{eq:June2.44on} on a compact 
subset, there exists a 
	constant $C_{\widetilde{U},f}>0$ (independent of $p$) such that for 
	$p\in\N$, $x\in\widetilde{U}$,
	\begin{equation}
		T^2_{f,p}(x,x)\leq C_{\widetilde{U},f}p^{n}.
		\label{eq:3.25}
	\end{equation}
		Combining \eqref{eq:holderpn} with the above inequalities, we infer that
	\begin{equation}
		\E\Big[\|\bb{S}_{f,p}\|^{2p^{n}}_{\mathcal{L}^{2}(\widetilde{U},L^{p}\otimes E)} \Big]\leq \left(\frac{C_{\widetilde{U}}\Vol(\widetilde{U})}{\nu_E^2}\right)^{p^{n}}\cdot\frac{p^{n}!}{\nu^{2p^{n+1}}}(p^{n})^{p^{n}}.
		\label{eq:3.26}
	\end{equation}
	By applying \eqref{eq:3.17bis} to $\bb{S}_{f,p}$, we get
	\begin{equation}
		\begin{split}
			\E\big[\mathcal{M}^{U}_{p}(\bb{S}_{f,p})^{2p^{n}}\big]&\leq 
		\Big(\frac{C}{\nu^{p}\nu_E}\Big)^{2p^{n}}\E\big[\|\bb{S}_{f,p}\|^{2p^{n}}_{\mathcal{L}^{2}(\widetilde{U},L^{p}\otimes E)}\big]
		\leq \frac{(\widetilde{C}p^{n})^{2p^{n}}}{\nu^{4p^{n+1}}},
		\end{split}
		\label{eq:3.18}
	\end{equation}
	where $\widetilde{C}>0$ is a constant independent of $p$.
Then \eqref{eq:5.2.2} follows from Chebyshev's inequality and the inequality
	$\frac{1}{\nu}\leq e^{\frac{\delta}{4}}$ from \eqref{eq:nuBd}.
\end{proof}

Now we consider the probabilities of small values of 
$\mathcal{M}^{U}_{p}(\bb{S}_{f,p})$, and we will adapt the 
ideas in \cite{SZZ}, \cite{DLM:21,DrLM:2023aa}: viewing $\bb{S}_{f,p}$ as a line bundle valued holomorphic Gaussian field on $X$ and studying its correlation function.

\begin{Proposition}\label{prop:3.17}
Assume Condition \ref{condt:main} to hold.
	 Fix $f\in \mathcal{L}^{\infty}_{\mathrm{c}}(X,\R_{\geq 0})$ and assume that on a small open ball $B$ of $X$, $f|_B$ is $\mathscr{C}^{m(B)+1}(B)$ (the quantity $m(B)$ is given in \eqref{eq:1.9intropart}) and nonvanishing. 	Let $U$ be a relatively compact open subset in $X$ such that $U\cap B\neq \varnothing$, then there exist constants $C_{U}>0, C_{U}'>0$ such 
	that for all $\delta>0$ and $p\in\N$, 
	\begin{equation}
		\P\big(\mathcal{M}^{U}_{p}(\bb{S}_{f,p})\leq 
		e^{-\delta p}\big)\leq e^{-C_{U}\delta p^{n+1}+C_{U}' 
		p^{n}\log{p}}\,.
		\label{eq:5.2.3}
	\end{equation}
\end{Proposition}
\begin{proof}
From the local uniformity of the Berezin--Toeplitz kernel expansions,
as explained in Section \ref{ss5.1a}, it follows that for every compact subset 
$K\subset B$
there exists $p(K)$ such that for all $p\geq p(K)$ and all $x\in K$
we have  $T^2_{f,p}(x,x)>0$. 

Now for any $x\in B$ we fix some $\lambda_{x}\in L_{x}$, $\mu_x\in E_x$ with 
 $|\lambda_{x}|_{h_L}=|\mu_x|_{h_E}=1$, and set
 \begin{equation}
 \xi_{x}=\frac{\langle \lambda_{x}^{\otimes p}\otimes \mu_x, 
 \bb{S}_{f,p}(x)\rangle_{h_{p}}}{\sqrt{T^2_{f,p}(x,x)}}\, ,
\label{eq:53.2.15}
 \end{equation}
 whenever $T^2_{f,p}(x,x)\neq 0$. Then $\xi_{x}$ is a complex Gaussian random variable. Moreover, for 
any two points $x,y\in B$, we have 
 \begin{equation}
	\big|\mathbb{E}[\xi_{x}\overline{\xi}_{y}]\big|
=N_{f,p}(x,y),
\label{eq:53.2.19}
 \end{equation}
 where $N_{f,p}(x,y)$ is defined by \eqref{eq:2.60nov}.

Since $U\cap B\neq \varnothing$, so we can find a small open ball $C:=B(x_0,r)\subset U\cap B$, then $f$ is $\mathscr{C}^{m(B)+1}$ and never vanishes near $C$. Then the asymptotic equations in \eqref{eq:5.2.4Nov} from Corollary \ref{cor:5.1.1} holds true for all $x,y\in C$ and for all $p\gg 0$. In particular, we have
\begin{equation}
		\P\big(\mathcal{M}^{U}_{p}(\bb{S}_{f,p})\leq 
		e^{-\delta p}\big)\leq \P\big(\mathcal{M}^{C}_{p}(\bb{S}_{f,p})\leq 
		e^{-\delta p}\big).
\end{equation}
Next step, by \eqref{eq:5.2.4Nov}, we may proceed using the similar arguments in \cite[Section 
	3.2]{SZZ} or \cite[Subsection 3.3 and Theorem 1.13]{DLM:21}, we can 
	prove a more general version of \eqref{eq:5.2.3} as follows: for a sequence of positive numbers 
$\{\lambda_{p}\}_{p\in\N}$,
\begin{equation}
	\P\big(\mathcal{M}^{U}_{p}(\bb{S}_{f,p})\leq 
	\lambda_{p}\big)\leq e^{Cp^{n}\log \lambda_{p} +C'p^{n}\log 
	p},\quad  p\gg 0.
	\label{eq:53.2.16}
\end{equation}
Then, for any 
$\delta>0$, choosing $\lambda_{p}=e^{-\delta p}$ in \eqref{eq:53.2.16}, we recover 
\eqref{eq:5.2.3}. This completes our 
proof.
\end{proof}
Combining Propositions \ref{prop:3.14ss} and \ref{prop:3.17}, we arrive at the following.
\begin{Corollary}\label{cor:3.18ss}
Let $(X,J,\Theta)$ be a connected Hermitian complex manifold and let $(L,h_{L})$, $(E,h_{E})$ be two holomorphic line bundles on $X$ with smooth Hermitian metrics. Assume Condition \ref{condt:main} to hold. Fix $f\in \mathcal{L}^{\infty}_{\mathrm{c}}(X,\R_{\geq 0})$ and assume that on a small open ball $B$ of $X$, $f|_B$ is $\mathscr{C}^{m(B)+1}(B)$ and nonvanishing. Let $U$ be a relatively compact open subset in $X$ such that $U\cap B\neq \varnothing$, then there exist constants $C=C_{U,f,\delta}>0$ such 
	that for all $\delta>0$ and $p\in\N$,
	\begin{equation}
		\P\left(\left|\log{\mathcal{M}^{U}_{p}(\bb{S}_{f,p})}\right|\geq \delta p\right) \leq e^{-Cp^{n+1}}.
		\label{eq:3.79ss}
	\end{equation}
\end{Corollary}

\subsection{Proof of Proposition \ref{thm:4.3Jan}}\label{Section3.3}
Now we give the proof of Proposition \ref{thm:4.3Jan}, which follows by combining from the arguments in \cite[Section 4.1]{SZZ}
with Corollary \ref{cor:3.18ss}. Since the kernels $T^2_{f,p}(x,y)$ behave in a more complicated way, depending on the values of $f$, than the Bergman kernels, some steps are required necessary modifications from that in \cite[Section 4.1]{SZZ}. In the same time, we also include several explicit function estimates in order to work out a rough formula for the quantity $r(U,U')$ in the statement of Proposition \ref{thm:4.3Jan}.
\begin{proof}[Proof of Proposition \ref{thm:4.3Jan}]
For $t>0$, set
\begin{equation}
	\log^{+}{t}=\max\{\log{t},0\},\; 
	\log^{-}{t}:=\log^{+}(1/t)=\max\{-\log{t},0\}.
	\label{eq:3.3.4}
\end{equation}
Then
\begin{equation}
	|\log{t}|=\log^{+}{t}+\log^{-}{t}.
	\label{eq:3.3.5}
\end{equation}

Let $U$ be the open subset as assumed in the proposition.
Then for any 
nonzero holomorphic section $s_{p}\in H^{0}(X,L^{p}\otimes E)$, we have that $\big|\log|s_{p}|_{h_{p}}\big|$ is integrable on 
$\overline{U}$ with respect to $\mathrm{dV}$. We now start with showing that for any $\delta >0$,
\begin{equation}
	\P\left(
	\int_{U}\log^{+}{|\bb{S}_{f,p}(x)|_{h_{p}}}\,\mathrm{dV}(x)\geq 
	\delta p\right)\leq e^{-C_{U,f,\delta}p^{n+1}}.
	\label{eq:3.3.6}
\end{equation}
For this purpose, observe that on $U$ we have
\begin{equation}
	\log^{+}{|\bb{S}_{f,p}|_{h_{p}}}\leq 
	\big|\log\mathcal{M}^{U}_{p}(\bb{S}_{f,p})\big|,
	\label{eq:3.3.7}
\end{equation}
which then supplies us with
\begin{equation}
		\P\left(
		\int_{U}\log^{+}{|\bb{S}_{f,p}(x)|_{h_{p}}}\,\mathrm{dV}(x)\geq 
	\delta p\right)
	\leq \P\left(
		\left|\log\mathcal{M}^{U}_{p}(\bb{S}_{f,p})\right|\geq 
		\frac{\delta}{\Vol(U)} p\right),
	\label{eq:3.3.9}
\end{equation}
where $\Vol(U)$ denotes the volume of $U$ with respect to 
$\mathrm{dV}$. 

Note that we assume that $f$ is of class $\mathscr{C}^{m+1}$ with $m=m(U')$ almost everywhere near $U'$ and with 
$$U\cap\, \esupp\, f\neq\varnothing\, ,$$
the set $\esupp\, f$, by its definition, can not be null measure in any given open subset, then there exists a small open ball $\mathbb{B}\subset \esupp\, f \cap U$ such that $f$ is of class $\mathscr{C}^{m(B)+1}$ (since $m(B)\leq m(U')$) and never vanishing on $\mathbb{B}$. Therefore, Corollary \ref{cor:3.18ss} applies to this open subset $\mathbb{B}$ of $U$ and the function $f$. Then combining \eqref{eq:3.3.9} with Corollary \ref{cor:3.18ss}, we obtain \eqref{eq:3.3.6} with any given $\delta>0$.

The next step is to study 
$\int_{U}\log^{-}{|\bb{S}_{f,p}(x)|_{h_{p}}}\,
\mathrm{dV}(x)$. At first, let us construct a particular 
family of holomorphic local charts to cover $U$. 
Since $U$ is relatively compact in $U'$, 
we apply Proposition \ref{prop:4.3jan} for 
$K=\overline{U}$ to choose the constants $(t_K, c_K)$ 
such that the $t_K$-neighbourhood of $U$ is still included 
in $U'$ and the line bundles $L$, $E$ can be locally 
trivialized on any small (geodesic) balls of radius $2t_K$ on $U$.

For each point $x\in \overline{U}$, let $(V_x, \phi_x)$
be the local chart as in Proposition \ref{prop:4.3jan}. Set 
\begin{equation}
A_x:=\phi_x^{-1}\left(\{\frac{1}{2} 
\frac{t_K}{c_K}<|z|<\frac{3}{4}\frac{t_K}{c_K}\}\right)\subset U'.
\end{equation}
Then the family $\{A_x\}_{x\in\overline{U}}$ 
is an open covering of $\overline{U}$, 
hence there exists a finite number of points 
$\{x_j\}_{j=1}^N\subset\overline{U}$ 
such that $\{A_{x_j}\}_{j=1}^N$ 
is already an open covering of $\overline{U}$. 
We also set
$$A_x\subset B_x:=\phi_x^{-1}\left(\left\{|z|<
\frac{t_K}{c_K}\right\}\right) \subset V_x.$$
By \eqref{eq:4.6Feb}, we have
\begin{equation}
r_{\C^n}(f\circ\phi_x^{-1},\phi_x(B_x))
\leq c_K r(f,V_x)\leq c_K r(f,U').
\end{equation}
For ${t_K}/{c_K}>\delta^{2n+2} > c_K r(f,U')$, 
any open coordinate ball (of $\C^n$) in $\{|z|<{t_K}/{c_K}\}$
of radius $\delta^{2n+2}$ intersects the 
essential support of $f\circ\phi_x^{-1}$.

Next we work on each $A_{x_j}$. 
In order to simplify the notation, after rescaling 
the coordinates of the local chart $(V_{x_j},\phi_{x_j})$, 
we may and we will consider the following model case: 
in the local coordinate $z=(z_1,\ldots, z_n)$, 
set $A=\mathbb{B}(2,3):=\{2<|z|<3\}$, and 
$B=\mathbb{B}(4):=\{|z|<4\}$. 
Assume that on the coordinate open ball 
$\mathbb{B}(5)\subset U'$, the holomorphic line bundles 
$L$ and $E$ can be trivialized by holomorphic local frames. 
In the following, we work on this coordinates $(z_1,\ldots,z_n)$
and let $r_{\C^n}(f,V)$ denote the number given as in 
\eqref{eq:4.3Jan2024} with respect 
to the standard Riemannian metric $g^{T\R^{2n}}$ on $\C^n$.

For $f\in \mathcal{L}^{\infty}_{\mathrm{c}}(X,\R_{\geq 0})$
which is of $\mathscr{C}^{m(U')+1}$ almost everywhere near $B$ 
with $r_{\C^n}(f,B)\leq \frac{1}{2}$, 
we always have $A\,\cap\,\esupp\, f \neq \varnothing$.

Let $e_{L}$, $e_E$ the be respective holomorphic local frames of the line bundles $L$, $E$ 
restricting to the open ball $B=\mathbb{B}(4)$. Set $\alpha_L(z)= 
\log|e_{L}(z)|^{2}_{h_L}$, $\alpha_E(z)= 
\log|e_{E}(z)|^{2}_{h_E}$. On this local chart, we can write
\begin{equation}
	\bb{S}_{f,p}|_{B}=F_{p}(z)e^{\otimes p}_{L}\otimes e_E,
	\label{eq:3.4.11DLM}
\end{equation}
where $F_{p}$ is a random holomorphic function on $B$. Then on $B$
\begin{equation}
	\log|\bb{S}_{f,p}|_{h_{p}}=\log|F_{p}|+\frac{p}{2}\alpha_L+\frac{1}{2}\alpha_E.
	\label{eq:3.4.12DLM}
\end{equation}

Fixing $\varepsilon>0$, we take
\begin{equation}
K_1:=2\varepsilon+\frac{1}{2}\sup_{z\in \mathbb{B}(4)}|\alpha_L(z)|.
\end{equation}
Then by \eqref{eq:3.3.6} and \eqref{eq:3.4.12DLM}, 
we get that, except the event 
$\big\{|\log\mathcal{M}^{\overline{A}}_{p}(\bb{S}_{f,p})|
\geq \varepsilon p\big\}$ of probability less than
$e^{-C_{A,f,\varepsilon}\, p^{n+1}}$ , we have for all $r\in [1,3]$, 
		\begin{equation}
		\int_{\{|z|=r\}} \log^{+}{|F_{p}|}(z)\,\mathrm{d}\sigma_r(z)\leq
	K_{1} p,
	\label{eq:3.4.13DLM}
		\end{equation}
		where $\mathrm{d}\sigma_r$ denote the invariant probability measure on the sphere $\{|z|=r\}$.

For $0<r<4$, consider the Poisson kernel for the ball $\mathbb{B}(r)$,
\begin{equation}
\mathscr{P}_r(\xi,z):=r^{2n-2}\frac{r^2-|\xi|^2}{|\xi-z|^{2n}},
\end{equation}
where $|z|=r$, $\xi\in \mathbb{B}(r)$. Since $\log|F_p|$ is a subharmonic function on $B$, by the sub-mean inequality in terms of Poisson kernel, we conclude that for any point $\xi\in \mathbb{B}(1)$, $1<r<3$,
\begin{equation}
\int_{\{|z|=r\}}\mathscr{P}_r(\xi, z)\big|\log|F_p|(z)\big|\mathrm{d}\sigma_r(z)\leq 2\int_{\{|z|=r\}}\mathscr{P}_r(\xi, z)\log^+|F_p|(z)\mathrm{d}\sigma_r(z)-\log|F_p|(\xi).
\label{eq:4.57nov}
\end{equation}
As a consequence, we have, for $r\in[2,3]$,
\begin{equation}
\begin{split}
&\int_{\{|z|=r\}}\big|\log|F_p|(z)\big|\mathrm{d}\sigma_r(z)\\
&\leq 12\left(\frac{9}{2}\right)^{2n-1}\int_{\{|z|=r\}}\log^+|F_p|(z)\mathrm{d}\sigma(z)-2\left(\frac{3}{2}\right)^{2n-1}\log|F_p|(\xi).
\end{split}
\label{eq:4.55Jan2024}
\end{equation}

Note that $\mathbb{B}(1)\subset B$, following from the condition $r_{\C^n}(f,B)\leq \frac{1}{2}$, Corollary \ref{cor:3.18ss} applies to $\mathbb{B}(1)$. We conclude that except for an event of probability less than $e^{-C_{\mathbb{B}(1),f,1}p^{n+1}}$, we can always find a point $\xi_0\in \overline{\mathbb{B}(1)}$ with the property:
\begin{equation}
\log|F_p|(\xi_0)>-p.
\end{equation}
Then combining it with \eqref{eq:3.4.13DLM} and \eqref{eq:4.55Jan2024}, we get that except for an event $\Omega_1$ of probability less than $e^{-C_{B,f,K_2}p^{n+1}}$, we have that for all $r\in[2,3]$,
\begin{equation}
\int_{\{|z|=r\}}\big|\log|F_p|(z)\big|\mathrm{d}\sigma_r(z)\leq K_2 p,
	\label{eq:3.4.14DLM}
\end{equation}
where $K_2>0$ is a sufficiently large constant given by
$$K_2=12\left(\frac{9}{2}\right)^{2n-1} K_1+2\left(\frac{3}{2}\right)^{2n-1}.$$

Given $\delta\in(0,\frac{1}{2}]$ (the case of $\delta > \frac{1}{2}$
will become a consequence of this one), we fix a decomposition of the 
unit sphere $\{|z|=1\}$ into a disjoint union of sets $I_1$, $I_2$, $\cdots$, 
$I_q$ with (Euclidean) diameter $\simeq \delta^{2n+2}$. 
The number $q$ is a large integer depending on $\delta$. 
The following inequality was established in \cite[(4.4)]{SZZ}: there exists a constant $C_n>0$ independent of $\delta\in\; 
(0,\frac{1}{2}]$ such that for any $r\in(1,3]\,$, any collection of points $\{\xi_k\}_{k=1}^q \subset \mathbb{B}(4)$ with $\mathrm{dist}(\xi_k, (r-\delta)I_k)< \delta^{2n+2}$ and any subharmonic function $u$ on $\mathbb{B}(4)$,
\begin{equation}
\int_{\{|z|=r\}}u(z)\mathrm{d}\sigma_r(z)\geq \sum_{j=1}^q \sigma_1(I_k)u(\xi_k)-C_n\delta \int_{\{|z|=r\}}|u(z)|\mathrm{d}\sigma_r(z).
\label{eq:4.59nov}
\end{equation}
The main idea to prove \eqref{eq:4.59nov} in \cite[Section 4.1]{SZZ}
is to use the Poisson kernel as in \eqref{eq:4.57nov} for each
term $u(\xi_k)$, and $C_n$ is a constant related to the estimate
on the derivatives of $\mathscr{P}_r$, and in our setting, 
we can take $C_n=36^{n}(4n+4)$.

From this point we proceed as in \cite[Section 
4.1, pp.\ 1992]{SZZ} but with some necessary modifications, and we will also produce a rough formula for the constant $r(U,U')$ needed in our proposition. Fix $\delta\in\; 
(0,\frac{1}{2}]$, and let $f\in \mathcal{L}^{\infty}_{\mathrm{c}}(X,\R_{\geq 0})$ which is of $\mathscr{C}^{m(U')+1}$ almost everywhere near $B$ such that $r_{\C^n}(f,B)\leq\frac{2}{3}\delta^{2n+2}$.

Set $Q:=\lceil \frac{9}{4}\delta^{-2n-2} \rceil$ and set $r_j=\frac{3}{2}+\frac{2}{3}\delta^{2n+2}j$, for $0\leq j\leq Q$. Then for $1\leq k \leq q$, $0\leq j\leq Q$, set
\begin{equation}
V_{kj}:=\{z\in B\;:\; \mathrm{dist}^{\C^n}(z,r_jI_k)< \frac{1}{2}\delta^{2n+2}\}.
\end{equation}
Then $\{V_{kj}\}_{k,j}$ form a covering for $\overline{A}$, moreover, each $V_{kj}$ always contains a coordinate open ball of raidus $\delta^{2n+2}$, so that Corollary \ref{cor:3.18ss} applies: there exists $C_{kj}>0$ such that for all $p\gg 0$
	\begin{equation}
		\P\left(\left|\log{\mathcal{M}^{V_{kj}}_{p}(\bb{S}_{f,p})}\right|> \delta p\right) \leq e^{-C_{kj}p^{n+1}}.
		\label{eq:4.60ppp}
	\end{equation}
	Set the event $\Omega_\delta:=\cup_{k,j} \left\{\left|\log{\mathcal{M}^{V_{kj}}_{p}(\bb{S}_{f,p})}\right|\geq \delta p\right\}$, then we have that for $p\gg 0$,
		\begin{equation}
		\P\left(\Omega_\delta\right) \leq e^{-C_{B,f,\delta}p^{n+1}}.
	\end{equation}

When $\bb{S}_{f,p}$ lies outside $\Omega_\delta$, we can always find $\xi_{kj}\in V_{kj}$ with 
\begin{equation}
\log |\bb{S}_{f,p}|_{h_p}(\xi_{kj})\geq -\delta p.
\end{equation}
Fix $r\in\,]2,3[\,$, there exists $j\in \{0,1,\ldots,Q\}$ such that 
$$|r-\delta-r_j|\leq \frac{2}{5}\delta^{2n+2}.$$ 
Then $\mathrm{dist}^{\C^n}(rI_k, \xi_{kj})<2\delta$, 
so that we have
\begin{equation}
\int_{\{|z|=r\}}\alpha_L(z)\mathrm{d}\sigma_r(z)\geq 
\sum_{j=1}^q \sigma_1(I_k)\alpha_L(\xi_{kj})-2\delta 
\sup_{z\in B} |d\alpha_L(z)|.
\end{equation}
A similar inequality holds also for $\alpha_E$.
Moreover, we have $\mathrm{dist}^{\C^n}((r-\delta)I_k,
\xi_{kj})<\delta^{2n+2}$, then we can apply \eqref{eq:4.59nov} 
for $\log |F_p|$.
We combine this with \eqref{eq:3.4.14DLM}, so 
when $\bb{S}_{f,p}$ lies outside $\Omega_1\cup \Omega_\delta$, we have
\begin{equation}
	\begin{split}
		&-\int_{\{|z|=r\}} \log|\bb{S}_{f,p}|_{h_{p}}(z)\mathrm{d}\sigma_r(z)\\
		&\qquad\leq 
		-\sum_{k=1}^{q} \sigma_r(I_k) 
		\log|\bb{S}_{f,p}|_{h_{p}}(\xi_{jk})+C_n\delta\int_{\{|z|=r\}} 
		\big|\log{|F_{p}|}\big|(z)\,\mathrm{d}\sigma_r(z)\\
		&\qquad\qquad+p\delta \sup_{z\in 
		B}|d\alpha_L(z)| + \delta \sup_{z\in 
		B}|d\alpha_E(z)|\\
		&\qquad\leq \delta p+C_n K_2 \delta p+p\delta \sup_{z\in 
		B}|d\alpha_L(z)| + \delta \sup_{z\in 
		B}|d\alpha_E(z)|.
	\end{split}
	\label{eq:4.65Jan24}
\end{equation}
We set
\begin{equation}
K_3:=1+C_n K_2+ \sup_{z\in 
		B}|d\alpha_L(z)| + \varepsilon >0.
\end{equation}
Let $b_0>0$ be such that for all $z\in A$, $\mathrm{dV}(z)\leq b_0 \mathrm{d}\lambda(z)$, where $\mathrm{d}\lambda(z)$ denotes the standard Lebesgue measure on $\C^n$. Then taking the integral with repsect to $r\in[2,3]$ and rescaling $\mathrm{d}\sigma_r$ to the standard spherical measure for the sphere of radius $r$, we conclude, by taking advantage of \eqref{eq:4.65Jan24}, that
\begin{equation}
	\P\left(
	-\int_{A}\log{|\bb{S}_{f,p}|_{h_{p}}}\,\mathrm{dV}\geq 
	b_0 K_3 \frac{\pi^n}{n!}(9^n-4^n)\delta p\right)\leq e^{-C_{B,f,\delta}p^{n+1}},\;\forall\;p\gg 0.
	\label{eq:3.4.16DLM}
\end{equation}
Combining the above estimate with \eqref{eq:3.3.6}, we get that for $p\gg 0$,
\begin{equation}
	\P\left(
	\int_{A}\Big|\log{|\bb{S}_{f,p}|_{h_{p}}}(z)\Big|\,\mathrm{dV}(z)> 
	b_0 K_3 \frac{\pi^n}{n!}(9^n-4^n)\delta p\right)\leq e^{-C'_{B,f,\delta}p^{n+1}}.
	\label{eq:4.67DLM3}
\end{equation}

Now we reformulate \eqref{eq:4.67DLM3} for $A_{x_j}$. By the choice of $(t_K,c_K)$ from Proposition \ref{prop:4.3jan} for $K=\overline{U}$, we can identify $A_{x_j}$ with $A$ and $B_{x_j}$ with $B$ via rescaling the coordinate by $\frac{4c_K}{t_K}$. Then we can take $b_0=\left(\frac{t_K}{4}\right)^{2n}$. Therefore, we get, for $p\gg 0$,
\begin{equation}
	\P\left(
	\int_{U}\Big|\log{|\bb{S}_{f,p}|_{h_{p}}}(z)\Big|\,\mathrm{dV}(z)> N
	b_0 K_3 \frac{\pi^n}{n!}(9^n-4^n)\delta p\right)\leq e^{-C^{\prime\prime}_{B,f,\delta}p^{n+1}}.
	\label{eq:4.67DLM4}
\end{equation}

The condition for \eqref{eq:4.67DLM4} to hold is that 
$r_{\C^n}(f,B)\leq \frac{2}{3}\delta^{2n+2}$
(with $\delta\leq \frac{1}{2}$). A sufficient condition for this is
\begin{equation}
r(f,U')\leq \min\left\{\frac{t_K}{6c_K^2}\delta^{2n+2},
\frac{t_K}{8c_K^2}\right\}.
\end{equation}
Finally, we set
\begin{equation}
r(U,U'):=\frac{t_K}{6c_K^2}\min\left\{3\cdot 4^n, \left(\frac{1}{N
	b_0 K_3 \frac{\pi^n}{n!}(9^n-4^n)}\right)^{2n+2}\right\}>0,
	\label{eq:4.70Jan}
\end{equation}
which fulfills our purpose as stated in Proposition \ref{thm:4.3Jan},
in particular, it implies that $U\cap\, \esupp\, f\neq\varnothing\,$.
This completes the proof.
\end{proof}

\subsection{Expected mass of random $\cLL$-holomorphic sections}\label{ss:mass}

 Recall that the volume form that is used to define the
 $\cLL$-inner products is $\mathrm{dV}:=\Theta^n/n!$. 
Similarly, since $c_1(L,h_L)$ is uniformly positive on $X$
we define a different volume form on $X$,
\begin{equation}
\mathrm{dV}^L:=\frac{c_1(L,h_L)^n}{n!},
\quad\mathrm{dV}^L(z)=\bb{b}_0(z) \mathrm{dV}(z),
\label{eq:volumeL}
\end{equation}
where the positive function 
$\bb{b}_0(z)=\det(\dot{R}_z^{L}/2\pi)$ on $X$ 
is given by \eqref{eq:b0}.

Recall that for a holomorphic section 
$s_p\in H^0(X,L^p\otimes E)$ we introduced 
the measure $M_p(s_p)$ on $X$
in Definition \ref{def:4.10mass}. 
The following proposition gives preliminary results on 
the expectation of mass distribution 
$M_p(\bb{S}_{f,p})$ of $\bb{S}_{f,p}$.
The proof follows from the asymptotic expansion
of $T^2_{f,p}$ given in Section \ref{ss2nov}.
\begin{Proposition}\label{prop:4.18p1}
Let $(X,J,\Theta)$ be a connected Hermitian complex manifold
and let $(L,h_{L})$, $(E,h_{E})$ be two holomorphic line bundles
on $X$ with smooth Hermitian metrics. 
Assume Condition \ref{condt:main} to hold. 
Fix a nontrivial 
$f\in\mathscr{C}^1_{\mathrm{c}}(X,\R_{\geq 0})$, or 
$f=\bb{1}_A$ for a relative compact domain $A\subset X$ 
with a continuous boundary, and let 
$\{\bb{S}_{f,p}\in H^0_{(2)}(X,L^p\otimes E)\}_p$ 
be the associated random $\cLL$-holomorphic sections defined 
via $\{T_{f,p}\}_p$. Then for any $z\in X$, we have
\begin{equation}
\E[|\bb{S}_{f,p}(z)|^2_{h_p}]=T^2_{f,p}(z,z).
\label{eq:5.42mass}
\end{equation}
Hence we have the locally uniform weak-star convergence of measures on $X$,
\begin{equation}
\E[M_p(\bb{S}_{f,p})]\rightarrow f^2\mathrm{dV}^L.
\label{eq:5.43mass}
\end{equation}

\end{Proposition}
\begin{proof}
The identity \eqref{eq:5.42mass} follows from the series
formula for $\bb{S}_{f,p}$ given in \eqref{eq:4.34series}
(see also \eqref{eq:3.12nov}). 
Then for any function $g\in\mathscr{C}^0_{\mathrm{c}}(X)$, we have
\begin{equation}
\E\left[\int_X g(z)M_p(\bb{S}_{f,p})(z)\right]=
\frac{1}{p^n}\int_X g(z)T^2_{f,p}(z,z)\mathrm{dV}(z).
\label{eq:5.44mass}
\end{equation}
When $f\in\mathscr{C}^1_{\mathrm{c}}(X,\R_{\geq 0})$, by Theorem \ref{thm:2.5Novfg}, we have 
\begin{equation}
T^2_{f,p}(z,z)=p^n f(z)^2 \bb{b}_0(z)+\mathcal{O}(p^{n-1/2}).
\end{equation}
We conclude \eqref{eq:5.43mass} from \eqref{eq:5.44mass}.

Now we consider the case $f=\bb{1}_A$. Fix a function $g\in\mathscr{C}^0_{\mathrm{c}}(X)$ and let $K$ be a compact subset with $\mathrm{supp}\; g\subset K$. By \eqref{eq:June2.38on}, we have for $x\in K$
\begin{equation}
\frac{1}{p^n}T^2_{f,p}(x,x)\leq C_K,
\end{equation}
with some constant $C_K\geq 1$.
For any $\varepsilon>0$, let $\delta>0$ be a small number such that 
\begin{equation}
\mathrm{Vol}(\{z\in X \;:\; \mathrm{dist}(z,\partial A)\leq \delta\})\leq \frac{\varepsilon}{3 C_K}.
\end{equation}
By the uniform expansion of $T^2_{f,p}(x,x)$ Theorem \ref{thm:2.5Novfg}
for the points away from $\partial A$, we get that there exits
$p_0\in \N$ such that for all $p\geq p_0$
\begin{equation}
\left|\frac{1}{p^n}\int_{z\in X, \mathrm{dist}(z,\partial A)\geq 
\delta} g(z)T^2_{f,p}(z,z)\mathrm{dV}(z)-
\int_{A(\delta)}g(z)\mathrm{dV}^L(z)\right|\leq 
\|g\|_{\mathcal{L}^\infty}\frac{\varepsilon}{3},
\end{equation}
where $A(\delta):=\{z\in A\;:\; \mathrm{dist}(z,\partial A)\geq \delta\}$.

Assembling the above estimates together we have for $p\geq p_0$,
\begin{equation}
\left| \frac{1}{p^n}\int_X g(z)T^2_{f,p}(z,z)\mathrm{dV}(z)- \int_A g(z)\mathrm{dV}^L(z)\right| \leq \|g\|_{\mathcal{L}^\infty}\varepsilon.
\end{equation}
Then we conclude that
\begin{equation}
\frac{1}{p^n}\int_X g(z)T^2_{f,p}(z,z)\mathrm{dV}(z)\rightarrow\int_A g(z)\mathrm{dV}^L(z),\quad p\to\infty.
\end{equation}
The proof is complete.
\end{proof}

We need the following lemma to prove Proposition \ref{prop:4.19p2}.
\begin{Lemma}\label{lm:massvar}
With the same assumptions as in Proposition \ref{prop:4.18p1}
for the geometric data, let $U$ be a relative compact open subset 
of $X$, and fix a nontrivial 
$f\in\mathscr{C}^{m(U)+1}_{\mathrm{c}}(U,\R_{\geq 0})$. 
Then for any $g\in\mathscr{C}^0_{\mathrm{c}}(U)$, set $Y_p^{g}
=\int_X g(z)M_{p}(\bb{S}_{f,p})(z)$, we have
\begin{equation}
\mathrm{Var}[Y_p^{g}]:=\E[|Y_p^{g}|^2]-|\E[Y_p^{g}]|^2
=\int_X |g(z)|^2 f(z)^4\mathrm{dV}^L(z)+o(1).
\label{eq:4.92mass}
\end{equation}
\end{Lemma}
\begin{proof}
Note that we have
\begin{equation}
\begin{split}
\E[|Y_p^{g}|^2]&=\frac{1}{p^{2n}}\int_{X\times X} g(z)\overline{g(w)}\E[|\bb{S}_{f,p}(z)|^2_{h_p}|\bb{S}_{f,p}(w)|^2_{h_p}]\mathrm{dV}(z)\mathrm{dV}(w)\\
&=\frac{1}{p^{2n}}\int_{X\times X} g(z)\overline{g(w)}T^2_{f,p}(z,z)T^2_{f,p}(w,w)\E\left[\frac{|\bb{S}_{f,p}(z)|^2_{h_p}}{T^2_{f,p}(z,z)}\frac{|\bb{S}_{f,p}(w)|^2_{h_p}}{T^2_{f,p}(w,w)}\right]\mathrm{dV}(z)\mathrm{dV}(w).
\end{split}
\label{eq:4.93mass}
\end{equation}
By \cite[Lemma 5.2]{MR3742811}, we have
\begin{equation}
\E\left[\frac{|\bb{S}_{f,p}(z)|^2_{h_p}}{T^2_{f,p}(z,z)}\frac{|\bb{S}_{f,p}(w)|^2_{h_p}}{T^2_{f,p}(w,w)}\right]=1+N_{f,p}(z,w)^2.
\label{eq:4.94mass}
\end{equation}
Since $f$ is assumed to be $\mathscr{C}^{m(U)+1}$, then $N_{f,p}(z,w)^2$ satisfies the asymptotics \eqref{eq:5.2.4Nov}, as a consequence, we have (cf. \eqref{eq:6.21feb24}--\eqref{eq:6.23feb24})
\begin{equation}
\begin{split}
\E[|Y_p^{g}|^2]=\left|\E[Y_p^{g}]\right|^2+\int_X |g(z)|^2|f(z)|^4\mathrm{dV}^L(z)+o(1).
\end{split}
\end{equation}
This proof is complete.
\end{proof}

\begin{proof}[Proof of Proposition \ref{prop:4.19p2}]
By Proposition \ref{prop:4.18p1}, we have
\begin{equation}
\E[Y^g_p]=\int_X g(z)f(z)^2\mathrm{dV}^L(z)+o(1).
\end{equation}
Then as $N\rightarrow+\infty$, we have
\begin{equation}
\frac{1}{N}\sum_{1\leq p \leq N} \E[Y^g_p]\rightarrow 
\int_X g(z)f(z)^2\mathrm{dV}^L(z).
\label{eq:4.97mass}
\end{equation}
Note that $\{Y^g_p\}_p$ is a sequence of independent variables, 
and by Lemma \ref{lm:massvar}, we have $\mathrm{Var}[Y^g_p]=
\mathcal{O}(1)$, so that by the strong law of large numbers for pairwise independent random variables
(see \cite[Theorem 1]{MR722846}), we have
\begin{equation}
\frac{1}{N}\sum_{1\leq p\leq N}\int_X g(z)M_{p}(\bb{S}_{f,p})(z)-\frac{1}{N}\sum_{1\leq p \leq N} \E[Y^g_p]\rightarrow 0, \text{almost\,\,surely.}
\end{equation}
Combining the above convergence with \eqref{eq:4.97mass}
we conclude the proof.
\end{proof}

When $X$ is compact, Bayraktar \cite{MR4159385}
considered the sub-Gaussian holomorphic sections and proved
the quantum ergodicity for their mass distribution
by using the Hanson-Wright inequality.
This result is stronger than Proposition \ref{prop:4.19p2}.

\section{Expectation and equidistribution on the support
of the symbol}
\label{section5:onsupport}
We continue our discussion of zeros of $\bb{S}_{f,p}$ in Section \ref{ss4Jan},
especially, we prove that the random zeros will be asymptotically uniformly distributed
with respect to $c_1(L,h_L)$ on parts inside the support of $f$ where it is
$\mathscr{C}^{m(U)+1}$ with $\supp f \subset U\Subset X$.
These results are applications of Theorem \ref{thm:4.6Jan2024} specializing to the case $r(f,U')=0$.
\subsection{Equidistribution of random zeros on the support; 
proof of Theorem \ref{thm:6.3}}\label{Section3.1nov}

In \cite[Section 5]{DrLM:2023aa}, the zeros of random sections $\bb{S}_{f,p}$
were studied for the case $0\neq f\in \mathscr{C}^{\infty}_{\mathrm{c}}(X,\R)$, 
and the equidistribution results were proven for their zeros as $p\rightarrow \infty$
on the support of $f$ with the assumption that $f$ vanishes up to order $2$.
Here, we remove the condition on the vanishing orders of $f$ for the equidistribution results.

\begin{proof}[Proof of Theorem \ref{thm:6.2}]
 By our assumption $U\subset \mathrm{ess.supp\,} f$, we get $r(f,U)=0$, in 
 Theorem \ref{thm:4.6Jan2024}, we have $\delta_0(f)=0$, thus \eqref{eq:6.0.7} 
 holds for all $\delta>0$. Similarly, taking $\delta_0(f)=0$ in \eqref{eq:4.73Jan}, 
 we get \eqref{eq:5.2Jan24}.
\end{proof}

\begin{proof}[Proof of Theorem \ref{thm:6.3}]
Note that the estimate \eqref{eq:1.6.14paris} is a direct consequence of 
\eqref{eq:6.1.10} by taking 
$\delta=n\Vol^{L}_{2n}(U')$. 
Hence, it is sufficient to prove \eqref{eq:6.1.10}. For this purpose, let 
$\bb{1}_{U'}$ denote the indicator function of $U'$ on $X$.
Let $\delta>0$ be arbitrary, and take $\psi_{1}$, $\psi_{2}\in 
\mathcal{C}^{\infty}_{0}(X,\R)$ with supports in $U$ such that
$0\leq \psi_{1}\leq\bb{1}_{U'}\leq \psi_{2}\leq 1,$ and
\begin{equation}
		\int_{X} \psi_{1}\frac{c_{1}(L,h_{L})^{n}}{n!}\geq 
		\Vol^{L}_{2n}(U')-\delta, \quad
		\int_{X} \psi_{2}\frac{c_{1}(L,h_{L})^{n}}{n!}\leq 
		\Vol^{L}_{2n}(U')+\delta\,.
	\label{eq:5.1.30paris}
\end{equation}
Note that the existence of such functions is guaranteed by the 
assumption that $\partial U'$ has measure $0$ with respect to 
$\mathrm{dV}$, hence also to $\frac{1}{n!}c_{1}(L,h_{L})^{n}$.
For $j\in \{1,2\},$ set $\varphi_{j}=\frac{1}{(n-1)!}{\psi_{j}}c_{1}(L,h_{L})^{n-1}$. 
By applying Theorem \ref{thm:6.2} to $\varphi_{j}$ separately, we 
get exactly \eqref{eq:6.1.10}. 
\end{proof}

Another interesting object that is closely related to this equidistribution result
is to study the asymptotic behavior of $\frac{1}{p}\mathbb{E}[[\Div(\bb{S}_{f,p})]]$
and compare it with $c_1(L,h_L)$, the expected limit. By Lemma \ref{lem:currentf}
and Theorem \ref{thm:3.6jan}, we have
\begin{equation}
\frac{1}{p}\mathbb{E}[[\Div(\bb{S}_{f,p})]]-c_1(L,h_L)=
\frac{1}{p} c_1(E,h_E)+\frac{\sqrt{-1}}{2\pi p} \partial\overline{\partial} \log {T_{f,p}^2(x,x)}.
\label{eq:5.8jan24}
\end{equation}

Now we show that, under a global finiteness on the geometry of the manifold, 
the convergence of $\frac{1}{p}\mathbb{E}[[\Div(\bb{S}_{f,p})]]$ to $c_1(L,h_L)$
on $\esupp\, f$ can be obtained as a consequence of the equidistribution 
result \eqref{eq:5.2Jan24}. Certainly, this is a weaker version of the convergence 
given in Theorem \ref{thm:6.2equi} via the method of pluripotential theory on $X$.
\begin{Proposition}\label{thm:5.5finite}
We assume the same geometric conditions on $X, L, E$ as in Theorem \ref{thm:6.2}. 
Furthermore, we assume one of the following conditions to hold:
\begin{itemize}
\item[(i)] $\Theta$ is K\"{a}hler and $\int_X  c_1(L,h_L) \wedge \Theta^{n-1}<\infty$, or
\item[(ii)] $\int_X c_1(L,h_L)^{n}<\infty$.
\end{itemize}
Fix $0\neq f\in \mathcal{L}^{\infty}_{\mathrm{c}}(X,\R)$. 
Let $U$ be an open subset of $X$ such that $U\subset \mathrm{ess.supp\,} f$, 
and we assume that $f$ is of class $\mathscr{C}^{m(U)+1}$ almost everywhere on $U$. 
Then we have
	\begin{equation}
	\begin{split}
	&\frac{1}{p}\mathbb{E}[[\Div(\bb{S}_{f,p})]]|_U \rightarrow c_1(L,h_L)|_U;\\
&\frac{\sqrt{-1}}{2\pi p} \partial\overline{\partial} \log {T_{f,p}^2(x,x)}|_U\rightarrow 0,
\end{split}
		\label{eq:5.11finite}
	\end{equation}
where the limits in \eqref{eq:5.11finite} are taken with respect 
to the convergence of $(1,1)$-currents on $U$.
\end{Proposition}
\begin{proof}
Applying Theorem \ref{thm:6.2}, we get that for any $\varphi\in 
	\Omega^{n-1,n-1}_{\mathrm{c}}(U)$, we have
	\begin{equation}
		\P\Big( \ 
		\lim_{p\rightarrow \infty}\big\langle\frac{1}{p}[\Div(\bb{S}_{f,p})],
		\varphi\big\rangle=\big\langle c_{1}(L,h_L),\varphi\big\rangle\, \Big)=1.
		\label{eq:5.15cpt}
	\end{equation}	

Now we assume that first case: $\Theta$ is K\"{a}hler and 
$\int_X  c_1(L,h_L) \wedge \Theta^{n-1}<\infty$. For all $\varphi\in 
	\Omega^{n-1,n-1}(X)$, $s_{p}\in H^{0}(X,L^{p})$, due to the 
	positivity of the current $[\Div(s_{p})]$, we have
	\begin{equation}
		|\langle[\Div(s_{p})],\varphi\rangle|\leq 
		|\varphi|_{\mathscr{C}^{0}(X)}\langle[\Div(s_{p})],
		\Theta^{n-1}\rangle=p|\varphi|_{\mathscr{C}^{0}(X)}\int_{X}c_{1}(L,h_{L})\wedge\Theta^{n-1}.
		\label{eq:3.25new23s}
	\end{equation}
 and the last identity follows from the Poincar\'{e}-Lelong formula. We 
	may set $C=\int_{X}c_{1}(L,h_{L})\wedge\Theta^{n-1}<\infty$, then
	\begin{equation}
		\frac{1}{p}|\langle[\Div(s_{p})],\varphi\rangle|\leq 
		C|\varphi|_{\mathscr{C}^{0}(X)}.
	\end{equation}
Then applying the dominated convergence theorem to the limit in \eqref{eq:5.15cpt}, we infer that 
	\begin{equation}
\lim_{p\rightarrow \infty}\E\left[\big\langle\frac{1}{p}[\Div(\bb{S}_{f,p})],
\varphi\big\rangle\right]=\big\langle c_{1}(L,h_L),\varphi\big\rangle.
	\end{equation}
	By the definition of $\E\left[[\Div(\bb{S}_{f,p})]\right]$,
	we can rewrite the above limit as
		\begin{equation}
		\lim_{p\rightarrow \infty}\big\langle\frac{1}{p}\E\left[[\Div(\bb{S}_{f,p})]\right],
		\varphi\big\rangle=\big\langle c_{1}(L,h_L),\varphi\big\rangle.
	\end{equation}
	This way, we get the first line of \eqref{eq:5.11finite}.

	By Theorem \ref{thm:3.6jan}, we have
	\begin{equation}
	\frac{1}{p}\E\left[[\Div(\bb{S}_{f,p})]\right]-c_1(L,h_L)=\frac{1}{p}c_1(E,h_E)+
	\frac{\sqrt{-1}}{2\pi p} \partial\overline{\partial} \log {T_{f,p}^2(x,x)}.
	\end{equation}
	Then the second line of \eqref{eq:5.11finite} follows from the first one.
	
For the second case $\int_X c_1(L,h_L)^{n}<\infty$, we just replace $\Theta$ by 
$c_1(L,h_L)$ in the above arguments, and certainly we have to change the 
$\mathscr{C}^{0}$-norm for test forms to the one induced by $c_1(L,h_L)$, 
then we conclude the same conclusion \eqref{eq:5.11finite}. This way, we complete our proof.
\end{proof}

\begin{Remark}
It is an interesting question to ask that if the finiteness assumptions listed in 
Proposition \ref{thm:5.5finite} together with \eqref{eq:mainassumptions} 
would imply that $d_{p}<+\infty$, it seems so for the case of Riemannian surfaces. 
But the authors have no awareness of such results in general.
\end{Remark}

\subsection{Proofs of Theorems \ref{thm:5.5jan24} and \ref{thm:6.2equi}}\label{Section:5.2march}
Now we give a general result on $\mathbb{E}[[\Div(\bb{S}_{f,p})]]$
that is obtained from the pluripotential theory on $X$. 

Let us recall some basic facts about the plurisubharmonic functions. 
Let $X$ be a connected complex manifold of complex dimension $n\geq 1$, 
and we fix a smooth K\"{a}hler form $\omega$ on $X$. 
Recall that a quasi-plurisubharmonic (quasi-psh for short) function 
on an open domain of $X$ is locally the difference of a plurisubharmonic function and a smooth one. 
For any open subset $U\subset X$, set $\mathrm{PSH}(U,\omega)$ the space of quasi-psh functions 
on $U$ with respect to $\omega$, that consists of the upper semi-coutinuous functions 
$\varphi: U\rightarrow [-\infty,+\infty[\,$ with the propoerty 
$dd^{\mathrm{c}}\varphi+\omega\geq 0$ as $(1,1)$-current on $U$, 
where $dd^{\mathrm{c}}=\frac{\sqrt{-1}}{\pi}\partial\overline{\partial}$. 

The following two propositions are well-known facts for the functions in 
$\mathrm{PSH}(U,\omega)$, which follow from the corresponding properties of 
plurisubharmonic functions on the domains of $\C^n$. For their proofs, we refer to 
\cite[Theorems 3.2.12, 3.2.13 and 4.1.8]{MR1301332}.
\begin{Proposition}\label{prop1appendix}
Fix an open domain $U$ in $X$, and a compact subset $K\subset U$.
 Let $\mathcal{F}\subset \mathrm{PSH}(U,\omega)$ be a family 
 of quasi-psh functions on $U$ such that for each 
 $\varphi\in \mathcal{F}$, $\max_{z\in K}\varphi(z)=0$, then $\mathcal{F}|_K$ 
 is a bounded subset in $\mathcal{L}^q(K,\R)$ ($\infty > q\geq 1$), 
 and it is also a relatively compact subset of $\mathcal{L}^1(K,\R)$.
\end{Proposition}

\begin{Proposition}\label{prop2appendix}
Let $U$ be a relatively compact open domain $U$ in $X$.
Let $\varphi_\ell$ ($\ell\in\N$) be a sequence of quasi-psh functions 
on $\overline{U}$ with $dd^{\mathrm{c}}\varphi\geq -\omega$ on 
$\overline{U}$, which have a uniform upper bound on $\overline{U}$. 
Then, only one the following situations holds
\begin{itemize}
	\item either $\varphi_\ell$ converge uniformly to $-\infty$ 
	on any compact subset of $U$, or
	\item there exists a subsequence $\varphi_{\ell_j}$ converges 
	in $\mathcal{L}^q$-norm ($\infty > q\geq 1$) on $U$ to a quasi-psh 
	function $\varphi$ on $U$ with $dd^{\mathrm{c}}\varphi\geq -\omega$.  
\end{itemize}
\end{Proposition}

The above propositions allow us to establish some estimates on the 
integrals of Berezin--Toeplitz kernel function $T^2_{f,p}(x,x)$, 
and then we can deduce the properties on $\mathbb{E}[[\Div(\bb{S}_{f,p})]]$. 
At first, we give the proof of Theorem \ref{thm:5.5jan24}.

\begin{proof}[Proof of Theorem \ref{thm:5.5jan24}]
With our assumption on $f$, by Lemma \ref{lemma:4.2jan}, we get $T_{f,p}\neq 0$ 
for all sufficiently large $p$, then we have the corresponding random 
$\cLL$-holomorphic section $\bb{S}_{f,p}$ which is not identically zero. 
By Theorem \ref{thm:3.6jan}, we have the identity of $(1,1)$-currents on $X$ for $p\gg 0$:
\begin{equation}
 0\leq \frac{1}{p}\E[[\mathrm{Div}(\bb{S}_{f,p})]]=
 \frac{\sqrt{-1}}{2\pi p}\partial\overline{\partial} \log T^2_{f,p}(x,x) + 
 c_1(L,h_L)+\frac{1}{p}c_1(E,h_E).
 \label{eq:5.22jan}
\end{equation}

Let $C_E>0$ be a constant such that $c_1(E,h_E)\leq C_E c_1(L,h_L)$ on 
$\overline{U}$. Then for all $p\geq \lceil C_E \rceil$, we have
\begin{equation}
c_1(L,h_L)+\frac{1}{p}c_1(E,h_E)\leq 2 c_1(L,h_L),\;\text{on}\; \overline{U}.
\label{eq:5.23jan}
\end{equation} 

Set 
\begin{equation}
u_p(x):=\frac{1}{p}\log T^2_{f,p}(x,x),
\label{eq:fctup}
\end{equation}
 and we work on the connected compact subset $K=\overline{U}$. 
 Then by \eqref{eq:5.22jan} and \eqref{eq:5.23jan},  
 $u_p$ is a quasi-psh function with respect to the K\"{a}hler form 
 $4 c_1(L,h_L)$ on $K$, that is, $dd^{\mathrm{c}}u_p+4 c_1(L,h_L)\geq 0$ on $K$.
Set now
\begin{equation}
m_p:=\max_{x\in K} u_p(x).
\end{equation}
Then applying \eqref{eq:June2.44on} on $K$ and Theorem \ref{thm:2.5Novfg} 
for the open subset $U\cap B\subset K$, there exist constants $C' \geq n$, 
$p_0\geq 2$ such that for all $p\geq p_0$,
\begin{equation}
0\leq m_p\leq \frac{C'\log {p}}{p}
\end{equation}
Therefore $\{u_p-m_p\}_{p\gg 0}$ is a family of quasi-psh functions on $K$ 
(actually on an open neighbourhood of $K$) such that each one has the maximum 
$0$ in $K$. By Proposition \ref{prop1appendix}, we conclude that  $\{u_p-m_p\}_{p\gg 0}$ 
is a bounded set in $\mathcal{L}^q(K,\R)$ with any $q\geq 1$. In particular, 
there exists a constant $C>0$ such that for all $p\gg 0$, we have
\begin{equation}
\|u_p-m_p\|_{\mathcal{L}^1(K,\R)}:=\int_K |u_p(x)-m_p|\mathrm{dV}(x)\leq C.
\end{equation}
As a consequence, for all $p\gg 0$, we get
\begin{equation}
\int_K |\log {T^2_{f,p}(x,x)}|\mathrm{dV}(x)\leq Cp+C'\mathrm{Vol}(K)\log{p},
\label{eq:5.19jan24}
\end{equation}
Then for any nonempty open subset $A\subset K$, we have
\begin{equation}
\int_A |\log {T^2_{f,p}(x,x)}|\mathrm{dV}(x)\leq Cp+C'\mathrm{Vol}(A)\log{p},
\end{equation}
We conclude that for all $p\gg 0$,
\begin{equation}
-Cp-C'\mathrm{Vol}(A)\log{p}\leq \int_A
\log {T^2_{f,p}(x,x)}\mathrm{dV}(x)\leq C'\mathrm{Vol}(A)\log{p},
\label{eq:5.29Jan}
\end{equation}
Then considering the probability space modelled as 
$(A,\frac{1}{\mathrm{Vol}(A)}\mathrm{dV}|_A)$ and applying Jensen's 
inequality with $\log(t)$, we get
\begin{equation}
 \frac{1}{\mathrm{Vol}(A)}\int_A \log {T^2_{f,p}(x,x)}\mathrm{dV}(x) \leq 
 \log\left(\frac{1}{\mathrm{Vol}(A)}\int_A T^2_{f,p}(x,x)\mathrm{dV}(x) \right)
 \label{eq:5.31feb24}
\end{equation}
Combining it with \eqref{eq:5.29Jan}, we get
\begin{equation}
\log\left(\frac{1}{\mathrm{Vol}(A)}\int_A T^2_{f,p}(x,x)\mathrm{dV}(x) \right)
\geq -\frac{Cp}{\mathrm{Vol}(A)}-C'\log{p}.
\end{equation}
Therefore, we get exactly \eqref{eq:4.11Toep}. 

The estimate \eqref{eq:5.13Toep24} follows from combining \eqref{eq:5.8jan24} 
and \eqref{eq:5.19jan24} for any test function in $\Omega^{n-1,n-1}_{\mathrm{c}}(U)$. 
Finally, since the quasi-psh functions $\{u_p\}_{p\gg 0}$ form a bounded set of 
$\mathcal{L}^1(K,\R)$ hence relatively compact by Proposition \ref{prop1appendix}, 
then there exists a subsequence $\{u_{p_j}\}_j$ which converges in 
$\mathcal{L}^1(K,\R)$ to a quasi-psh function $\widehat{f}$. 
Then applying \eqref{eq:5.22jan}, we get \eqref{eq:5.14limit}.
In particular, since $u_p$ converges uniformly to $0$ on $B$, 
we also have $\widehat{f}|_{B\cap U}\equiv 0$. This way, we finish our proof.
\end{proof}

When $f^2\geq c>0$ with some constant $c>0$, 
the lower bound in \eqref{eq:4.11Toep}
is trivial due to the asymptotic expansions of Bergman kernel function $P_p(x,x)$. 
When $f$ vanishes identically on a nonempty open subset $V$ of $U$ 
(which is given as in Theorem \ref{thm:5.5jan24}), 
by Theorem \ref{thm:2.4June}, we have
\begin{equation}
T^2_{f,p}(x,x)=\mathcal{O}_V(p^{-\infty}), 
\text{\,uniformly\,for\,all\,} x\in \overline{V}.
\end{equation}
In this case, the lower bound in \eqref{eq:4.11Toep} is nontrivial and provides
a better understanding of the decay of $T^2_{f,p}(x,x)$ on $U$.

Following the properties introduced in the proof of Theorem \ref{thm:5.5jan24}, 
we can generalize Proposition \ref{thm:5.5finite} to obtain the convergence of 
$\frac{1}{p}\mathbb{E}[[\Div(\bb{S}_{f,p})]]$ on the support of $f$. 
\begin{proof}[Proof of Theorem \ref{thm:6.2equi}]
We use the same notation as in the proof of Theorem \ref{thm:5.5jan24}. 
Note that $\overline{U}$ is compact in $X$. Since $f$ is of class $\mathscr{C}^1$ 
almost everywhere near $U$ and it has nonvanishing points,
then the family $\{u_p\}_p$ defined in \eqref{eq:fctup} is a bounded set in
$\mathcal{L}^2(K,\R)$ following Proposition \ref{prop1appendix}. 
Let $C_2>0$ be the constant such that for all $p\geq p_0$,
\begin{equation}
\|u_p\|^2_{\mathcal{L}^2(\overline{U},\R)}:=
\int_{\overline{U} }|u_p(x)|^2\mathrm{dV}(x)\leq C_2.
\label{eq:5.36feb24}
\end{equation}
As a consequence, for any $A\subset \overline{U}$ measurable, we have
\begin{equation}
\|u_p\|_{\mathcal{L}^1(A,\R)}\leq
\left(\int_{A }|u_p(x)|^2\mathrm{dV}(x)\right)^{1/2}
 \mathrm{Vol}(A)^{1/2}\leq \sqrt{C_2}\mathrm{Vol}(A)^{1/2}.
\end{equation}
Combining this with the fact that $\{u_p\}_p$ is a bounded in 
$\mathcal{L}^2(K,\R)$, 
we get that the family $\{u_p\}_p$ is uniformly integrable on $U$ 
(see \cite[Definition 4.5.1 and Proposition 4.5.3]{MR2267655}).

By Theorem \ref{thm:2.5Novfg}, 
the assumption on $f$ 
implies that $u_p$ converges to $0$ almost everywhere on $U$ with 
respect to the Lebesgue measure (or $\mathrm{dV}$).
Finally by the Lebesgue-Vitali theorem (see \cite[Theorem 4.5.4 and 
Corollary 4.5.5]{MR2267655}), 
we get that $u_p$ converges to $0$ in $\mathcal{L}^1$-norm on $U$, and 
\eqref{eq:6.0.7equi} follows.
Combining \eqref{eq:5.22jan} with \eqref{eq:6.0.7equi}, 
we conclude the convergence \eqref{eq:5.35feb24}. The proof is complete.
\end{proof}

\begin{Remark}
(i) If we use the bound of $\mathcal{L}^2$-norm for $u_p$
such as \eqref{eq:5.36feb24} in the proof of Theorem \ref{thm:5.5jan24}, 
then we replace the exponent $-{Cp}{\mathrm{Vol}(A_p)}^{-1}-C'\log{p}$ 
in the right-hand side of \eqref{eq:4.11Toep}
by $-{C_2p}/{\sqrt{\mathrm{Vol}(A_p)}}$.

\noindent
(ii) In Theorem \ref{thm:6.2equi}, if we assume a stronger condition on $f$, 
for example, if the function $|\log f^2|$ is integrable on $U$ with respect to 
$\mathrm{dV}$, 
then one can ask if the following convergence of $(1,1)$-currents on $U$ as $p\rightarrow +\infty$ holds,
\begin{equation}
\mathbb{E}[[\Div(\bb{S}_{f,p})]]|_U- p c_1(L,h_L)|_U\rightarrow 
\frac{\sqrt{-1}}{2\pi }\partial\overline{\partial} \log f^2 +c_1(E,h_E)|_U\, ?
\end{equation}
\end{Remark}

\subsection{Random zeros and lowest Toeplitz eigenvalues on
compact manifolds}\label{ss:compactandeigen}
Now we apply the results from previous Sections to the case
of compact complex manifolds. 
In this Section, we always assume $(X,\Theta)$ to be a 
connected compact Hermitian manifold and $(L,h_L)$ 
to be a positive holomorphic line bundle on $X$. 
A random section $\bb{S}_{f,p}$ as in \eqref{eq:1.5intropart}
can be equivalently defined by 
\eqref{eq:5.35jan24}.

As in the Definition \ref{def:1.6m} we can set by
the compactness of $X$,
\begin{equation}
\kappa(R^L,X):=\sup\{\max \mathrm{spec}\,(\dot{R}^L_x),\; x\in X\}\geq \varepsilon_0,
\label{eq:1.09cpt}
\end{equation}
and then set
\begin{equation}
m(X):=\left\lceil (6n+6)\frac{\kappa(R^L, X)}{\varepsilon_0}\right\rceil\in\N.
\label{eq:1.9cpt}
\end{equation}
At first, we apply Theorems \ref{thm:6.2} and 
\ref{thm:6.2equi} to the case $\supp{f}=X$.
\begin{Theorem}\label{thm:6.2compact}
Let $(X,J,\Theta)$ be a connected, compact complex Hermitian manifold and 
let $(L,h_{L})$, $(E,h_E)$ be holomorphic line bundles on $X$ with smooth 
Hermitian metrics. Assume $h_L$ to be positive. Fix a real bounded function 
$f$ which is of class $\mathscr{C}^{m(X)+1}$ almost everywhere on $X$ and 
with $\esupp\,{f}=X$ (or equivalently $\supp{f}=X$). Then we have $\P$-a.s.\ that
\begin{equation}
\lim_{p\rightarrow \infty}	
\left\|\frac{1}{p}[\Div(\bb{S}_{f,p})]- c_1(L,h_L)\right\|_{X,-2}
= 0,
\label{eq:6.0.7cptU}
\end{equation}
If $f$ is a real bounded function which is $\mathscr{C}^{1}$ almost everywhere 
on $X$ and with $\esupp\,{f}=X$, then as $p\rightarrow +\infty$, we have
\begin{equation}
\left\|\frac{1}{p} \log {T_{f,p}^2(x,x)}\right\|_{\mathcal{L}^1(X,\R)}\rightarrow 0,
\label{eq:6.0.7equicpt}
\end{equation}
and
\begin{equation}
\left\|\frac{1}{p}\mathbb{E}[[\Div(\bb{S}_{f,p})]]
- c_1(L,h_L)\right\|_{X,-2}\rightarrow 0.
\label{eq:5.35feb24cpt}
\end{equation}
\end{Theorem}

\begin{Remark}
The convergence of $\frac{1}{p} \log {T_{f,p}^2(x,x)}$ 
to the zero function in \eqref{eq:6.0.7equicpt} can fail if we take the 
pointwise limit of $\frac{1}{p} \log {T_{f,p}^2(x,x)}$. 
An example is given in Section \ref{ss5.3Jan}, where the smallest eigenvalue 
$\min \mathrm{Spec}(T_{f,p})$ is exponentially small in $p$.
\end{Remark}

A consequence of Theorem \ref{thm:6.2compact} is an improvement of the lower bound \eqref{eq:4.11Toep}.
\begin{Proposition}
Let $(X,J,\Theta)$ be a connected, compact complex Hermitian manifold and let $(L,h_{L})$, $(E,h_E)$ be holomorphic line bundles on $X$ with smooth Hermitian metrics. Assume $h_L$ to be positive. Fix a real bounded function $f$ (it can take negative values) which is of class $\mathscr{C}^1$ almost everywhere on $X$ and with $\esupp\,{f}=X$ (or equivalently $\supp{f}=X$). Then there exists a constant $c>0$ and a decreasing sequence of strictly positive numbers $\{r_p\}_{p\geq 1}$ with limit $\lim_{p\rightarrow \infty}r_p=0$ which depend only on $X$, $L$ and $f$ such that for all $p\geq 1$, any sequence of nonempty open subsets $\{A_p\}_{p\geq 1}$ of $X$, we have
	\begin{equation}
		\frac{1}{\mathrm{Vol}(A_p)}\int_{A_p} T^2_{f,p}(x,x)\mathrm{dV}(x)\geq \exp\left(-r_p p/\mathrm{Vol}(A_p)\right).
		\label{eq:4.14nov}
	\end{equation}
\end{Proposition}
\begin{proof}
By \eqref{eq:6.0.7equicpt} (which only require $f$ to be $\mathscr{C}^1$ almost everywhere on $X$), we fix a decreasing sequence $r_p>0$, $p\in\N$, such that
\begin{equation}
\left\|\frac{1}{p} \log {T_{f,p}^2(x,x)}\right\|_{\mathcal{L}^1(X,\R)}\leq r_p \rightarrow 0.
\end{equation}
Similarly to \eqref{eq:5.29Jan}, we get for any $A\subset X$,
\begin{equation}
-r_p p\leq \int_A \log {T^2_{f,p}(x,x)}\mathrm{dV}(x)\leq r_p p,
\label{eq:5.51feb24}
\end{equation}
When $\mathrm{Vol}(A)\neq 0$, by Jensen's inequality with $\log(t)$, we get
\begin{equation}
 \log\left(\frac{1}{\mathrm{Vol}(A)}\int_A T^2_{f,p}(x,x)\mathrm{dV}(x) \right) \geq \frac{1}{\mathrm{Vol}(A)}\int_A \log {T^2_{f,p}(x,x)}\mathrm{dV}(x) \geq -\frac{r_p p}{\mathrm{Vol}(A)}
 \label{eq:5.52feb24}
\end{equation}
This way, we conclude \eqref{eq:4.14nov}.
\end{proof}

Theorems \ref{thm:6.2} and \ref{thm:6.2equi} seem to indicate
is that the asymptotic behavior of the random zeros of 
$\frac{1}{p}[\mathrm{Div}(T_{f,p}\bb{S}_p)]$ 
depends mainly on the support of the function $f$
rather than the precise values of $f$. 
Following this thread, we formulate 
some interesting questions. 
For the sake of simplicity, we will restrict ourselves to 
to consider only non-negative functions $f$, unless explicitly stated otherwise.

At first, we consider the following equivalence relation on the
non-negative functions on $X$: for 
$f_1, f_2\in\mathcal{L}^\infty(X,\R_{\geq 0})$ 
(resp.\ $\mathscr{C}^k(X,\R_{\geq 0})$), 
we say they are comparable if there exists a constant 
$C=C(f_1, f_2)\geq 1$ such that on whole $X$,
\begin{equation}
\frac{1}{C} f_1\leq f_2 \leq C f_1.
\label{eq:5.44feb24}
\end{equation}
Let $\mathcal{EL}^\infty(X,\R_{\geq 0})$ 
(resp.\ $\mathcal{EC}^k(X,\R_{\geq 0})$) be the set of all equivalent classes of comparable functions in $\mathcal{L}^\infty(X,\R_{\geq 0})$ 
(resp.\ $\mathscr{C}^k(X,\R_{\geq 0})$), then it has a semigroup
structure given by the pointwise addition of functions 
with the zero element given by the zero function.
\begin{Lemma}
If $f_1, f_2\in\mathscr{C}^0(X,\R_{\geq 0})$ are comparable,
then they have exactly the same set-theoretic support,
that is $\{x\in X\;:\; f_1(x)>0\}=\{x\in X\;:\; f_2(x)>0\}$, 
therefore $\mathrm{supp}\, f_1=\mathrm{supp}\, f_2$. 
\end{Lemma}

On can then enquire about the relation between the asymptotic  
behaviors of the random zero sets $[\mathrm{Div}(T_{f_1,p}\bb{S}_p)]$
and $[\mathrm{Div}(T_{f_2,p}\bb{S}_p)]$ in the case
when $f_1$, $f_2$ are comparable functions on $X$, 
which might not be fully supported on $X$.
\begin{Question}\label{qn:5.15}
Let $(X,J,\Theta)$ be a connected, compact complex Hermitian manifold and let $(L,h_{L})$, $(E,h_E)$ be holomorphic line bundles on $X$ with smooth Hermitian metrics. Assume $h_L$ to be positive. Fix the sequence of the standard Gaussian random holomorphic sections $\{\bb{S}_{p}\}_p$. If two nontrivial functions $f_1, f_2\in\mathcal{L}^\infty(X,\R_{\geq 0})$ are comparable, is it true that $\P$-a.s.,
\begin{equation}
\lim_{p\rightarrow +\infty} \left\|\frac{1}{p}[\mathrm{Div}(T_{f_1,p}\bb{S}_p)]-\frac{1}{p}[\mathrm{Div}(T_{f_2,p}\bb{S}_p)]\right\|_{X,-2}=0,
  \end{equation}
 and that
\begin{equation}
\lim_{p\rightarrow +\infty} \left\|\E\left[\frac{1}{p}[\mathrm{Div}(T_{f_1,p}\bb{S}_p)]-\frac{1}{p}[\mathrm{Div}(T_{f_2,p}\bb{S}_p)]\right]\right\|_{X,-2}=0 
\end{equation}
 hold true?
\end{Question}

By Theorem \ref{thm:6.2compact}, 
when $f_1,f_2 \in \mathscr{C}^{m(X)+1}$
 have full support on $X$, then the above two equations hold true. 
So the difficult case is when $f_j$s are not fully supported. 
Note that if the results hold in general, they allow us to classify
the asymptotic behavior of the random zeros associated the function
$f$ via their equivalence class in $\mathcal{EL}^\infty(X,\R_{\geq 0})$.

A potential approach to a positive answer to Question \ref{qn:5.15}
is via establishing the following pointwise estimates for a pair 
of comparable functions $(f_1, f_2)$ for all $x\in X$,
\begin{equation}
\frac{1}{C^2}T^2_{f_1,p}(x,x)\leq T^2_{f_2,p}(x,x) \leq {C^2}T^2_{f_1,p}(x,x).
\label{eq:5.46feb24}
\end{equation}
Note that the condition \eqref{eq:5.44feb24} only implies 
\begin{equation}
\frac{1}{C}T_{f_1,p}(x,x)\leq T_{f_2,p}(x,x) \leq {C}T_{f_1,p}(x,x),
\label{eq:5.47feb24}
\end{equation}
which, in general, could not infer \eqref{eq:5.46feb24}. Moreover, near the nonvanishing smooth point $x$ of $f_j$, inequality \eqref{eq:5.46feb24} always holds true for sufficiently large $p$ after taking a bit large constant $C$ (by the asymptotic expansion \eqref{eq:2.43June}).

Another interesting question closely related to \eqref{eq:5.46feb24}
or \eqref{eq:5.47feb24} is to understand the asymptotic spectra of 
$T_{f,p}$ or $T^2_{f,p}$ as $p\rightarrow +\infty$;
this is already partially addressed in Section \ref{ssintro:minimal} 
and the last part Section \ref{ss:2.2feb}. 
Then in the sequel, we will try to understand more how the lowest
eigenvalues of $T_{f,p}$ decay to $0$ as $p\rightarrow+\infty$
when $f$ has vanishing points.

We start with a lemma for comparison of eigenvalues of Hermitian matrices.
\begin{Lemma}[Weyl's inequality, see {\cite[Theorem 4.3.7]{MR832183}}]
For $m\in\N_{\geq 1}$. Let $A$, $B$ be two Hermitian  $m\times m$-matrices, and let $\lambda_1 \leq \lambda_2 \leq \ldots \leq \lambda_m$ be the increasingly ordered eigenvalues of a Hermitian matrix, then we have the inequalities for the eigenvalues, for $j,\ell \in \{1,\ldots,m\}$,
\begin{equation}
\lambda_{j+\ell-m}(A+B)\leq \lambda_j(A)+\lambda_\ell(B)\leq \lambda_{j+\ell-1}(A+B),
\label{eq:5.46march}
\end{equation}
whenever the subscript makes sense.
\end{Lemma}
As a consequence, we have the following comparison results for the eigenvalues of Toeplitz operators.
\begin{Lemma}\label{lm:5.17feb24}
Let $(X,J,\Theta)$ be a connected, compact complex Hermitian manifold and let $(L,h_{L})$, $(E,h_E)$ be holomorphic line bundles on $X$ with smooth Hermitian metrics. Assume $h_L$ to be positive. Let $f_1, f_2\in \mathcal{L}^\infty(X,\R_{\geq 0})$ be such that each of them is not identically zero and $f_1\leq f_2$. For $j=1,2$, let $0<\lambda^p_\mathrm{min}(f_j)=\lambda^p_1(f_j)\leq \lambda^p_2(f_j)\leq\ldots\leq \lambda^p_{d_p}(f_j)= \lambda^p_\mathrm{max}(f_j)$ denote the ordered eigenvalues of $T_{f_j,p}$ respectively. Then we have for all $p\in\N$, $j=1,\ldots,d_p$,
\begin{equation}
\lambda^p_j(f_1)\leq \lambda^p_j(f_2).
\label{eq:5.47march}
\end{equation}

As a consequence, if $f_1, f_2\in \mathcal{L}^\infty(X,\R_{\geq 0})$ are comparable as in \eqref{eq:5.44feb24}, then we have for all $p\in\N$, $j=1,\ldots,d_p$,
\begin{equation}
\frac{1}{C}\lambda^p_j(f_1)\leq \lambda^p_j(f_2) \leq {C}\lambda^p_j(f_1).
\label{eq:5.50march24}
\end{equation}
\end{Lemma}
\begin{proof}
The comparison result \eqref{eq:5.50march24} is a direct consequence of the first part of this lemma, which is actually an application of \cite[Corollary 4.3.3]{MR832183}, we give a proof in the sequel.

Set $h:=f_2-f_1\geq 0$, then $T_{h,p}= T_{f_2, p} - T_{f_1,p}$ is a positive semidefinite Hermitian operator. Take $A=T_{f_1,p}$, $B=T_{h,p}$ and $\ell=1$ in \eqref{eq:5.46march}, we conclude 
\begin{equation}
\lambda^p_j(f_1)+\lambda^p_1(h)\leq \lambda^p_j(f_2).
\end{equation}
Since $\lambda^p_1(h)\geq 0$, we get \eqref{eq:5.47march}. The lemma is proved.
\end{proof}

Among all the eigenvalues of $T^2_{f,p}$, we now concentrate on the lowest nonzero one. These lowest eigenvalues are also interesting in the sense of their correspondence to the lowest energy for the ground state in the Berezin--Toeplitz quantization
as in the work of Deleporte \cite{Del2019,Del2020}. He obtained the asymptotic expansions for $\lambda^p_{\mathrm{min}}$ as well as for the corresponding eigensections when $f\geq 0$ is a smooth function with only nondegenerate vanishing points of order $2$. Then the general results about the lowest eigenvalues are expected, so that we put it as the following conjecture/question:

\begin{Question}\label{question13}
Let $(X,J,\Theta)$ be a connected, compact complex Hermitian manifold and let $(L,h_{L})$, $(E,h_E)$ be holomorphic line bundles on $X$ with smooth Hermitian metrics. Assume $h_L$ to be positive. For $f\in\mathscr{C}^\infty(X,\R)$ (it can take negative values), set $\kappa(f):=\sup_{x\in X}\mathrm{ord}_x(f)\in \N\cup\{+\infty\}$ (where $\mathrm{ord}_x(f)$ denotes the vanishing order of $f$ at $x$), show that as $p\rightarrow +\infty$
\begin{equation}
\min \mathrm{Spec}^*(T^2_{f,p})=\begin{cases} \simeq p^{-\kappa(f)} &,\text{\, if \,} \kappa(f)<+\infty; \\ \simeq e^{-c_f \sqrt{p}}&, \text{\, if \,} \kappa(f)=+\infty \text{\,and\,}  \supp{f}=X \\
\simeq e^{-c_f p} &, \text{\, if \,} \kappa(f)=+\infty \text{\,and\,} X\setminus \supp{f} \neq\varnothing \end{cases}
\label{eq:5.18Jan}
\end{equation} 
where $c_f>0$ is some constant depending on $f$, and $\mathrm{Spec}^*(T^2_{f,p}):=\mathrm{Spec}(T^2_{f,p})\setminus\{0\}$, the sign $\simeq$ means up to a multiplication of some nonzero constant. Moreover, a related interesting question is to describe the corresponding eigensections on $X$. Note that in Section \ref{ss5.3Jan}, we give the examples on $\mathbb{CP}^1$ where the lowest eigenvalues of $T^2_{f,p}$ are the cases listed in \eqref{eq:5.18Jan}. Following the results on the off-diagonal decay of $P_p(x,y)$ given in \cite{MR1979933, MR3903324, MR3903323, MR4167085,MR4236547}, it is possible that we need to assume the analyticity on $\Theta$ and $h_L$, $h_E$ to get the three nice cases in \eqref{eq:5.18Jan}.
\end{Question}

The results in Proposition \ref{lm:upperbounds} and Theorem \ref{thm:4.7Toeplitz} shows a partial answer to the above question for a non-negative function $f$. In the sequel, we give their proofs as consequence of Theorems \ref{thm:offdiagonal}, \ref{thm:2.2June}, \ref{thm:2.5Nov} and the pluripotential theory on $X$.

\begin{proof}[Proof of Proposition \ref{lm:upperbounds}]
In both case, we set $\lambda^p_{\mathrm{min}}(f)=\min \mathrm{Spec}(T_{f,p})$, since $f\geq 0$ is nontrivial, we always have $\lambda^p_{\mathrm{min}}(f)>0$. In the same time, we have the global point-wise estimate on $X$ for all $p$,
\begin{equation}
T_{f,p}(x,x)\geq \lambda^p_{\mathrm{min}}(f) P_p(x,x),
\end{equation}
that is
\begin{equation}
\lambda^p_{\mathrm{min}}(f)\leq \frac{T_{f,p}(x,x)}{P_p(x,x)}.
\label{eq:5.60march24}
\end{equation}

Assume the case (i). Now we apply Theorem \ref{thm:2.5Nov} for the point $x_0$ together with the information on the $Q_{f,x_0}(f)$ given in the last of Theorem \ref{thm:2.2June}, we get for all $p\gg 0$
\begin{equation}
\left|\frac{1}{p^n}T_{f,p}(x_0,x_0)-(Q_{2N,x_0}(f)\mathcal{P}_{x_0})(0,0)p^{-N}\right|\leq Cp^{-N-1/2}.
\label{eq:5.61march24}
\end{equation}
As a consequence, we have $T_{f,p}(x_0,x_0)\leq Cp^{n-N}$. 
In the same time, we have $P_p(x_0,x_0)=\bb{b}_0(x_0)p^n+\mathcal{O}(p^{n-1})$. Taking these estimates in \eqref{eq:5.60march24}, we obtain \eqref{eq:5.57march24}.

For case (ii), since $f$ is assumed to vanish with vanishing order $+\infty$, by Theorem \ref{thm:2.2June} and the same arguments for \eqref{eq:5.61march24}, we get that for any $\ell\in\N$, we have a constant $C_\ell$ such that $T_{f,p}(x_0,x_0)\leq C_\ell p^{n-\ell}$. This way, we get \eqref{eq:5.58march24}. The proof is completed.
\end{proof}

\begin{Remark}
We also have a uniform lower bound for $T_{f,p}(x,x)$ by an explicit calculation on $(Q_{2N,x_0}(f)\mathcal{P}_{x_0})(0,0)$ in \eqref{eq:5.61march24}: in fact, we have
\begin{equation}
(Q_{2N,x_0}(f)\mathcal{P}_{x_0})(0,0)=\bb{b}_0(x_0)\sum_{\alpha\in \N^{2n},|\alpha|=2N,} C_\alpha \frac{\partial^{2N}f}{\partial v^\alpha}(x_0)>0,
\end{equation}
where $C_\alpha$ is a constant depending only on $\alpha$ which is $>0$ if each component of $\alpha$ is even, and is $0$ otherwise.
As a consequence, if $f\in\mathscr{C}^{2N+1}(X,\R_{\geq 0})$ and $f$ vanishes on $X$ with the vanishing order at most $2N$, then there exists $C>0$ such that for all $p>0$, $x\in X$,
	\begin{equation}
	T_{f,p}(x,x)\geq Cp^{n-N}.
	\label{eq:5.62march24}
	\end{equation}
See also \cite[Proposition 5.6]{DrLM:2023aa}. 
Using instead the expansion \eqref{eq:2.43June} and Theorem \ref{thm:2.5Novfg}, one can prove the analogous results of \eqref{eq:5.62march24} for $T^2_{f,p}(x,x)$ with a general function $f\in\mathscr{C}^\ell(X,\R)$ which might take negative values, see also \cite[Proposition 5.16]{DrLM:2023aa}.
\end{Remark}

\begin{proof}[Proof of Theorem \ref{thm:4.7Toeplitz}]
Let $B$ denote a very small (nonempty) open ball in $X$, and we consider the case $f=\bb{1}_B$, the indicator function for $B$. In this case, $d_p=\mathcal{O}(p^n)$, and we know $T_{f,p}$ is injective and positive for all $p\gg 0$.

Now we set $\lambda^p_{\mathrm{min}}=\lambda^p_{\mathrm{min}}(f):= \min \mathrm{Spec}(T_{f,p}) >0$. For each $p\gg 0$, there exists a section $s_p^{\mathrm{min}}\in H^0(X,L^p\otimes E)$ with $\|s_p^{\mathrm{min}}\|_{\mathcal{L}^2(X,L^p\otimes E)}=1$ and
\begin{equation}
\lambda^p_{\mathrm{min}}=\langle T_{\bb{1}_B,p} s_p^{\mathrm{min}}, s_p^{\mathrm{min}}\rangle_{\mathcal{L}^2(X,L^p\otimes E)}=\int_B |s_p^{\mathrm{min}}(x)|^2_{h_p}\mathrm{dV}(x).
\label{eq:5.49feb24}
\end{equation}

Now we consider the family of integrable real functions $v_p(x):=\frac{1}{p}\log |s_p^{\mathrm{min}}(x)|^2_{h_p}$ on $X$. As a consequence of Poincar\'{e}-Lelong formula, $\{v_p\}_{p\gg 0}\subset\mathrm{PSH}(X, 4c_1(L,h_L))$, and this family of quasi-psh functions is bounded uniformly from above because of the uniform asymptotic expansion of $P_p(x,x)$ on $X$. 

Now we prove that $\{v_p\}_{p\gg 0}$ is a bounded subset of $\mathcal{L}^1(X,\R)$. If it is not true, we can assume that there exists a subsequence $v_{p_j}$ such that
\begin{equation}
\|v_{p_j}\|_{\mathcal{L}^1(X,\R)}\rightarrow +\infty, \text{\,as\,}j\rightarrow +\infty.
\end{equation}
Then we apply Proposition \ref{prop2appendix} to this subsequence and $X$, we conclude that $v_{p_j}$ has to converge uniformly to $-\infty$ on $X$, but this contracts to our condition $\|s_p^{\mathrm{min}}\|_{\mathcal{L}^2(X,L^p\otimes E)}=1$. Therefore, there exists a constant $C>0$ such that for all $p\gg 0$,
\begin{equation}
\|v_{p}\|_{\mathcal{L}^1(X,\R)}\leq C.
\end{equation}
Hence 
\begin{equation}
\int_B \log |s_p^{\mathrm{min}}(x)|^2_{h_p} \mathrm{dV}(x) \geq -Cp.
\end{equation}

Similar to \eqref{eq:5.31feb24}, we get
\begin{equation}
 \frac{1}{\mathrm{Vol}(B)}\int_B \log {|s_p^{\mathrm{min}}(x)|^2_{h_p}}\mathrm{dV}(x) \leq \log\left(\frac{1}{\mathrm{Vol}(B)}\int_B |s_p^{\mathrm{min}}(x)|^2_{h_p}\mathrm{dV}(x) \right)
 \label{eq:5.53feb24}
\end{equation}
Then combining the above inequality with \eqref{eq:5.49feb24}, we conclude
\begin{equation}
\lambda^p_{\mathrm{min}}\geq \mathrm{Vol}(B) e^{-Cp/\mathrm{Vol}(B)}.
\label{eq:5.54feb24}
\end{equation}
This shows that \eqref{eq:5.48feb24} holds with $f=\bb{1}_B$. In general, by our assumption on $f$ from the statement, we can always find a small nonempty open ball $B$ and a constant $c>0$ such that $f\geq c\bb{1}_B$, then \eqref{eq:5.48feb24} follows from Lemma \ref{lm:5.17feb24} for minimal eigenvalues and \eqref{eq:5.54feb24}. 

Now we prove \eqref{eq:5.48m2}, we use the technique of choosing coherent sections and the estimate \eqref{eq:2.36feb24}. Assume now $X\setminus \esupp f \neq \varnothing$, then we fix $x_0\not\in\esupp f $ and a sufficiently small $\delta>0$ such that the small geodesic ball $B(x_0,2\delta)\subset X\setminus \esupp f$. We also assume that the line bundle $L$ and $E$ can be trivialized over $B(x_0,2\delta)$. Let $e_{L,x_0}$, $e_{E,x_0}$ denote the respective unit (smooth) frames of $L_{x_0}$, $E_{x_0}$, a coherent section $S^p_{x_0}\in H^0(X,L^p\otimes E)$ associated to $e_{L,x_0}$, $e_{E,x_0}$ is the unique section $S^p_{x_0}$ such that for any $s_p\in H^0(X,L^p\otimes E)$, we have
\begin{equation}
s_p(x_0)=\langle s, S^p_{x_0}\rangle_{\cLL(X,L^p\otimes E)}e_{L,x_0}^{\otimes p}\otimes e_{E,x_0}.
\end{equation}
In fact, we have $S^p_{x_0}(x)=P_p(x,x_0)e_{L,x_0}^{\otimes p}\otimes e_{E,x_0}$. Then we have
\begin{equation}
\|S^p_{x_0}\|^2_{\cLL(X,L^p\otimes E)}=P_p(x_0,x_0)\sim \bb{b}_0(x_0)p^n,
\end{equation}
and 
\begin{equation}
|S^p_{x_0}(y)|_{h_{p,y}}=|P_p(y,x_0)|_{h_{p,y}\otimes h_{p,x_0}^*}.
\end{equation}

As a consequence, we get
\begin{equation}
\begin{split}
\langle T_{f,p} S^p_{x_0},S^p_{x_0}
\rangle_{\mathcal{L}^2} &=
\int_X f(y)|S^p_{x_0}(y)|^2_{h_{p,y}}\mathrm{dV}(y)\\
&\leq \|f\|_{\mathcal{L}^\infty}
\int_{\mathrm{dist}(y,x_0)\geq \delta}
|P_p(y,x_0)|^2_{h_{p,y}\otimes h_{p,x_0}^*}\mathrm{dV}(y)\\
&\leq C_{\delta,A}\mathrm{Vol}(X) 
\|f\|_{\mathcal{L}^\infty} e^{-A\sqrt{p\log{p}}},
\end{split}
\end{equation}
where the last line follows from the estimate \eqref{eq:2.36feb24}.

Finally, we get
\begin{equation}
\lambda^p_{\mathrm{min}}\leq 
\frac{\langle T_{f,p} S^p_{x_0},S^p_{x_0}
\rangle_{\mathcal{L}^2}}{\|S^p_{x_0}\|^2_{\cLL(X,L^p\otimes E)}}
\leq C\frac{e^{-A\sqrt{p\log{p}}}}{p^n},
\end{equation}
then \eqref{eq:5.48m2} follows. The proof is complete.
\end{proof}

\begin{Remark}\label{Rm:5.23m}
As in \cite{MR4167085, MR4236547}, when $\Theta$, $h_L$, $h_E$ are analytic, we have a sharper version of \eqref{eq:2.36feb24} as follows
\begin{equation}
\sup_{\mathrm{dist}(x,y)\geq \delta}|P_p(x,y)|_{h_{p,x}\otimes h^*_{p,y}}\leq C_{\delta} e^{-c_\delta p}.
\end{equation}
Then in this case, when $f$ is not fully supported on $X$, we have
\begin{equation}
 \min \mathrm{Spec}(T_{f,p})\leq C' e^{-c'p},
 \label{eq:5.78m2}
 \end{equation}
 which fits exactly the third case in \eqref{eq:5.18Jan} under the analyticity condition.
\end{Remark}

As a consequence, we get the uniform point-wise lower bound for $\log T^2_{f,p}(x,x)$ on a compact manifold $X$.
\begin{Corollary}
With the same geometric assumption in Theorem \ref{thm:4.7Toeplitz} or Lemma \ref{lm:5.17feb24}, for $f\in \mathcal{L}^{\infty}(X,\R_{\geq 0})$ which is not identically zero and is continuous near a nonvanishing point, there exists a constant $C_f>0$ such that for all $p\gg 0$, $x\in X$,
\begin{equation}
\log T^2_{f,p}(x,x)\geq -C_f p.
 \label{eq:5.77feb24}
 \end{equation}
\end{Corollary}

\begin{Definition}
For any nontrivial $f\in \mathscr{C}^0(X,\R_{\geq 0})$, we set
\begin{equation}
c(f):=\liminf_{p\rightarrow +\infty} \frac{1}{p}\log \lambda_{p,\mathrm{min}}(f)\in\R_{\leq 0}.
\end{equation}
By Theorem \ref{thm:4.7Toeplitz}, the limit $c(f)$ always exists (that is, finite) and is nonpositive. For $f\equiv 0$, we set $c(0):=-\infty$.
\end{Definition}

The following result is an easy consequence of Lemma \ref{lm:5.17feb24}.
\begin{Proposition}
When $f>0$ is a continuous function on $X$ and never vanishes, we have $c(f)=0$. If $f_1$, $f_2\in \mathscr{C}^0(X,\R_{\geq 0})$ are comparable in the sense of \eqref{eq:5.44feb24}, then we have
\begin{equation}
c(f_1)=c(f_2).
\end{equation}
This way, we define a function $c: \mathcal{EC}^0(X,\R_{\geq 0})\rightarrow [-\infty,0]$.  
\end{Proposition}

\subsection{Examples and simulations on Riemann sphere}
\label{ss5.3Jan}

In this Section, we present the examples of the Toeplitz operators on $\mathbb{CP}^1$, in particular, we give the explicit computations on their spectra. In the last part, we present the random zeros of the sections given by certain Toeplitz operators acting on the $\mathrm{SU}(2)$-random polynomials on $\C$.

Let us start with our basic settings. We consider one standard chart 
$U_0\simeq \C$ for $\mathbb{CP}^{1}$. In this chart, the Fubini-Study metric is 
given by
\begin{equation}
	\Theta=\omega_{\mathrm{FS}}=\frac{\sqrt{-1}}{2\pi}\frac{dz\wedge 
	d\bar{z}}{(1+|z|^{2})^{2}}.
\end{equation}
It is also the volume form on $\C$ that we will use.

For $p\in \mathbb{N}_{>0}$, let $\mathscr{O}(p)$ denote the 
holomorphic line bundle on $\mathbb{CP}^{1}$ of degree $k$. Then 
$\mathscr{O}(p)=\mathscr{O}(1)^{\otimes p}$, we equip 
$\mathscr{O}(1)$ with the standard Fubini-Study Hermitian metric 
$h_{\mathrm{FS}}$ and equip $\mathscr{O}(p)$ with the induced 
Hermitian metric. Then 
$\Theta=\omega_{\mathrm{FS}}=c_{1}(\mathscr{O}(1),h_{\mathrm{FS}})$.

On this chart, the global holomorphic sections of 
$\mathscr{O}(p)$ are given by all the polynomials in $z$ with degree 
$\leq p$, i.e.,
\begin{equation}
	H^{0}(\mathbb{CP}^{1},\mathscr{O}(p))=\mathrm{Span}_{\C}\{1,z,\cdots, z^{p}\}.
\end{equation}
The canonical orthonormal basis of 
$H^{0}(\mathbb{CP}^{1},\mathscr{O}(p))$ with respect to the 
$\cLL$-inner product is given by
\begin{equation}
	S^{p}_{j}(z)=\sqrt{(p+1)\binom{p}{j}} z^{j},\;\; j=0,1,\cdots, p.
\end{equation}
Note that $\mathrm{SU}(2)$ acts on $\mathbb{CP}^{1}$ holomorphically 
and transitively.  On the chart $\C$ the action is given by
\[g\cdot z = \frac{a+b z}{c+d z}\in\C, \quad
g=\begin{bmatrix}
	a & b\\
	c & d
\end{bmatrix}\in \mathrm{SU}(2).\]

The action of $\mathrm{SU}(2)$ always lifts to $\mathscr{O}(p)$
and preserves the metrics 
$\Theta$ and $h_{\mathrm{FS}}$. In particular, $\mathrm{SU}(2)$ acts 
on $H^{0}(\mathbb{CP}^{1},\mathscr{O}(p))$ isometrically with respect 
to $\cLL$-inner product. 

Note that for the explicit expressions of eigenvalues of $T_{f,p}$ in Example \ref{ex:5.27} for $f_k$ and in Section \ref{ssb:5.5.2} for $\bb{1}_r$ ($r=1,1/4$) were already computed in \cite{MR3822756,MR3613959}.

\subsection{Lowest eigenvalues of Toeplitz operators for fully supported functions}\label{ssb:5.5.1}
In this part, we compute the spectra of Toeplitz operators associated to some special functions, as examples for the first two cases listed in \eqref{eq:5.18Jan}.

\begin{Example}\label{ex:5.27}
For $k\in\N_{\geq 1}$, on $U_0\simeq \C$, we take the function 
\begin{equation}
f_k(z):=\frac{|z|^{2k}}{(1+|z|^2)^{k}}
\end{equation}
Then $f_k$ is a smooth nonnegative function on $\mathbb{CP}^1$ 
with the only vanishing point at $z=0$, and the vanishing order is $2k$.

At first, we have for $0\leq j\neq \ell\leq p$, then
\begin{equation}
\langle T_{f_k,p} S^{p}_{j}, S^{p}_{\ell}\rangle_{\cLL(\mathbb{CP}^1,\mathscr{O}(p))}=0
\label{eq:5.39Jan}
\end{equation}
Then we have
\begin{equation}
\mathrm{Spec}\; T_{f_k,p}=\{\langle T_{f_k,p} S^{p}_{j}, S^{p}_{j}\rangle_{\cLL(\mathbb{CP}^1,\mathscr{O}(p))} >0 \;:\; j=0,1,\ldots,p\}.
\end{equation}
By elementary techniques on the integrals, we get
\begin{equation}
\langle T_{f_k,p} S^{p}_{j}, S^{p}_{j}\rangle_{\cLL(\mathbb{CP}^1,\mathscr{O}(p))} = \frac{(j+k)!}{j!}\frac{(p+1)!}{(k+p+1)!}.
\end{equation}
It is clear that this quantity increases as $j$ grows. Then we have the following smallest ($j=0$) and biggest ($j=p$) eigenvalues of $T_{f_k,p}$:
\begin{equation}
\lambda^p_\mathrm{min}=\frac{k!(p+1)!}{(k+p+1)!},\; \lambda^p_\mathrm{max}=1-\frac{k}{p+k+1}.
\end{equation}
In particular, we have, as $p\rightarrow +\infty$,
\begin{equation}
\lambda^p_\mathrm{min}=k!p^{-k}(1+\frac{k(k+3)}{2p}+\mathcal{O}(p^{-2})).
\label{eq:5.43Jan}
\end{equation}
The corresponding eigensection is $S^p_0$.

The asymptotic expansion in \eqref{eq:5.43Jan} gives an example for the first case listed in \eqref{eq:5.18Jan}.
\end{Example}

\begin{Example}\label{ex:5.22}
On $U_0\simeq \C$, we take the function 
\begin{equation}
f(z):=e^{-\frac{1}{|z|^2}}.
\end{equation}
Then $f_k$ is a smooth nonnegative function on $\mathbb{CP}^1$ with the only vanishing point at $z=0$, and the vanishing order is $+\infty$, but we have $\supp {f}=\mathbb{CP}^1$.

Similar to \eqref{eq:5.39Jan}, for $0\leq j\neq \ell\leq p$, we have
\begin{equation}
\langle T_{f,p} S^{p}_{j}, S^{p}_{\ell}\rangle_{\cLL(\mathbb{CP}^1,\mathscr{O}(p))}=0
\end{equation}
Then we have
\begin{equation}
\mathrm{Spec}\; T_{f,p}=\{\langle T_{f,p} S^{p}_{j}, S^{p}_{j}\rangle_{\cLL(\mathbb{CP}^1,\mathscr{O}(p))} >0 \;:\; j=0,1,\ldots,p\}.
\end{equation}
By elementary techniques on the integrals, we get
\begin{equation}
\lambda^p_j:=\langle T_{f,p} S^{p}_{j}, S^{p}_{j}\rangle_{\cLL(\mathbb{CP}^1,\mathscr{O}(p))} = \frac{(p+1)!}{j!}\mathrm{HyperU}(p+2,j+2,1),
\end{equation}
where $\mathrm{HyperU}(a,b,z)$ denotes the Tricomi confluent hypergeometric function, that is defined as
\begin{equation}
\mathrm{HyperU}(a,b,z):=\frac{1}{\Gamma(a)}\int_0^{+\infty} e^{-zt}t^{a-1}(1+t)^{b-a-1}dt.
\end{equation}
Then we can rewrite
\begin{equation}
\lambda^p_j=\int_0^{+\infty} e^{-t}\left(\frac{t}{1+t}\right)^{p+1} \frac{(1+t)^j}{j!}dt.
\end{equation}
A direct computation shows that 
\begin{equation}
\lambda^p_{\mathrm{min}}:=\lambda^p_0< \lambda^p_1<\cdots<\lambda^p_p=:\lambda^p_{\mathrm{max}}.
\end{equation}
For the biggest ($j=p$) eigenvalues of $T_{f,p}$, we have the following explicit formula
\begin{equation}
\begin{split}
\lambda^p_\mathrm{max}&=1+\sum_{j=1}^p (-1)^j\frac{(p-j)!}{p!}+\frac{(-1)^{p+1}}{p!}\int_0^{+\infty}\frac{e^{-t}}{1+t}dt\\
&=1-\frac{1}{p}+\mathcal{O}(p^{-2}),\qquad\text{as } p\rightarrow +\infty.
\end{split}
\end{equation}
For the smallest eigenvalues, we have
\begin{equation}
\lambda^p_\mathrm{min}=\int_0^{+\infty}e^{-t}\left(\frac{t}{1+t}\right)^{p+1}dt=(p+1)\int_0^{+\infty}e^{-t}\frac{t^p}{(1+t)^{p+2}}dt.
\end{equation}
Now let us give the asymptotic behavior of $\lambda^p_\mathrm{min}$ as $p\rightarrow +\infty$. At first, we have
\begin{equation}
\begin{split}
(p+1)\int_0^{+\infty}e^{-t}\frac{t^p}{(1+t)^{p+2}}dt &\leq (p+1)\int_0^{\sqrt{p}/2}\frac{t^p}{(1+t)^{p+2}}dt \\
&\qquad + (p+1)\int_{\sqrt{p}/2}^{+\infty}e^{-t}\left(\frac{t}{1+t}\right)^{p}\frac{dt}{(1+t)^2}\\
&\leq \left(\frac{\sqrt{p}/2}{1+\sqrt{p}/2}\right)^{p+1}+\frac{p+1}{\sqrt{p}/2+1}e^{-2\sqrt{p}+\mathcal{O}(1)},
\end{split}
\end{equation}
where the second term follows from finding the maximum of $e^{-t}\left(\frac{t}{1+t}\right)^{p}$.  Using the fundamental inequality for $t>0$,
\begin{equation}
\left(\frac{t}{1+t}\right)^{1+t}\leq e^{-1}\leq \left(\frac{t}{1+t}\right)^{t},
\label{eq:5.54Jan}
\end{equation}
we get (for all sufficiently large $p$)
\begin{equation}
\begin{split}
(p+1)\int_0^{+\infty}e^{-t}\frac{t^p}{(1+t)^{p+2}}dt \leq 2\sqrt{p}e^{-2\sqrt{p}+\mathcal{O}(1)}.
\end{split}
\label{eq:5.55Jan}
\end{equation}
For the lower bound, we get
\begin{equation}
\begin{split}
(p+1)\int_0^{+\infty}e^{-t}\frac{t^p}{(1+t)^{p+2}}dt &\geq (p+1)e^{-\sqrt{p}}\int_0^{\sqrt{p}}\frac{t^p}{(1+t)^{p+2}}dt \\
&\geq e^{-\sqrt{p}}\left(\frac{\sqrt{p}}{1+\sqrt{p}}\right)^{p+1}\\
&\geq e^{-2\sqrt{p}-1/\sqrt{p}},
\end{split}
\label{eq:5.56Jan}
\end{equation}
where the last estimate follows from \eqref{eq:5.54Jan}.

Finally, we combining \eqref{eq:5.55Jan} and \eqref{eq:5.56Jan}, we get an asymptotic formula for $\lambda^p_\mathrm{min}$:
\begin{equation}
\lambda^p_\mathrm{min}=e^{-2\sqrt{p}(1+o(1))},
\label{eq:5.48Jan}
\end{equation}
or equivalently
\begin{equation}
\lim_{p\rightarrow \infty}\frac{\log \lambda^p_\mathrm{min}}{\sqrt{p}}=-2.
\end{equation}
This gives an example for the second case listed in \eqref{eq:5.18Jan}.

\end{Example}

\subsection{Spectrum of Toeplitz operators for indicator functions}
\label{ssb:5.5.2}
For $r>0$, set the function $\bb{1}_{r}$ on $\C$ as
\begin{equation}
	\bb{1}_{r}(z)=\bb{1}_{\mathbb{D}(0,r)}(z),
\end{equation}
where $\bb{1}_{\mathbb{D}(0,r)}$ denotes the indicator function 
for the open disc $\mathbb{D}(0,r)$ of radius $r$. In the sequel, for 
$U\subset \C$ or $\mathbb{CP}^{1}$, let $\mathrm{Vol}(U)$ denote the 
Fubini-Study volume of $U$. In particular, we have
\begin{equation}
	\mathrm{Vol}(\mathbb{D}(0,r))=\frac{r^{2}}{1+r^{2}}.
\end{equation}
This disc $\mathbb{D}(0,r)$ is a geodesic ball in $\mathbb{CP}^{1}$ 
with radius $\frac{1}{\sqrt{\pi}}\arctan{r}$ (with respect to 
Fubini-Study metric).

Now we consider the Toeplitz operators $T_{\bb{1}_{r},p}$ acting on 
$H^{0}(\mathbb{CP}^{1},\mathscr{O}(p))$. 
\begin{Proposition}
	The spectrum of $T_{\bb{1}_{r},p}$ is given by $\{\lambda^{p}_{j}(r), 
	j=0,1,\cdots, p\}$ where
	\begin{equation}
		\lambda^{p}_{j}(r)=(1+r^{2})^{-p-1}\sum_{i=0}^{p-j} 
		\binom{p+1}{p-i-j}r^{2i+2j+2}.
		\label{eq:5.91ss5.3}
	\end{equation}
	Moreover, we have
	\begin{equation}
		\begin{split}
			&\lambda^{p}_{\mathrm{max}}(r):=\max\{\lambda^{p}_{j}(r), 
	j=0,1,\cdots, p\}=\lambda^{p}_{0}(r)=1-\frac{1}{(1+r^{2})^{p+1}},\\
	&\lambda^{p}_{\mathrm{min}}(r):=\min\{\lambda^{p}_{j}(r), 
	j=0,1,\cdots, 
	p\}=\lambda^{p}_{p}(r)=\left(\frac{r^{2}}{1+r^{2}}\right)^{p+1}=\mathrm{Vol}(\mathbb{D}(0,r))^{p+1}.
		\end{split}
		\label{eq:5.7July}
	\end{equation}
	For $j=0,1,\cdots, p$, the $\lambda^{p}_{j}(r)$-eigenspace is 
	spanned by $S^{p}_{j}$ (hence is $1$-dimensional).
\end{Proposition}
\begin{proof}
By definition, we have
\begin{equation}
\begin{split}
\lambda^{p}_{j}(r)&=\langle T_{\bb{1}_{r},p} S^p_j, S^p_j\rangle_{\mathcal{L}^2}\\
&=(p+1)\binom{p}{j}\int_0^{r^2} \frac{t^j}{(1+t)^{p+2}}dt\\
&=(p+1)\binom{p}{j}\int_{0}^1 \frac{t^j}{(1+r^2 t)^{p+2}}dt\\
&=r^{2j+2}\frac{p+1}{j+1}\binom{p}{j} {}_2F_1(p+2,j+1;j+2;-r^2),
\end{split}
\end{equation}
where ${}_2F_1(a,b;c;z)$ denotes the hypergeometric function ${}_2F_1$. Using the fact that $p$ and $j$ are non-negative integers, the value of ${}_2F_1(p+2,j+1;j+2;-r^2)$ can be worked out explicitly so that we get exactly \eqref{eq:5.91ss5.3}. The rest part is clear.
\end{proof}

Using the $\mathrm{SU}(2)$-symmetry, we get the following corollary.
\begin{Corollary}
	Let $\mathbb{B}$ be a nontrivial geodesic ball in $\mathbb{CP}^{1}$ with 
	$\mathrm{Vol}(\mathbb{B})\in\;]0,1]$, then we have
	\begin{equation}
		\min\mathrm{Spec}(T_{\bb{1}_{\mathbb{B}},p})=\mathrm{Vol}(\mathbb{B})^{p+1}.
		\label{eq:5.8July}
	\end{equation}
	In particular, if $\mathrm{Vol}(\mathbb{B})<1$, then 
	$\min\mathrm{Spec}(T_{\bb{1}_{\mathbb{B}},p})$ decays 
	exponentially to $0$ as $p\rightarrow +\infty$.
\end{Corollary}

Finally, as a consequence of Lemma \ref{lm:5.17feb24}, we get
\begin{Proposition}
	For a $f\in \mathcal{L}^{\infty}(\mathbb{CP}^{1},\R_{\geq 0})$, if there 
	exists two geodesic ball $\mathbb{B}\subset \mathbb{B}'$ and two 
	constants $C\geq c>0$ such that
	\begin{equation}
		C\bb{1}_{\mathbb{B}'}\geq f \geq c\bb{1}_{\mathbb{B}},
	\end{equation}
	Then we have for $p\geq 1$,
	\begin{equation}
		C\mathrm{Vol}(\mathbb{B}')^{p+1}\geq 
		\min\mathrm{Spec}(T_{f,p}) \geq 
		c\mathrm{Vol}(\mathbb{B})^{p+1}.
		\label{eq:5.10July}
	\end{equation}
	In particular, when $\mathrm{Vol}(\mathbb{B}')<1$, 
	$\min\mathrm{Spec}(T_{f,p})$ decays 
	exponentially to $0$ as $p\rightarrow +\infty$.
\end{Proposition}

\subsection{Simulations of random zeros on the support
of the symbol}\label{ss:5.4.3new}
In this part, we present some simulation results for the zeros of $\bb{S}_{f,p}$ on $\mathbb{CP}^1$, where the function $f$ is given as in Sections \ref{ssb:5.5.1} and \ref{ssb:5.5.2}. 

We now explain our computation model. Fix a concrete choice of $f$ as above, we already know the precise spectrum of $T_{f,p}$: for $j=0,\ldots, p$, let $\lambda^p_j(f)$ be the eigenvalues of $T_{f,p}$ with the eigensection $S^p_j$. Then we write on $U_0\simeq \C$,
\begin{equation}
\bb{S}_{f,p}(z):=\sum_{j=1}^p \eta_j \lambda^p_j(f) \sqrt{(p+1)\binom{p}{j}} z^{j},
\label{eq:5.91feb24}
\end{equation}
where $\{\eta_j\}$ is a sequence of i.i.d.\ standard complex Gaussian random variables, such random variables can be simulated properly by the mathematical computating softwares. Our random section $\bb{S}_{f,p}$ now becomes a random polynomial on $\C$ of degree $p$. In the special case where $f\equiv 1$, $\bb{S}_{1,p}=\bb{S}_p$ is exactly the $\mathrm{SU}(2)$-polynomial.

In the sequel simulations, we always compare our random zeros with the expected distribution $\omega_{\mathrm{FS}}$ on $U_0\simeq \C$. In order to visualize the comparison, we will classify all the random zeros for each simulation according to their Fubini-Study norms (defined below as $r_{\mathrm{FS}}$), then we draw the corresponding density histograms. Such density histograms are viewed as the approximations to the probability density function with respect to the Fubini-Study norms of zeros of the polynomial $\bb{S}_{f,p}(z)$; in the same time, the corresponding probability density function of $\omega_{\mathrm{FS}}$ in the Euclidea norm $r=|z|$ is given by the function
\begin{equation}
\R_{\geq 0}\ni r\mapsto \frac{2r}{(1+r^2)^2}.
\label{eq:5.91k24}
\end{equation}
We consider the polar coordinate on $U_0\simeq \C$ with respect to $\omega_{\mathrm{FS}}$, so that for any $z\in\C$, $r=|z|$, it has the $\omega_{\mathrm{FS}}$-Riemannian distance $r_\mathrm{FS}=\frac{1}{\sqrt{\pi}}\arctan {r}$ from $z=0$, in particular, $r_{\mathrm{FS}}=\frac{\sqrt{\pi}}{2}$ corresponds to the North Pole $\{\infty\}\in\mathbb{CP}^1$. So that the diameter of $(\mathbb{CP}^1, \omega_{\mathrm{FS}})$ is $\frac{\sqrt{\pi}}{2}\simeq 0.88622\ldots\,$, and the equator is given by $r=|z|=1$ or the circle with $r_\mathrm{FS}=\frac{\sqrt{\pi}}{4}\simeq 0.44311\ldots\,$. 
For $z\in \C\simeq U_0$, $r_{\mathrm{FS}}=\frac{1}{\sqrt{\pi}}\arctan {|z|}$ is called the Fubini-Study norm of $z$. By \eqref{eq:5.91k24}, in terms of the Fubini-Study norm $r_{\mathrm{FS}}\in [0,\frac{\sqrt{\pi}}{2}]$, the density function of $\omega_{\mathrm{FS}}$ on $\mathbb{CP}^1$ is given by the function
\begin{equation}
\psi(r_{\mathrm{FS}})=\sqrt{\pi } \sin \left(2 \sqrt{\pi } r_{\mathrm{FS}}\right).
\label{eq:FScurve}
\end{equation}
Then by plotting the graph of $\psi$ along with the density histograms of the Fubini-Study norms of zeros from the simulations, we can visualize directly the differences between $\omega_{\mathrm{FS}}$ and the simulated random zeros.

Through all the simulations, we also keep the total number of random zeros to be a fixed large number such as $20000$ for each comparison. This means, if we consider the degree $p$ (such as $p=20,50,100,200$, which is a factor of $20000$), we will repeat $\mathrm{RN}:=20000/p$ times of realizing $\bb{S}_{f,p}(z)$ and computing its zeros ($\mathrm{RN}$ is the short for {\it repeating number}), this way, we will get in total $20000$ points as random zeros, then we draw the density histogram with respect to their Fubini-Study norms described as above.

We use MATLAB (version 23.2, R2023b) to perform the aforementioned simulations. Note that for large $p$ (such as $100$, $200$ or bigger), the combinatorial numbers in the coefficients of \eqref{eq:5.91feb24} can be extremly big, and the software can only compute properly their numerical approximations due to the precision limit on a laptop. Such simulations are not suitable for a precise quantitative analysis but they are good enough for our purpose of visualizing the random zeros and comparing them with the expected one. The following Figure \ref{fig:su2example} shows one simulation for the zeros of $\mathrm{SU}(2)$-polynomial. The picture in the left-hand side plots the roots that lie inside the square of side length $12$ among all $20000$ roots obtained via $\mathrm{RN}=1000$ realizations of $\mathrm{SU}(2)$-polynomials of degree $p=20$, and the one in the right-hand side compares the density histogram of the Fubini-Study norms of all these $20000$ roots with the expected density function $\psi(r_{\mathrm{FS}})$ (plotted as the red curve). From this comparison, we can see how well they fit to each other when $p=20$.
\begin{figure}[h]
\centering
\includegraphics[width=0.8\textwidth]{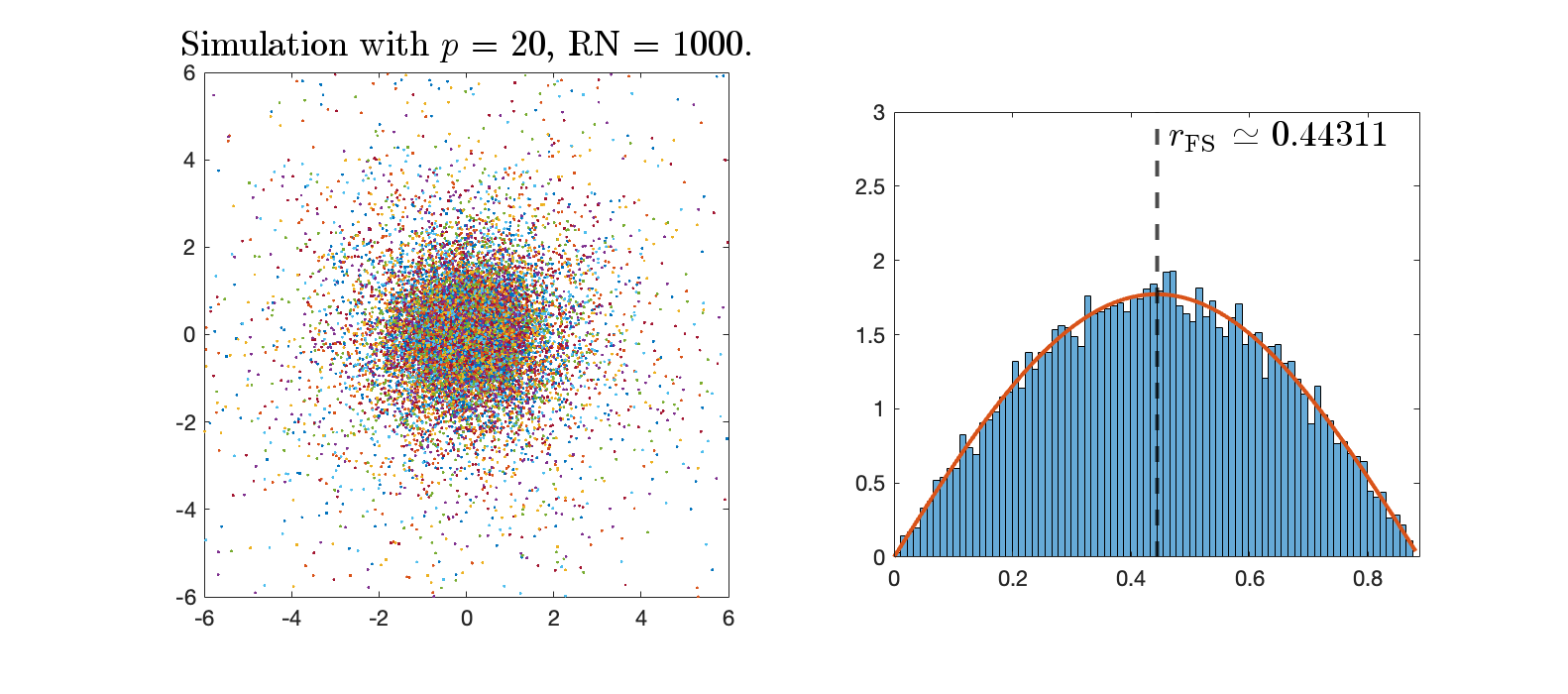}
\caption[Random zeros of $\mathrm{SU}(2)$-polynomial]{Comparison of zeros of $\mathrm{SU}(2)$-polynomial with $\omega_{\mathrm{FS}}$ on $\mathbb{CP}^1$. The density function $\psi(r_{\mathrm{FS}})$ is plotted as the red curve in the right-hand side.}
\label{fig:su2example}
\end{figure}

Now we consider two examples: $f_3$ from Section \ref{ssb:5.5.1} and $\bb{1}_1$ from Section \ref{ssb:5.5.2}.  
\begin{Example}[The function $f_3(z)=\frac{|z|^6}{(1+|z|^2)^3}$]
Following Example \ref{ex:5.27}, let us consider the function $f_3(z)=\frac{|z|^6}{(1+|z|^2)^3}$. The only vanishing point for $f_3$ is $z=0$ (corresponding to $r=0$), and the vanishing order is $6$. Since $\supp\, f_3=\mathbb{CP}^1$, by Theorem \ref{thm:6.2compact}, the random zeros $\frac{1}{p}[\Div(\bb{S}_{f_3,p})]$ will converge to $\omega_{\mathrm{FS}}$ as $p\rightarrow +\infty$. We shall expect a nice approximation as shown in the right-hand side of Figure \ref{fig:su2example} for sufficently large $p$.

\begin{figure}[H]
\centering
\includegraphics[width=0.8\textwidth]{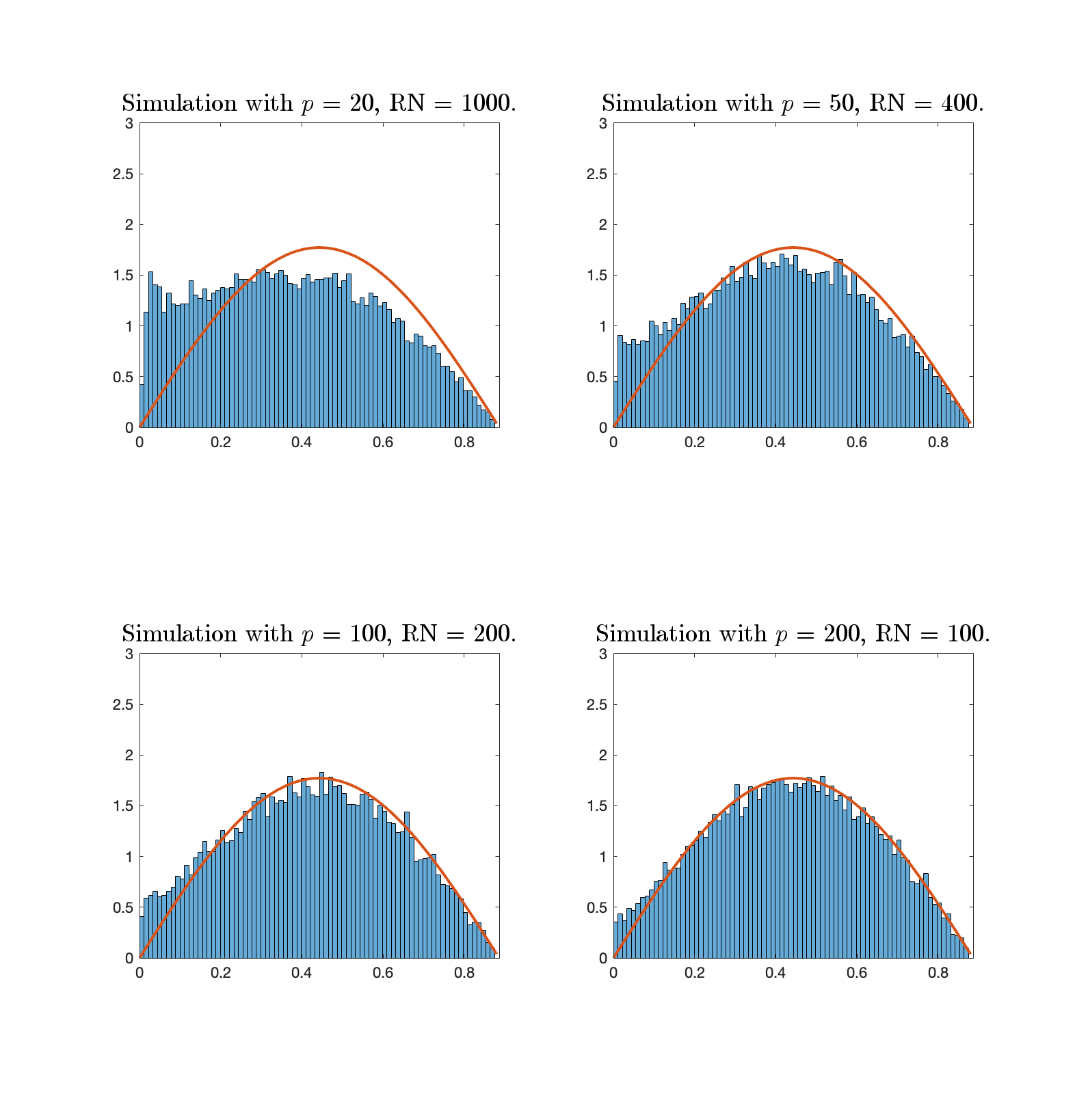}
\caption[Simulations for random zeros of $\bb{S}_{f_3, p}$]{Comparison of zeros of $\bb{S}_{f_3 , p}$ ($p=20,50,100,200$)  with $\omega_{\mathrm{FS}}$ on $\mathbb{CP}^1$. The density function $\psi(r_{\mathrm{FS}})$ is plotted as the red curve in each picture.}
\label{fig:f3example}
\end{figure}

We did $4$ simulations for different degrees: $p=20,\, 50,\, 100,\, 200$. The results are displayed as $4$ pictures in Figure \ref{fig:f3example}. It is a straighforward observation that the simulated random zeros approximate better and better to the expected one represented by $\psi(r_{\mathrm{FS}})$ as the degree $p$ grows. This is exactly the main point proved in Theorem \ref{thm:6.2compact}. In the same time, we see also that, different from the case of $\mathrm{SU}(2)$-polynomial in Figure \ref{fig:su2example}, the result with degree $p=20$ in Figure \ref{fig:f3example} does not fit nicely with $\psi(r_{\mathrm{FS}})$. The main deviation happens around $z=0$ ($r_\mathrm{FS}=0$), which is exactly the unique vanishing point of the function $f_3(z)$. Roughly speaking, this vanishing point is a sort of {\it global minimizer} of $\log T_{f,p}^2(x,x)$ on $\mathbb{CP}^1$, such that for small $p$, one should observe more zeros around this {\it minimizer}.

Analogous to this example, one can also consider the function $f(z)=e^{-1/|z|^2}$ discussed in Example \ref{ex:5.22}, then the point $z=0$ will still behave abnormally as in Figure \ref{fig:f3example}, and since the vanishing order is $+\infty$, we will need a very large $p\gg 200$ to observe a nice approximation of the simulation result to $\psi(r_\mathrm{FS})$.

\end{Example}

\begin{Example}[The function $\bb{1}_1(z)=\bb{1}_{\mathbb{D}(0,1)}(z)$]
A simulation result for the indicator function $\bb{1}_1$ was shown in Figure \ref{fig:chi1p20}. As in previous example, we now increase the degrees: in Figure \ref{fig:chi1example4plots}, we show how the simulated zeros change as the degree $p$ goes from $20$ to $200$. As in Theorems \ref{thm:6.2} and \ref{thm:6.2equi}, the random zeros approximate to $\omega_{\mathrm{FS}}$ on the support of $\bb{1}_1(z)$ (the part $\{|z|\leq 1\}=\{r_{\mathrm{FS}}\leq 0.44311\ldots\}$). But outside the support, the simulated zeros do not converge to $c_1(L,h_L)$, and the observation is that the density of random zeros has a rapid drop outside but near the support, and has a concentration near the farthest point from the support. It is an interesting question to investigate such phenomenon in general.

\begin{figure}[H]
\centering
\includegraphics[width=0.8\textwidth]{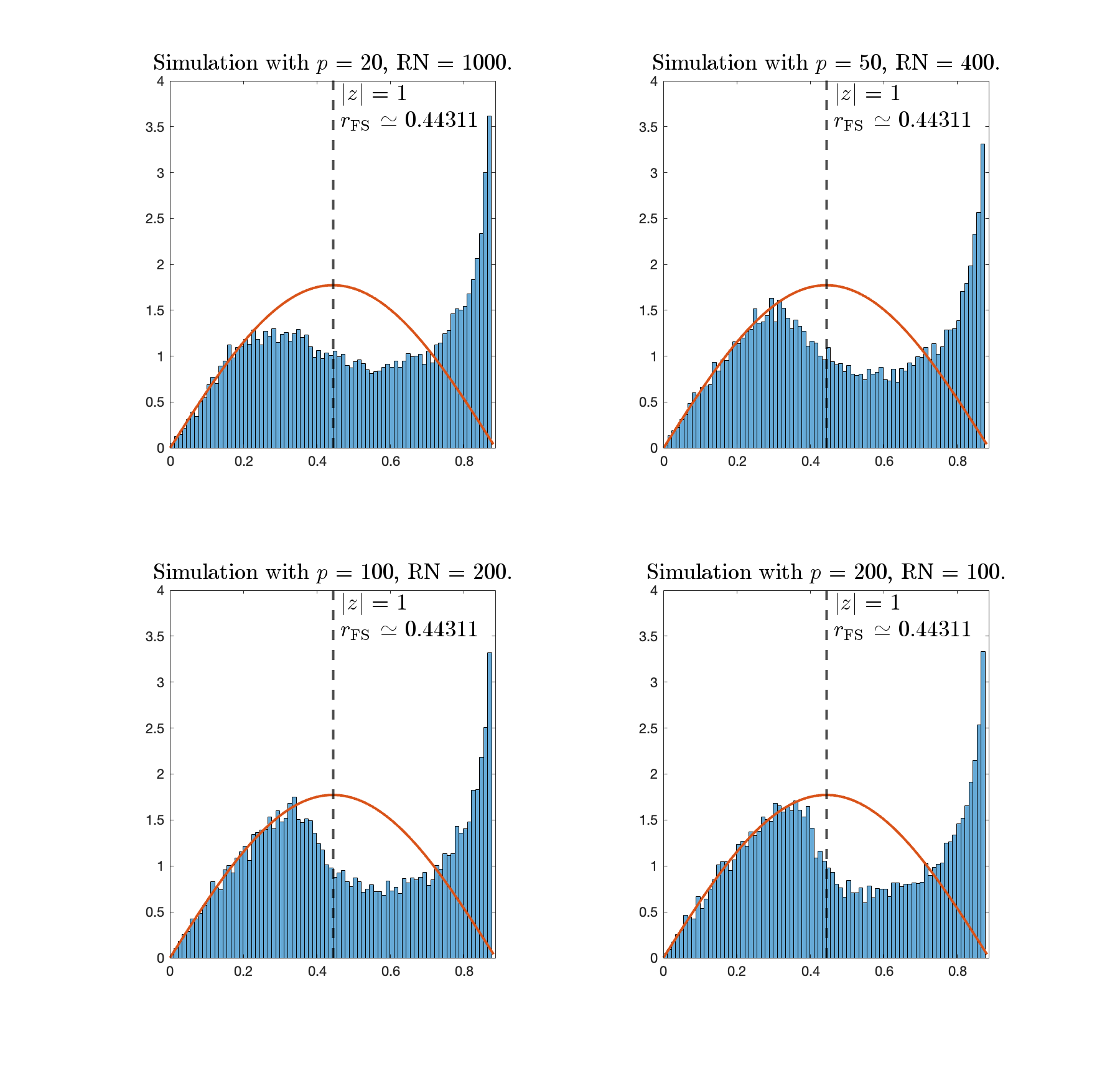}
\caption[Simulations for random zeros of $\bb{S}_{\bb{1}_1 , p}$]{Comparison of zeros of $\bb{S}_{\bb{1}_1 , p}$ ($p=20,50,100,200$) with $\omega_{\mathrm{FS}}$ on $\mathbb{CP}^1$. The density function $\psi(r_{\mathrm{FS}})$ is plotted as the red curve in each picture, and the region $\supp\, \bb{1}_1= \{|z|\leq 1\}=\{r_{\mathrm{FS}}\leq 0.44311\ldots\}$.}
\label{fig:chi1example4plots}
\end{figure}

\begin{figure}[H]
\centering
\includegraphics[width=0.8\textwidth]{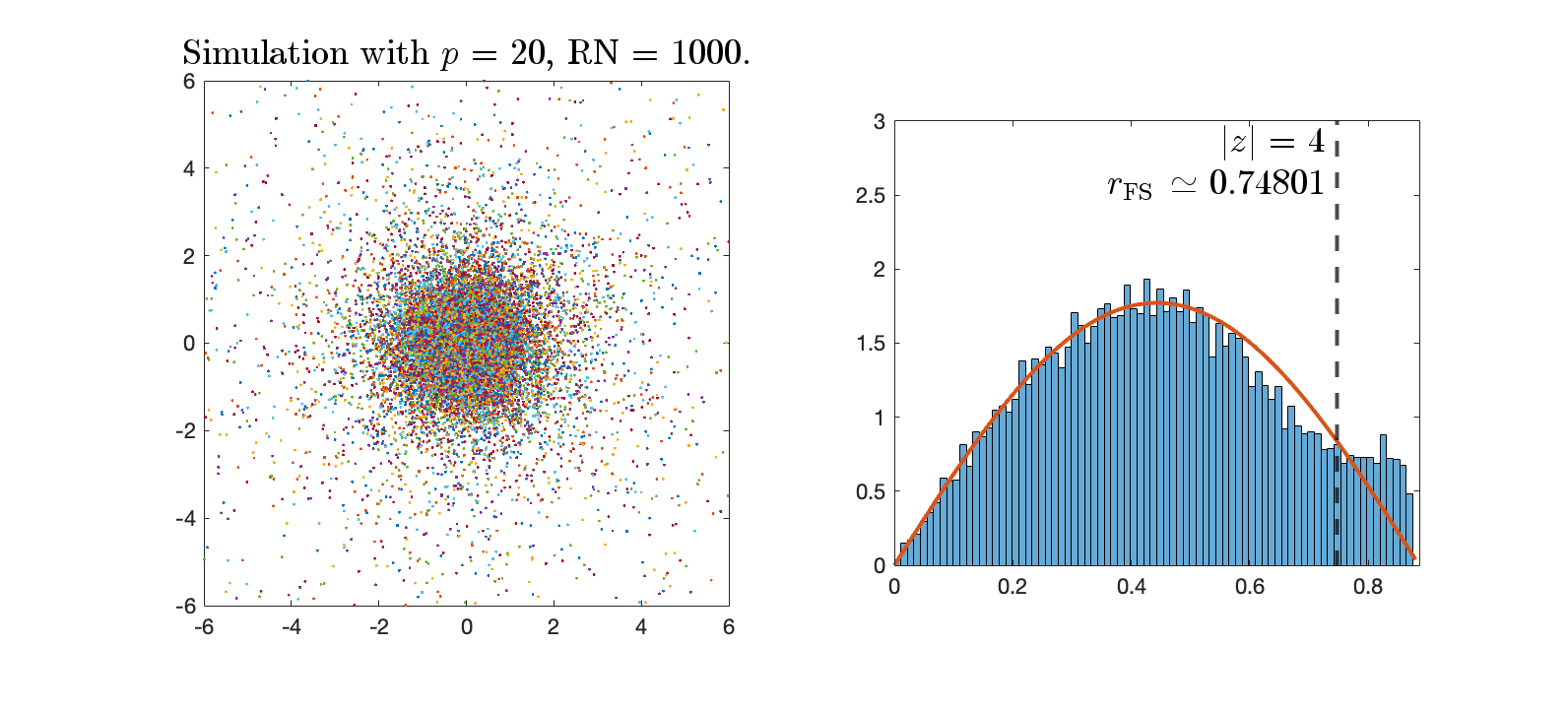}
\caption[Random zeros of $\bb{S}_{\bb{1}_4 , 20}$]{Comparison of zeros of $\bb{S}_{\bb{1}_4 , 20}$ with $\omega_{\mathrm{FS}}$ on $\mathbb{CP}^1$. The density function $\psi(r_{\mathrm{FS}})$ is plotted as the red curve in the right-hand side, and the region $\{|z|\leq 4\}=\{r_{\mathrm{FS}}\leq 0.74801\ldots\}$ corresponds to the support of $\bb{1}_4$.}
\label{fig:chi4p20}
\end{figure}

If we increase the radius $r$, that is, enlarging the support of the function $\bb{1}_r$ to define the Toeplitz operator, we will observe that the simulation of random zeros at degree $p=20$ behaves better than in Figure \ref{fig:chi1p20} when we compare it with the density function $\psi$. In Figure \ref{fig:chi4p20}, we show a simulation result for $r=4$ and $p=20$, $\mathrm{RN}=1000$, and such result is supported by Theorem \ref{thm:4.6Jan2024} under the condition that the support of $\bb{1}_r$ nearly fills in the whole part of $\mathbb{CP}^1$.

\end{Example}

\section{Smooth statistics for random zeros of {$\cLL$}-holomorphic sections}\label{section:statistics}
 In this section, we always consider the same geometric setting in Sections \ref{ss2nov} and \ref{ss4Jan}, in particular, we assume the Condition \ref{condt:main} to hold.

We also concentrate on the case $f\in \mathscr{C}^{\infty}_{\mathrm{c}}(X,\R)$ and the associated random sections $\bb{S}_{f,p}$, $p\in\N$. The goal of this section is to explain some partial extensions of the seminal results on the smooth statistics of random zeros obtained by Shiffman--Zelditch \cite{MR2465693, MR2742043} on compact K\"{a}hler manifolds to the zeros of $\bb{S}_{f,p}$ inside the support of $f$. In particular, we focus on the number variance of $\langle [\Div(\bb{S}_{f,p})],\varphi\rangle$ with a given test form $\varphi$ supported in $\supp f$ and the corresponding central limit theorem (see also \cite{STr}). Note that, except allowing $X$ to be noncompact , another difference in our geometric setting from \cite{MR2465693, MR2742043} is that we do not assume the connection between the Hermitian form $\Theta$ with $c_1(L,h_L)$. The proofs essentially follow from the arguments
presented in \cite{MR2465693, MR2742043}, 
and we will point out the necessary modifications for our setting.

\subsection{Number variance on the support of the symbol}
For $f\in \mathscr{C}^{\infty}_{\mathrm{c}}(X,\R)$
which is not identically zero, we have studied in
Section \ref{section5:onsupport} the random $(1,1)$-currents 
$[\mathrm{Div}(\bb{S}_{f,p})]$ ($p\gg 0$) on $\supp{f}$, 
especially, the expectations 
$\E\left[[\mathrm{Div}(\bb{S}_{f,p})]\right]$. 
Now we are going to study the variance of $[\mathrm{Div}(\bb{S}_{f,p})]$.

Following \cite[\S 3]{MR2465693}, we now introduce the variance current of $[\mathrm{Div}(\bb{S}_{f,p})]$. Let $\pi_1, \pi_2: X\times X\rightarrow X$ denote the projections to the first and second factors. Then if $S$ and $T$ are two currents on $X$ with respective degree $r$ and $q$, then we define a current of degree $r+q$ on $X\times X$ as follows
\begin{equation}
S\boxtimes T:= \pi_1^* S\wedge \pi_2^* T.
\end{equation} 
In particular, $[\mathrm{Div}(\bb{S}_{f,p})]\boxtimes [\mathrm{Div}(\bb{S}_{f,p})]$ defines a random $(2,2)$-current on $X\times X$.

In the same time, we introduce the following notation: for a current $T$ on $X\times X$, we write
\begin{equation}
\partial T=\partial_1 T+\partial_2 T,
\end{equation}
where $\partial_1$, $\partial_2$ denote the corresponding $\partial$-operators on the first and second factors of $X\times X$. Let $(z_1,\ldots, z_n)$ be a local complex coordinate on the first factor of $X\times X$, and let $(w_1,\ldots, w_n)$ be a local complex coordinate on the second factor of $X\times X$, then we can write locally
\begin{equation}
\partial_1=\sum dz_j \frac{\partial}{\partial z_j},\; \partial_2=\sum dw_j \frac{\partial}{\partial w_j}.
\end{equation}
Similarly, we also write $\overline{\partial} T=\overline{\partial}_1 T+\overline{\partial}_2 T$.

\begin{Definition}
The variance current of $[\mathrm{Div}(\bb{S}_{f,p})]$, denoted as $\mathbf{Var}[\bb{S}_{f,p}]$, is a $(2,2)$-current on $X\times X$ defined by
\begin{equation}
\mathbf{Var}[\bb{S}_{f,p}]:=\E\big[[\mathrm{Div}(\bb{S}_{f,p})]\boxtimes [\mathrm{Div}(\bb{S}_{f,p})]\big]-\E\left[[\mathrm{Div}(\bb{S}_{f,p})]\right]\boxtimes \E\left[[\mathrm{Div}(\bb{S}_{f,p})]\right]
\end{equation}
\end{Definition}

In order to simplify the notation, it is enough to consider only the real test forms $\Omega_{\mathrm{c}}^{n-1,n-1}(X,\R)$. For $\varphi\in \Omega_{\mathrm{c}}^{n-1,n-1}(X,\R)$, we have
\begin{equation}
\mathrm{Var}\big[\langle [\mathrm{Div}(\bb{S}_{f,p})],\varphi\rangle\big]=\langle \mathbf{Var}[\bb{S}_{f,p}], \varphi\boxtimes \varphi\rangle. 
\end{equation}

Shiffman and Zelditch \cite{MR2465693, MR2742043} established the framework to compute such variance current on a compact K\"{a}hler manifold, in particular, they obtained a pluri-bipotential for it. Their method can be easily adapted to our setting. Let us start with recalling the main ingredients from their results.

For $t\in [0,1]$, we set the function
\begin{equation}
\widetilde{G}(t):=-\frac{1}{4\pi^2} \int_0^{t^2} \frac{\log(1-s)}{s}ds=\frac{1}{4\pi^2}\sum_{j=1}^\infty \frac{t^{2j}}{j^2}.
\label{eq:6.6geb}
\end{equation}
This is an analytic function with radius of convergence $1$. Moreover, for $t\sim 0$, we have $\widetilde{G}(t)=\mathcal{O}(t^2)$.

Set $W_p=\{z\in X\;:\; T^2_{f,p}(z,z)=0\}\subset X$.
Recall that the function $N_{f,p}(z,w)$ on $X\times X$
is defined in \eqref{eq:2.60nov}. 
This is a smooth function on $X\times X\setminus (W_p\times X \cup X\times W_p)$ with values in $[0,1]$. In particular, for $z\in X\setminus W_p$, $N_{f,p}(z,z)=1$.
\begin{Definition}[{cf. \cite[Theorem 3.1]{MR2465693}}]
For $(z,w)\in X\times X\setminus (W_p\times X \cup X\times W_p)$, define
\begin{equation}
Q_{f,p}(z,w):=\widetilde{G}(N_{f,p}(z,w))=-\frac{1}{4\pi^2} \int_0^{N_{f,p}(z,w)^2} \frac{\log(1-s)}{s}ds.
\end{equation}
Then $Q_{f,p}(z,w)$ is a continuous function on $(z,w)\in X\times X\setminus (W_p\times X \cup X\times W_p)$.
\end{Definition}

Since the near-diagonal behavior of $N_{f,p}(z,w)$ depend on if there points $z,w$ lie in the support of $f$ or not, which are different from the case for Bergman kernel (such as in \cite{MR2465693, MR2742043} or \cite[Section 1.5]{DLM:21}). Following the computations in \cite[\S 3.1]{MR2465693} and we use our results proved in Theorem \ref{thm:5.1.1} and Lemma \ref{lm:2.14feb24}, we have the following results for $Q_{f,p}(z,w)$.
\begin{Proposition}[{cf.\ \cite[\S 3.1]{MR2465693}}]\label{prop:6.3feb}
Let $U$ be an open subset of $X$ such that $\overline{U}\subset \{f\neq 0\}$ (hence $\overline{U}$ is compact).
\begin{itemize}
\item[(i)] Then there exists an integer $p_0\in \N$ such that for all $p\geq p_0$, $T^2_{f,p}(z,z)$ never vanishes on $\overline{U}$. Moreover, for all $p\geq p_0$, the function $Q_{f,p}(z,w)$ is smooth on the region $U\times U\setminus \Delta_U$ ($\Delta_U$ denotes the diagonal) and it is $\mathscr{C}^1$ on $U\times U$.
\item[(ii)] Fix $b\gg 0$ and $\epsilon>0$, then for all sufficiently large $p$ and for $x\in U$, $v\in T_x X$ with $\|v\|\leq b\sqrt{\log{p}}$, we have
\begin{equation}
Q_{f,p}(x,\exp_x(v'/\sqrt{p}))=\widetilde{G}\left(\exp(-\Phi_{x}(0,v')^{2}/4)\right)+\mathcal{O}(p^{-1/2+\epsilon}),
\end{equation}
where $\Phi_{x}(0,v')$ is defined in \eqref{eq:2.16July2}.
\item[(iii)] For given $k, \ell\in \N$, there exist a sufficiently large $b>0$ such that there exist a constant $C>0$ such that for all $z,w\in U$, $\mathrm{dist}(z,w)\geq b\sqrt{\log{p}/p}$, we have
\begin{equation}
|\nabla^\ell_{z,w}Q_{f,p}(z,w)|\leq C p^{-k}.
\label{eq:6feb24}
\end{equation}
\end{itemize}
\end{Proposition}

For a real $(n-1,n-1)$-form  $\varphi$ on $X$ with $\mathscr{C}^3$-coefficients, recall that $L(\varphi)\in \mathscr{C}^1(X,\R)$ is defined by
\begin{equation}
 \sqrt{-1}\partial\overline{\partial} \varphi =
 L(\varphi) \frac{c_1(L,h_L)^n}{n!}\,\cdot
\end{equation}
Recall also that we have two volume forms 
$\mathrm{dV}=\Theta^n/n!$ and $\mathrm{dV}^L:=c_1(L,h_L)^n/n!$ 
(see \eqref{eq:volumeL}). Moreover, we have
\begin{equation}
\mathrm{dV}^L(z)=\bb{b}_0(z) \mathrm{dV}(z),
\end{equation}
where the positive function 
$\bb{b}_0(z)=\det(\dot{R}_z^{L}/2\pi)$ on $X$ is given in \eqref{eq:b0}.
\begin{Theorem}\label{thm:6.4number}
Let $(X,J,\Theta)$ be a connected complex Hermitian manifold and let $(L,h_{L})$, $(E,h_{E})$ be two holomorphic line bundles on $X$ with smooth Hermitian metrics.
We assume Condition \ref{condt:main} (see also \eqref{eq:mainassumptions}) to hold. Fix $f\in \mathscr{C}^{\infty}_{\mathrm{c}}(X,\R)$ which is not identically zero, and let $U$ be an open subset of $X$ such that $\overline{U}\subset \{f\neq 0\}$. Then for sufficiently large $p$, we have the identity of $(2,2)$-currents on $U\times U$,
\begin{equation}\label{eq:6.11feb24}
\mathbf{Var}[\bb{S}_{f,p}]|_{U\times U}=-\partial_1\overline{\partial}_1\partial_2\overline{\partial}_2 Q_{f,p}(z,w)|_{U\times U}=(\sqrt{-1}\partial\overline{\partial})_z  (\sqrt{-1}\partial\overline{\partial})_w Q_{f,p}(z,w)|_{U\times U}.
\end{equation}

Let $\varphi$ be a real $(n-1,n-1)$-form on $X$ with $\mathscr{C}^3$-coefficients and $\supp\varphi\subset U$, then we have the formula for $p\gg 0$,
\begin{equation}\label{eq:6.12feb24}
\mathrm{Var}\big[\langle [\mathrm{Div}(\bb{S}_{f,p})],\varphi\rangle\big]=p^{-n}\left(\frac{\zeta(n+2)}{4\pi^2}\int_U |L(\varphi)(z)|^2 \mathrm{dV}^L(z)+\mathcal{O}(p^{-1/2+\epsilon})\right),
\end{equation}
where
$$\zeta(n+2)=\sum_{k=1}^\infty \frac{1}{k^{n+2}}.$$ 
\end{Theorem}
\begin{proof}
We sketch the proof based at the proofs of \cite[Theorem 3.1]{MR2465693} and \cite[\S 3.1]{MR2742043}. At first, fix a test form $\varphi\in\Omega^{n-1,n-1}_{\mathrm{c}}(U,\R)$, a routine calculation (see also \cite[Proof of Theorem 3.7]{DrLM:2023aa}) shows that
\begin{equation}\label{pb3.25}
\begin{split} 
&\mathrm{Var}\big[\left\langle[\mathrm{Div}(\bb{S}_{f,p})], \varphi\right\rangle\big]\\
&= -\frac{1}{\pi^2 } \int_{U\times U} (\partial\overline{\partial}\varphi(z))\wedge
(\partial\overline{\partial}\varphi(w)) \\
&\qquad\qquad\times\E\left[ \log|
T^2_{f,p}(z,z)^{-1/2}\bb{S}_{f,p}(z)|_{h_{p}}\cdot\log|
T^2_{f,p}(w,w)^{-1/2}\bb{S}_{f,p}(w)|_{h_{p}}\right].
\end{split}
\end{equation}
Then by Proposition \ref{prop:6.3feb} (i) and \cite[Lemma 3.3]{MR2465693}, on $U\times U$ and for all $p\gg 0$, we have 
\begin{equation}
\E\left[ \log|
T^2_{f,p}(z,z)^{-1/2}\bb{S}_{f,p}(z)|_{h_{p}}\cdot\log|
T^2_{f,p}(w,w)^{-1/2}\bb{S}_{f,p}(w)|_{h_{p}}\right]=\frac{\gamma^2}{4}+\pi^2\widetilde{G}(N_{f,p}(z,w)),
\end{equation}
where $\gamma$ is the Euler's constant. Then we can rewrite \eqref{pb3.25} as
\begin{equation}\label{eq:6.15feb}
\begin{split} 
\mathrm{Var}\big[\left\langle[\mathrm{Div}(\bb{S}_{f,p})], \varphi\right\rangle\big] &= - \int_{U\times U} (\partial\overline{\partial}\varphi(z))\wedge
(\partial\overline{\partial}\varphi(w))\widetilde{G}(N_{f,p}(z,w))\\
&=\langle (\sqrt{-1}\partial\overline{\partial})_z  (\sqrt{-1}\partial\overline{\partial})_w Q_{f,p}(z,w),\varphi\boxtimes\varphi\rangle.
\end{split}
\end{equation}
This way, we get \eqref{eq:6.11feb24}. In fact, \eqref{eq:6.15feb} still holds if $\varphi$ is with $\mathscr{C}^3$-coefficients, now we show \eqref{eq:6.12feb24} by using Proposition \ref{prop:6.3feb}-(ii) and (iii).

The first step is to rewrite the integral in the form
\begin{equation}
\int_{U\times U}\cdots=\int_{z\in U}\int_{\{z\}\times U}\cdots.
\end{equation}
As in \cite[\S 3.1]{MR2742043}, we set for $z\in U$,
\begin{equation}
\begin{split}
\mathcal{I}_p(z)&=\int_{\{z\}\times U}Q_{f,p}(z,w)
(\sqrt{-1}\partial\overline{\partial}\varphi(w))\\
&=\int_{\{z\}\times U} Q_{f,p}(z,w)L(\varphi)(w)
\bb{b}_0(w)\mathrm{dV}(w).
\end{split}
\label{eq:6.18feb24}
\end{equation}
Let $b>0$ be a fixed number which is sufficiently large as in
Proposition \ref{prop:6.3feb} (iii) with $k=n+1$. 
Then we have for $p\gg 0$,
\begin{equation}
\begin{split}
\mathcal{I}_p(z)=\int_{v\in T_z X, \|v\|\leq 
b\sqrt{\log{p}}} Q_{f,p}(z,\exp_z(v/\sqrt{p}))
(L(\varphi)\bb{b}_0)(\exp_z(v/\sqrt{p}))
\mathrm{dV}(\exp_z(v/\sqrt{p}))&\\
+\mathcal{O}(p^{-n-1}).&
\end{split}
\label{eq:6.19feb24}
\end{equation}

As in \eqref{eq:2.18July}, let 
$\mathrm{dV_{Eucl}}_{,z}(v)$ denote the Euclidean 
volume form on the real vector space $(T_z X, g^{TX}_z)$. 
Then for $v\in T_z X, \|v\|\leq b\sqrt{\log{p}}$, we have
\begin{equation}
\mathrm{dV}(\exp_z(v/\sqrt{p}))=
\frac{1}{p^n} \left(\mathrm{dV_{Eucl}}_{,z}(v)+
\mathcal{O}(\sqrt{p^{-1}\log{p}})\right).
\label{eq:6.20feb24}
\end{equation}
Since $L(\varphi)$ is $\mathscr{C}^1$ on $\overline{U}$, 
then by Proposition \ref{prop:6.3feb}-(ii) and 
\cite[(34)--(37)]{MR2742043}, 
we get for a fix small $\epsilon>0$,
\begin{equation}
\begin{split}
\mathcal{I}_p(z)=\frac{(L(\varphi)\bb{b}_0)(z)}{p^n}
\int_{v\in T_z X, \|v\|\leq b\sqrt{\log{p}}}
\widetilde{G}\left(\exp(-\Phi_{x}(0,v')^{2}/4)\right)
\mathrm{dV_{Eucl}}_{,z}(v)+\mathcal{O}(p^{-n-1/2+\epsilon}).
\end{split}
\label{eq:6.21feb24}
\end{equation}
By \eqref{eq:5.1.15paris}, we can fix $b$ to be large enough such that 
\begin{equation}
\int_{v\in T_z X, \|v\|\geq b\sqrt{\log{p}}}
\widetilde{G}\left(\exp(-\Phi_{x}(0,v')^{2}/4)\right)
\mathrm{dV_{Eucl}}_{,z}(v) =\mathcal{O}(p^{-n-1}).
\label{eq:6.22feb24}
\end{equation}
In the same time, by the formula \eqref{eq:2.16July2}
for $\Phi_{x}(0,v')$, then for $k\in\N_{\geq 1}$, we have
\begin{equation}
\int_{v\in T_z X}\left(\exp(-k\Phi_{x}(0,v')^{2}/2)\right)
\mathrm{dV_{Eucl}}_{,z}(v) = \frac{1}{\bb{b}_0(z) k^{n}}.
\label{eq:6.23feb24}
\end{equation}
Finally, combining the Taylor series \eqref{eq:6.6geb}
with \eqref{eq:6.18feb24} -- \eqref{eq:6.23feb24}, we get for $z\in U$,
\begin{equation}
\begin{split}
\mathcal{I}_p(z)=p^{-n}\left(\frac{\zeta(n+2)}{4\pi^2}
L(\varphi)(z)+\mathcal{O}(p^{-1/2+\epsilon})\right),
\end{split}
\label{eq:6.24feb24}
\end{equation}
after taking the integration with respect to $z\in U$, we conclude exactly \eqref{eq:6.12feb24}. This way, we complete our proof.
\end{proof}

\subsection{Asymptotic normality of random zeros;
proof of Theorem \ref{thm:3.5.1ss}}
\label{ss:6.2CLT}
The asymptotic normality of the zeros of random holomorphic functions 
or sections has been introduced and proved by Sodin--Tsirelson 
\cite{STr} for certain random holomorphic functions on 
$\C$ or $\mathbb{D}$ and by Shiffman--Zelditch 
\cite[Theorem 1.2 and \S 4]{MR2742043} for the random holomorphic sections of 
line bundles on a compact K\"{a}hler manifold. One key ingredient in 
their approaches is the normalized Bergman kernel which is
the covariance function of the corresponding Gaussian holomorphic fields on $\C$ or $X$, analogous to the construction in the proof of Proposition \ref{prop:3.17}. Then the 
problem is reduced to the seminal result proved by Sodin and Tsirelson 
in \cite[Theorem 2.2]{STr} for the non-linear functionals of Gaussian process.

Let us recall the main result of \cite[\S 2.1]{STr}. Let $(T,\mu)$ be a measure space with a finite positive measure $\mu$ (with $\mu(T)>0$). We also fix a sequence of measurable functions $A_k: T\rightarrow \C$, $k\in\N$ such that on $T$,
\begin{equation}
\sum_k |A_k(t)|^2\equiv 1.
\label{eq:6.25feb24}
\end{equation}
We consider a complex-valued Gaussian process on $T$ defined as
\begin{equation}
W(t):=\sum_k \eta_k A_k(t).
\label{eq:6.26feb24}
\end{equation}
Then $\{\eta_k\}$ is a sequence of i.i.d.\ standard complex Gaussian random variables. Then for each $t\in T$, $W(t)\sim \mathcal{N}_{\C}(0,1)$. The covariance function for $W$ is $\rho_W: T\times T\rightarrow \C$ given by
\begin{equation}
\rho_W(s,t):=\E[W(s)\overline{W(t)}]=\sum_k A_k(s)\overline{A_k(t)}.
\end{equation}

Let $\{W_p\}_{p\in \N}$ be a sequence of independent Gaussian processes on $T$ described as above, and let $\rho_p(s,t)$ ($p\in\N$) denote the corresponding covariance functions. We also fix a nontrivial real function $F\in\mathcal{L}^2(\R_+, e^{-r^2/2}r dr)$, and a bounded measurable function $\psi: T\rightarrow \R$, set
\begin{equation}
Z_p:=\int_T F(|W_p(t)|)\psi(t)d\mu(t).
\label{eq:6.28feb24}
\end{equation}

Sodin and Tsirelson proved the following results.

\begin{Theorem}[{\cite[Theorem 2.2]{STr}}]\label{thm:STr2004}
With the above constructions, and suppose that
\begin{itemize}
	\item[(i)] $$\liminf_{p\rightarrow +\infty} \frac{\int_{T}\int_{T} |\rho_p(s,t)|^{2\alpha}\psi(s)\psi(t)d\mu(s)d\mu(t)}{\sup_{s\in T}\int_T |\rho_p(s,t)|d\mu(t)}>0,$$
	for $\alpha=1$ if $f$ is monotonically increasing, or for all $\alpha\in \N$ otherwise;
	\item[(ii)] $$\lim_{p\rightarrow+\infty} \sup_{s\in T} \int_T |\rho_p(s,t)|d\mu(t)=0.$$
\end{itemize}
Then the distributions of the random variables
\begin{equation}
\frac{Z_p-\E[Z_p]}{\sqrt{\mathrm{Var}[Z_p]}}
\end{equation}
converge weakly to  the (real) standard Gaussian distribution $\mathcal{N}_\R(0,1)$ as $p\rightarrow +\infty$.
\end{Theorem}

Now we are ready to present the proof of Theorem \ref{thm:3.5.1ss}.
\begin{proof}[Proof of Theorem \ref{thm:3.5.1ss}]
The proof is an easy modification of \cite[\S 4 Proof of Theorem 1.2]{MR2742043}, together with the results in Theorem \ref{thm:5.1.1} due to the assumption $\overline{U}\subset\{f\neq 0\}$.

We take $F(r)=\log{r}$, $(T,\mu)=(U,\mathrm{dV}^L|_U)$, 
$\psi(z)=\frac{1}{\pi} L(\varphi)(z)$ which satisfies the 
conditions in Theorem \ref{thm:STr2004}. 
Let $\sigma:\overline{U}\rightarrow L$ be a continuous section 
such that $|\sigma(z)|_{h_L}\equiv 1$ on $\overline{U}$. 
For each $p$, fix an orthonormal basis 
$\{S^p_j\}_{j=1}^{d_p(f)}$ consisting of the eigensections of $T_{f,p}$
for nonzero eigenvalues $\{\lambda^p_j\}_j$. 
Then on $\overline{U}$, we write
\begin{equation}
S^p_j(z)=a^p_j(z)\sigma^{\otimes p}(z).
\end{equation}
Then we can set $A^p_j(z)={\lambda^p_j a^p_j(z)}/{\sqrt{T^2_{f,p}(z,z)}}$,
which forms a sequence of measurable functions on $U$ satisfying \eqref{eq:6.25feb24}. Then by \eqref{eq:4.34series}, we have the identity on $U$
\begin{equation}
\frac{\bb{S}_{f,p}(z)}{\sqrt{T^2_{f,p}(z,z)}}=W_p(z) \sigma^{\otimes p}(z),
\label{eq:6.33feb24}
\end{equation}
where $W_p$ is the Gaussian process on $U$ constructed
as in \eqref{eq:6.26feb24}. The covariance function 
$\rho_p(z,w)$ for $W_p$ is given by
\begin{equation}
|\rho_p(z,w)|=N_{f,p}(z,w).
\end{equation}
Let $Z_p$ be the random variable defined
as in \eqref{eq:6.28feb24} from $W_p$. 
Then \eqref{eq:6.33feb24} implies that
\begin{equation}
Z_{f,p}(\varphi)=Z_p + C_p,
\end{equation}
where $C_p$ is a deterministic constant. Thus the asymptotic 
normality of $Z_{f,p}(\varphi)$ is equivalent to that of 
$Z_p$, which follows by checking the Conditions (i) and (ii) 
in Theorem \ref{thm:STr2004} for $N_{f,p}(z,w)$. 
Finally, we apply Theorem \ref{thm:5.1.1} and proceed
as in the last part of \cite[\S 4]{MR2742043}. This completes the proof.
\end{proof}

\begin{center}
	{\large Declarations}
\end{center}

\noindent
{\bf Conflict of interest.} On behalf of all authors, the corresponding author states that there is
no conflict of interest.

\noindent
{\bf Data sharing.} Data sharing not applicable to this article as no datasets were generated or
analyzed during the current study.

\bibliographystyle{siam}

\bibliography{DLM24}

\end{document}